\documentclass[a4paper, 12pt]{extarticle}
\usepackage{mathtext}
\usepackage{graphicx}
\usepackage{amscd}
\usepackage[T2A]{fontenc}
\usepackage[english]{babel}
\usepackage{latexsym,amsfonts,amssymb,amsmath,longtable,amsthm}
\usepackage[pic]{xy}
\usepackage{bbm,makeidx}
\usepackage[centerlast,small]{caption}
\usepackage{amsmath}
\usepackage[cp1251]{inputenc}
\newtheorem{lem}{{\scshape Lemma}}
\newtheorem{cor}{{\scshape Corollary}}
\newtheorem{rmk}{{\scshape Remark}}

\newtheorem{ttt}{{\scshape Theorem}}
\newtheorem{prp}{{\scshape Proposition}}

\begin{document}
\title{The classification of fused links}
\author{Timur R. Nasybullov\footnote{The author is supported by  Russian Science Foundation (project 14-21-00065)}~\footnote{e-mail:~timur.nasybullov@mail.ru, ntr@math.nsc.ru}}
\maketitle
\begin{abstract}
We construct the complete invariant for fused links. It is proved that the set of equivalence classes of $n$-component fused links is in one-to-one correspondence with the set of elements of the abelization $UVP_n/UVP_n^{\prime}$ up to conjugation by the elements from the symmetric group $S_n<UVB_n$.

~
\end{abstract}
\section{Introduction}
Knot invariants are functions of knots that do not change under isotopies.
The study of knot invariants is at the core of knot theory. Indeed, the isotopy
class of a knot is, tautologically, a knot invariant. During last years different authors constructed vast number of knot invariants: the (self) linking number, the unknotting number, the knot group, the knot quandle, the Jones polynomial, the Conway polynomial and so on (see, for example, \cite{Man, DCM, Joy}). The disadvantage of a large number of easy countable invariants is that they are not complete, i.~e. they do not distinguish some knots or links.

 At the same time there are few examples of complete knot invariants, but usually it is difficult to understand if the value of the invariant on the knot coincides with the value of this invariant on the another knot. The complement of a knot itself (as a topological space) is known to be a complete invariant of the knot by the theorem of C.~Gordon and J.~Luecke \cite{GL} in the sense that it distinguishes the given knot from all other knots up to ambient isotopy and mirror image. Some invariants associated with the knot complement include the knot group which is just the fundamental group of the complement. The knot quandle is also a complete invariant  in this sense \cite{Joy} but it is difficult to determine if two quandles are isomorphic.

 Thus an important problem in knot theory is to construct a complete knot invariant which can be easily found and used.

 Recently some generalizations of classical knots and links were defined and studied: singular links \cite{Vas, Bir}, virtual links \cite{GuPoVi, Kauff}, welded
links \cite{FRiRu} and fused links \cite{Kad, BBD, Kalu}. The problem of constructing invariants is also important for all of this knot theories.

One of the ways of studying classical links is to study the braid group.
Singular braids \cite{Bir, Ba}, virtual braids \cite{Kauff,Ver}, welded braids \cite{FRiRu} and unrestricted virtual braids \cite{Kauf,Kad} were defined
similar to the classical braid group adding the extra generators and relations. Theorem of A.~A.~Markov \cite[\S2.2]{Bir2} reduces the problem of classification of links to some algebraic problems of the theory of braid groups. There are generalizations
of Markov theorem for virtual links, welded links and fused links \cite{SK}.

\begin{figure}[bh]
\noindent\centering{
\includegraphics[width=50mm]{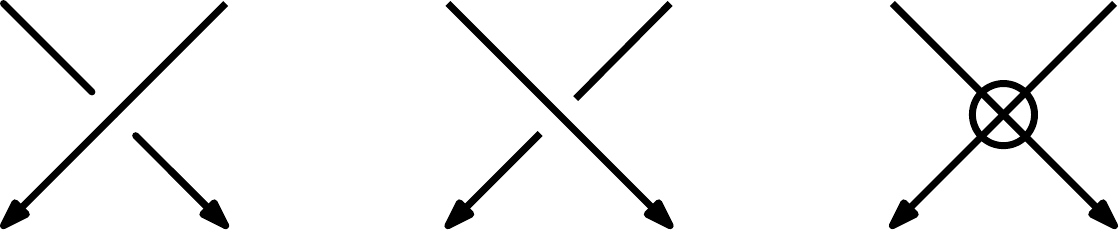}
}
\caption{Crossings in the double welded knot diagram}
\label{gend1}
\end{figure}

\begin{figure}[bh]
\noindent\centering{
\includegraphics[width=120mm]{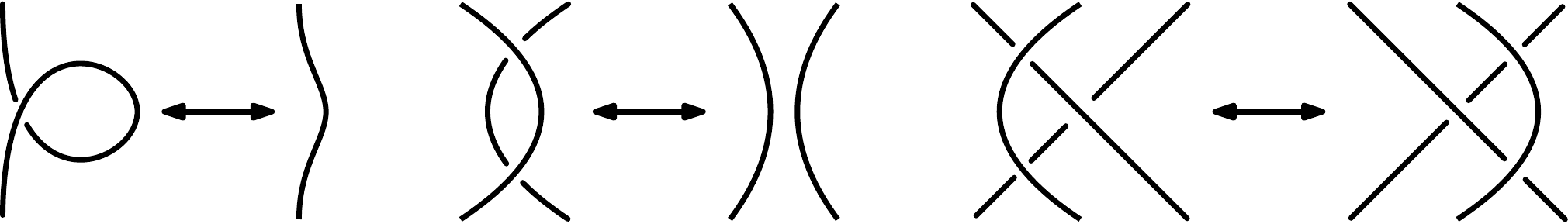}
}
\caption{Classical Reidemmeister moves}
\label{gend2}
\end{figure}

\begin{figure}[bh]
\noindent\centering{
\includegraphics[width=120mm]{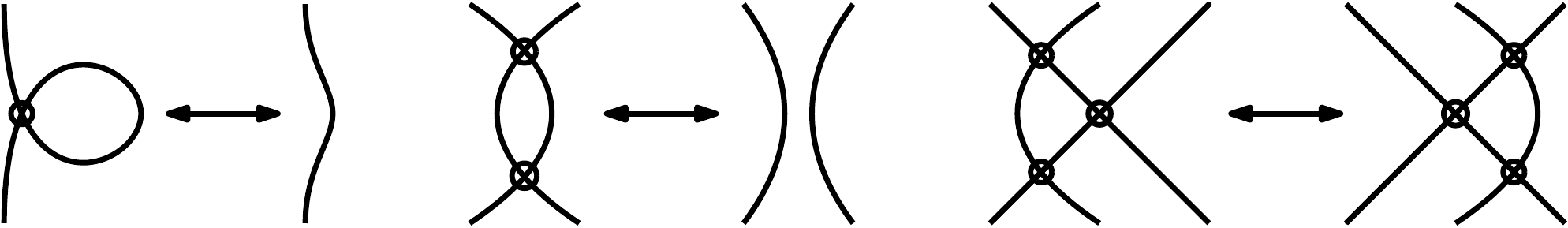}
}
\caption{Virtual Reidemeister moves}
\label{gend3}
\end{figure}

\begin{figure}
\noindent\centering{
\includegraphics[width=50mm]{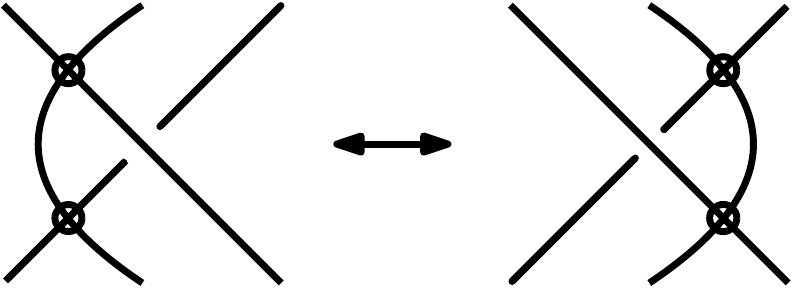}
}
\caption{Mixed Reidemeister moves}
\label{gend4}
\end{figure}

\begin{figure}
\noindent\centering{
\includegraphics[width=100mm]{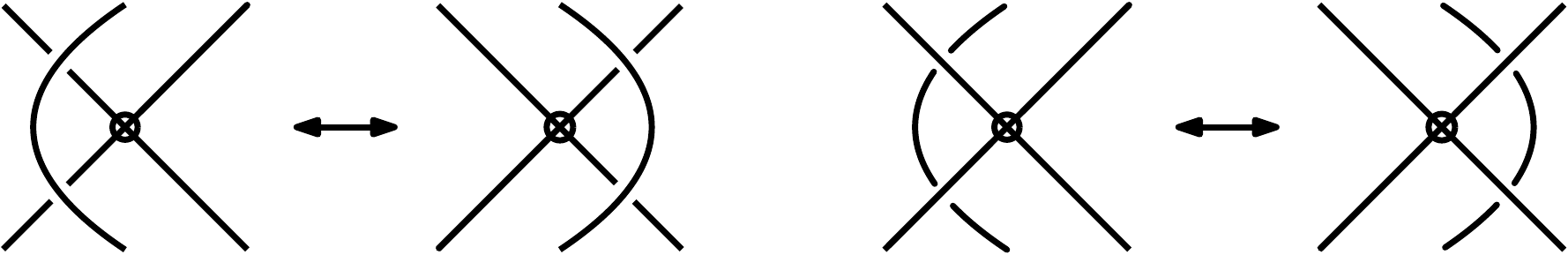}}
\caption{Forbidden moves}
\label{gend5}
\end{figure}

 In the paper we study fused links, which were defined by L.~Kauffman and S.~Lambropoulou in \cite{Kalu}, and their invariants. Fused links are represented as generic immersions of circles
in the plane (fused link diagrams) where double points can be classical (with the usual
information on overpasses and underpasses) or virtual (see Fig. \ref{gend1}).
\\
Fused link diagrams are equivalent
under ambient isotopy and some types of local moves (generalized Reidemeister moves):
classical Reidemeister moves (see Fig. \ref{gend2}), virtual Reidemeister moves (see Fig. \ref{gend3}), mixed Reidemeister
moves (see Fig. \ref{gend4}) and Forbidden moves (see. Fig. \ref{gend5}).

\begin{figure}[bh]
\noindent\centering{
\includegraphics[width=120mm]{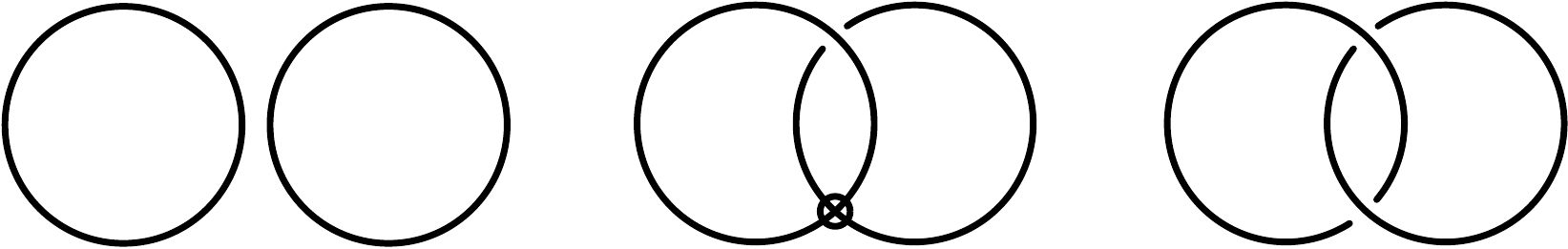}
}
\caption{Different 2-component fused links}
\label{hopf}
\end{figure}

 In the theory of fused links every knot is equivalent to the unknot \cite{K}. However not every link is equivalent to the trivial link. For example, trivial 2-component link, Hopf link and Hopf link with one virtual crossing and with one classical crossing all are different (see Fig. \ref{hopf}).
 \\
 The full classification of fused links is not
(completely) trivial. In particular, A.~Fish and E.~Keyman proved that the fused link with classical crossings only is completely determined by the linking
numbers of each pair of components \cite{FiKe}. It means that the set of linking
numbers for each pair of components of fused links is a full invariant for fused links which have only classical crossings.

In the present paper we find all non-equivalent classes of fused links and construct an easy computable full invariant for fused links. We use the following proposition, which is implicitly formulated in \cite{BBD}.\\
{\scshape\textbf{Proposition.}} There exists a map $\varrho^*: UVB_{\infty}\to UVP_{\infty}$, such that the closures of the braids $\beta$ and $\varrho^*(\beta)$ are equivalent as fused links.

Denote by $T$ the set of coset representatives of $UVP_n/UVP_n^{\prime}$ and for the element $\alpha\in UVP_n$ denote by $\overline{\alpha}\in T$ the unique coset representative $\alpha UVP_n^{\prime}=\overline{\alpha}UVP_n^{\prime}.$
Then the main result of the paper can be formulated in the following form.\\
{\scshape\textbf{Theorem.}} Let $\alpha$ and $\beta$ be unrestricted virtual braids. Then their closures $\widehat{\alpha}$ and $\widehat{\beta}$ are equivalent as fused links if and only if $\overline{\varrho^*(\alpha)}$ and $\overline{\varrho^*(\beta)}$ are conjugated by the element from $S_n<UVB_n$.

Thus the map $\widehat{\alpha}\to \overline{\varrho^*(\alpha)}^{S_n}$ is a complete invariant for fused links and in order to understand that two fused links $\widehat{\alpha}$ an $\widehat{\beta}$ are equivalent or not  we just need to compare two finite sets  $\overline{\varrho^*(\alpha)}^{S_n}$ and $\overline{\varrho^*(\beta)}^{S_n}$, or equivalently, to understand that the element $\overline{\varrho^*(\alpha)}$ belongs to the finite set $\overline{\varrho^*(\beta)}^{S_n}$ or not.

At the almost same time  B.~Audoux, P.~Bellingeri, J.-B.~Meilhan and E.~Wagner found the full classification of fused links independently from the author \cite[Theorem 2.5]{suki}. Using the different from the author's methods they proved that every fused link is completely determined by the set of virtual linking numbers for each pair of components. The result is the same (however in different formulation) and these two approaches seem to us complementary and both interesting.

The author is grateful to Valeriy Bardakov for multiple helpful advices during the work on the paper. Also the author would like to thank the University of Bologna (Italy) where a part of this work was completed and especially Prof. Michele Mulazzani for his help and support.

\section{Definitions and results}

\begin{figure}[bh]
\noindent\centering{
\includegraphics[width=120mm]{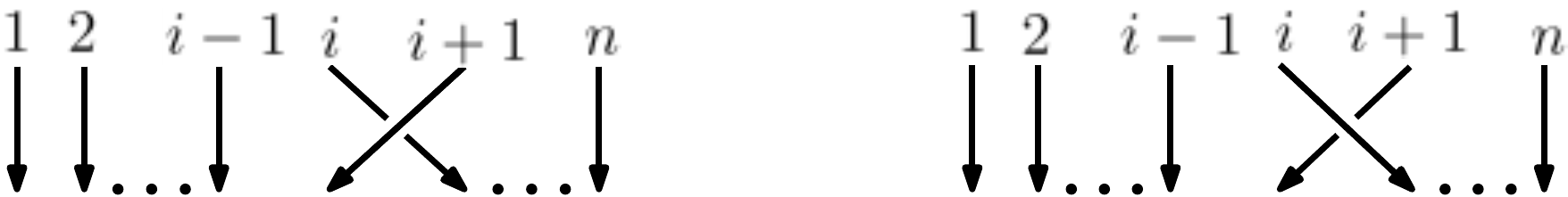}
}
\caption{Geometric braids representing $\sigma_i$ (on the left) and $\sigma_i^{-1}$ (on the right)}
\label{gen1}
\end{figure}

In this section we fix notation and recall basic definitions and known results about different generalizations of braid groups.
The classical braid group $B_n$ on $n$ strands ($n>1$) is the group generated by the elements $\sigma_1,\dots,\sigma_{n-1}$ (see Fig. \ref{gen1})
\\
with the following defining relations.
\begin{align}
\sigma_i\sigma_{i+1}\sigma_i&=\sigma_{i+1}\sigma_i\sigma_{i+1}&i=1,\dots,n-2;\tag{$B_1$}\label{B1}\\
\sigma_i\sigma_j&=\sigma_{j}\sigma_i&|i-j|\geq2.\tag{$B_2$}\label{B2}
\end{align}
There exists a homomorphism $\iota: B_n\to S_n$ from the braid group $B_n$ onto the symmetric group $S_n$ on $n$ letters. This homomorphism maps the generator $\sigma_i$ to the transposition $\tau_i=(i,~i+1)$ for $i=1, 2, \dots, n - 1$. The kernel of this homomorphism is called pure
braid group on $n$ strands and denoted by $P_n$.

\begin{figure}[bh]
\noindent\centering{
\includegraphics[width=45mm]{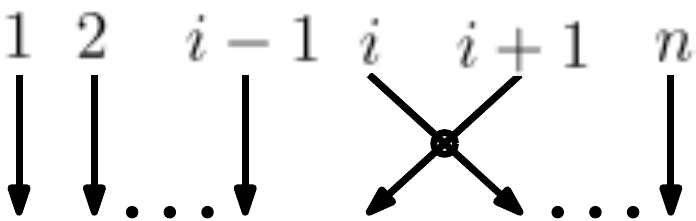}
}
\caption{Geometric braid representing $\rho_i$}
\label{gen2}
\end{figure}

The virtual braid group $VB_n$ is a group obtained from $B_n$ adding new generators $\rho_1,\dots,\rho_{n-1}$ (see Fig. \ref{gen2})
\\
and additional relations
\begin{align}
\rho_i\rho_{i+1}\rho_i&=\rho_{i+1}\rho_i\rho_{i+1}&i=1,\dots,n-2;\tag{$P_1$}\label{P1}\\
\rho_i\rho_j&=\rho_{j}\rho_i&|i-j|\geq2;\tag{$P_2$}\label{P2}\\
\rho_i^2&=1&i=1,\dots,n-1;\tag{$P_3$}\label{P3}
\end{align}
\begin{align}
\sigma_i\rho_j&=\rho_{j}\sigma_i&|i-j|\geq2;\tag{$M_1$}\label{M1}\\
\rho_i\rho_{i+1}\sigma_i&=\sigma_{i+1}\rho_i\rho_{i+1}&i=1,\dots,n-2.\tag{$M_2$}\label{M2}
\end{align}
It is easy to verify that the elements $\rho_1,\dots,\rho_{n-1}$ generate the symmetric group $S_n$. Also it is known that the elements $\sigma_1,\dots,\sigma_{n-1}$ generate
the braid group $B_n$.
In the paper \cite{GuPoVi} it is proved that the relations
\begin{align}
\rho_i\sigma_{i+1}\sigma_i&=\sigma_{i+1}\sigma_i\rho_{i+1}&i=1,\dots,n-2;\tag{$F_1$}\label{F1}\\
\rho_{i+1}\sigma_i\sigma_{i+1}&=\sigma_i\sigma_{i+1}\rho_i&i=1,\dots,n-2;\tag{$F_2$}\label{F2}
\end{align}
do not hold in the group $VB_n$. According to \cite{FRiRu} the welded braid group $WB_n$ on $n$ strands is a quotient of the group $VB_n$ by the forbidden relation (\ref{F1}), i.~e. it is a group with the generators $\sigma_1,\dots,\sigma_{n-1},\rho_{1},\dots,\rho_{n-1}$ and relations (\ref{B1})--(\ref{F1}). If we add to the group $WB_n$ the second forbidden relation (\ref{F2}) the we get the unrestricted virtual braid group $UVB_n$.

The elements $\sigma_1,\dots,\sigma_{n-2},\rho_1,\dots,\rho_{n-2}$ generate the subgroup $UVB_{n-1}$ in the group $UVB_n$. Then we have the following chain of inclusions.
$$UVB_{2}<UVB_{3}<UVB_{4}<\dots<UVB_{\infty}=\bigcup_{n\geq2}UVB_{n}$$
The homomorphism $\iota$ can be extended to the homomorphism $UVB_n\to S_n$ by the rule $\iota:\sigma_i, \rho_i\mapsto\tau_i=(i~i+1)$. The kernel of this homomorphism is called pure unrestricted virtual braid group  and is denoted by $UVP_n$. The group $UVP_{\infty}$ is defined analogically to the group $UVB_{\infty}$.
$$UVP_{2}<UVP_{3}<UVP_{4}<\dots<UVP_{\infty}=\bigcup_{n\geq2}UVP_{n}$$

The symmetric group $S_n$ acts on the set $\{1,\dots,n\}$. By the symbol $\pi(k)$ we denote an image of the integer $k\in\{1,\dots,n\}$ under the permutation $\pi$. If for the braid $\alpha$ we have $\iota(\alpha)(k)=k$, then we say that $\alpha$ fixes $k$ or $k$ is fixed by $\alpha$.

We say that the braid $\alpha$ does not involve the strand $j$ if $\alpha$ belongs to the group $\langle\sigma_1,\dots,\sigma_{j-2},\sigma_{j+1},\dots, \sigma_{n-1},\rho_1,\dots,\rho_{j-2},\rho_{j+1},\dots, \rho_{n-1}\rangle$. It is obvious that the braid $\alpha\in\langle\rho_1,\dots,\rho_{n-1}\rangle$ which fixes the strands $n-1,n$ does not involve an $n$-strand.

We say that the braid $\beta\in UVB_{n+1}$ is obtained from the braid $\alpha\in UVB_n$ by right stabilization of positive (negative, virtual) type if $\beta=\alpha\sigma_{n}$ ($\beta=\alpha\sigma_{n}^{-1}$, $\beta=\alpha\rho_{n}$ respectively). In this case we say that the braid $\alpha$ is obtained from the braid $\beta$ by the opposite to the right stabilization of positive (negative, virtual) type transformation. S.~Kamada proved an analogue of Markov theorem for welded links in  \cite[Theorem 2]{SK}.
\begin{ttt}\label{tkam} Closures of two welded braids $\alpha$ and $\beta$ are equivalent as welded links if and only if they are related by the finite sequence of the following transformations.
\begin{enumerate}
\item A conjugation in the welded braid group,
\item A right stabilization of positive, negative or virtual type,
\item An opposite to a right stabilization of positive, negative or virtual type transformation.
\end{enumerate}
\end{ttt}
This theorem also holds for fused links since every relation of welded braid group is fulfilled in the unrestricted virtual braid group.

\begin{figure}[bh]
\noindent\centering{
\includegraphics[width=140mm]{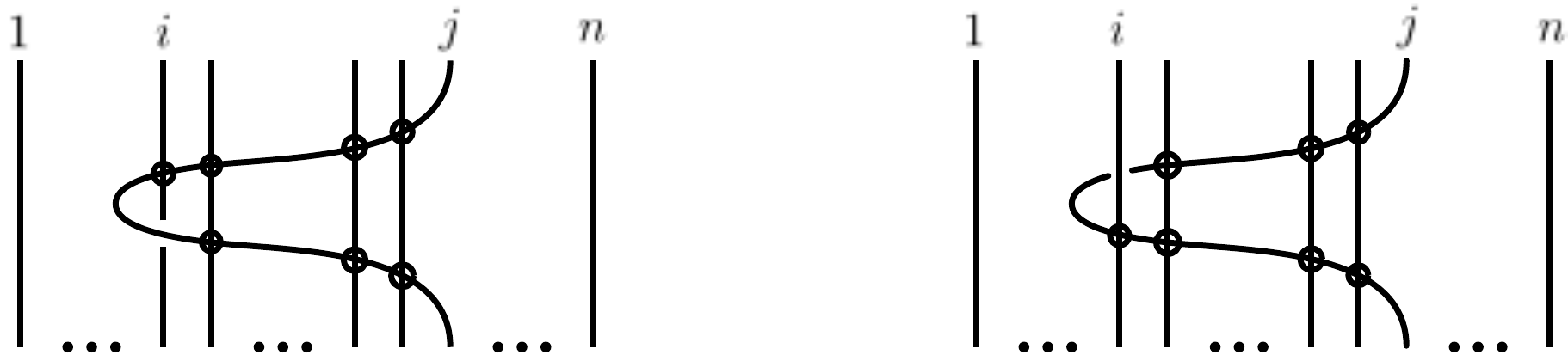}
}
\caption{Geometric braids representing $\lambda_{i,j}$ (on the left) and $\lambda_{j,i}$ (on the right)}
\label{lambdan}
\end{figure}

Let us define some element of $UVB_n$. For $i=1,\dots,n-1:$
\begin{align}\notag\lambda_{i,i+1}&=\rho_i\sigma_i^{-1},\\
\notag\lambda_{i+1,i}&=\rho_i\lambda_{i,i+1}\rho_i=\sigma_i^{-1}\rho_i.
\end{align}
For $1\leq i<j-1\leq n-1$:
$$\begin{matrix}\lambda_{i,j}&=&\rho_{j-1}\rho_{j-2}\dots\rho_{i+1}\lambda_{i,i+1}\rho_{i+1}\dots\rho_{j-2}\rho_{j-1},\\
\lambda_{j,i}&=&\rho_{j-1}\rho_{j-2}\dots\rho_{i+1}\lambda_{i+1,i}\rho_{i+1}\dots\rho_{j-2}\rho_{j-1}.
\end{matrix}$$
The elements $\lambda_{i,j}$ and $\lambda_{j,i}$ for $i<j$ belong to the pure unrestricted virtual braid group $UVP_n$ and have the following geometric interpretation (see Fig. \ref{lambdan}).
\\
The following lemma is proved in \cite[Lemma 1]{B} for the corresponding elements in $VB_n$ and therefore is also true in the quotient  $UVB_n$.
\begin{lem}\label{lem1} The following conjugating rules are fulfilled in $UVB_n$:
\begin{enumerate}
\item for $1\leq i<j\leq n$ and $k< i-1$ or $i<k<j-1$ or $k>j$:
$$\rho_k\lambda_{i,j}\rho_k=\lambda_{i,j},~~~\rho_k\lambda_{j,i}\rho_k=\lambda_{j,i};$$
\item for $1\leq i<j\leq n$:
$$\rho_{i-1}\lambda_{i,j}\rho_{i-1}=\lambda_{i-1,j},~~~\rho_{i-1}\lambda_{j,i}\rho_{i-1}=\lambda_{j,i-1};$$
\item for $1\leq i<j-1\leq n$:
$$\rho_i\lambda_{i,i+1}\rho_i=\lambda_{i+1,i},~~~\rho_i\lambda_{i,j}\rho_i=\lambda_{i+1,j},$$
$$\rho_i\lambda_{i+1,i}\rho_i=\lambda_{i,i+1},~~~\rho_i\lambda_{j,i}\rho_i=\lambda_{j,i+1};$$
\item for $1\leq i+1<j\leq n$:
$$\rho_{j-1}\lambda_{i,j}\rho_{j-1}=\lambda_{i,j-1},~~~\rho_{j-1}\lambda_{j,i}\rho_{j-1}=\lambda_{j-1,i};$$
\item for $1\leq i<j\leq n$:
$$\rho_j\lambda_{i,j}\rho_j=\lambda_{i,j+1},~~~\rho_j\lambda_{j,i}\rho_j=\lambda_{j+1,i}.$$
\end{enumerate}
\end{lem}

The following result on the structure of the pure unrestricted virtual braid group  $UVP_n$ is presented in \cite[Theorem 2.7]{BBD}.
\begin{ttt}\label{t1} The group $UVP_n$ has a presentation with generators $\lambda_{k,l}$ for $1\leq k\neq l \leq n$, and defining relations: $\lambda_{i,j}$ commutes with $\lambda_{k,l}$ if and only if $\{i,j\}\neq\{k,l\}$.
\end{ttt}

Let $U_j$ be the subgroup of $UVP_n$ generated by the  elements $\{\lambda_{i,j}, \lambda_{j,i}~|~1\leq i<j\}$. Then by Theorem \ref{t1} we have
\begin{equation}\label{eq2}
U_j=\langle\lambda_{1,j},\lambda_{j,1}\rangle\times\langle\lambda_{2,j},\lambda_{j,2}\rangle\times\dots\times\langle\lambda_{j-1,j},\lambda_{j,j-1}\rangle
=F_2\times\dots\times F_2,\end{equation}
\begin{equation}\label{eq4}
UVP_n=U_2\times U_3\times\dots\times U_n.
\end{equation}

The structure of the unrestricted virtual braid group $UVB_n$ follows from Lemma \ref{lem1} and the theorem \ref{t1} and is also given in \cite[Theorem 2.4]{BBD}.
\begin{ttt}\label{t2} The group $UVB_n$ is isomorphic to the semidirect product $UVP_n\leftthreetimes S_n$, where the symmetric group $S_n$ acts by permutations on the indices of generators $\lambda_{i,j}$ of $UVP_n$.
\end{ttt}

Denote by $B_{i,j}=\rho_{j-1}\rho_{j-2}\dots\rho_{i+1}\rho_i$ if $i<j$ and $B_{i,j}=1$ in other cases. Using simple calculations in the group $S_n$ it is easy to see that for $i<j$ the image $\iota(B_{i,j})$ is a cycle $(i~i+1\dots j)$.

\section{Construction of $\varrho^*$}

In this section we construct the map $\varrho^*:UVB_{\infty}\to UVP_{\infty}$, which is implicitly constructed in \cite{BBD}. At first we prove the following lemma.
\begin{lem}\label{lem2}Let $\alpha\in UVB_n$ and $s\in\{1,\dots,n\}$ be a maximal number, such that $\iota(\alpha)(s)\neq s$.  Then $\alpha$ can be uniquely expressed in the form
$$\alpha=\gamma x_sB_{k_s,s}x_{s+1}x_{s+2}\dots x_n, $$
where $x_i\in U_i$ and the braid $\gamma\in UVB_n$ does not involve the strands $s,s+1,\dots,n$.
\end{lem}
\textbf{Proof.} By the theorem \ref{t2} for certain elements  $\beta\in UVP_n$ and $\pi\in S_n$ we have
\begin{equation}\label{eq1}
\alpha=\beta\pi.
\end{equation}
Let $k_s=\pi(s)=\iota(\beta\pi)(s)=\iota(\alpha)(s)\neq s$. Since $\iota(\alpha)$ is a bijection of $\{1,\dots,n\}$ and $\iota(\alpha)(s)=k_s$, then we have $\iota(\alpha)(k_s)\neq k_s$. Since $s$ is a maximal number with the condition $\iota(\alpha)(s)\neq s$ then  $k_s<s$ and therefore $B_{k_s,s}\neq1$.

Since the permutations $\pi$ and $B_{k_s,s}^{-1}$ act identically on the integers $s+1,\dots,n$, then the permutation $\delta=\pi B_{k_s,s}^{-1}$ also acts identically on the integers $s+1,\dots,n$.  Moreover the image of the integer $s$ under the permutation $\delta$ is equal to $s$
$$ \begin{CD}
s @>~\pi~>> k_s @>~B_{k_s,s}^{-1}~>> s,
\end{CD}$$
therefore $\delta=\pi B_{k_s,s}^{-1}$ acts identically on the integers $s,\dots,n$.

 By the equality (\ref{eq4}) for certain elements $y_j\in U_j$ we have $\beta=y_2y_3\dots y_n$, and therefore we can rewrite the equality (\ref{eq1}).
\begin{equation}\label{eq3}
\alpha=\beta\pi= y_2y_3\dots y_n\delta B_{k_s,s}=\delta (y_2y_3\dots y_n)^{\delta} B_{k_s,s}
\end{equation}

Since $(y_2\dots y_n)$ is a pure braid, then $(y_2y_3\dots y_n)^{\delta}$ is a pure braid and therefore by the equality (\ref{eq4}) we have $(y_2y_3\dots y_n)^{\delta}=z_2\dots z_n$ for certain elements $z_j\in U_j$. The braids $z_2,\dots,z_{s-1}$ do not involve the strands $s,s+1,\dots,n$, therefore the braid $\gamma=\delta z_2\dots z_{s-1}$ does not involve the strands $s,s+1,\dots,n$. Then the equality (\ref{eq3}) can be rewritten in the following form.
\begin{align}
\notag\alpha&=\delta (y_2y_3\dots y_n)^{\delta} B_{k_s,s}=\delta z_2\dots z_{s-1}z_sz_{s+1}\dots z_n B_{k_s,s}\\
\notag&=\gamma z_sz_{s+1}\dots z_n B_{k_s,s}=\gamma z_sB_{k_s,s}(z_{s+1}\dots z_n)^{B_{k_s,s}}=\gamma z_sB_{k_s,s}z_{s+1}^{B_{k_s,s}}\dots z_n^{B_{k_s,s}}
\end{align}

Since $z_j\in U_j=\langle\lambda_{1,j},\lambda_{j,1}\rangle\times\langle\lambda_{2,j},\lambda_{j,2}\rangle\times\dots\times\langle\lambda_{j-1,j},\lambda_{j,j-1}\rangle$, then by the theorem \ref{t2} for every $j=s+1,\dots,n$ we have  $z_{j}^{B_{k_s,s}}\in U_j$. Therefore, if we denote $x_s=z_s$, $x_{s+1}=z_{s+1}^{B_{k_s,s}}$, $\dots$, $x_n=z_n^{B_{k_s,s}}$ then we have
$$\alpha=\gamma x_sB_{k_s,s}x_{s+1}\dots x_n.$$
 and we proved that the braid $\alpha$ can be written in the form from the formulation of the lemma.

To prove that such a representation is unique we consider another representation of $\alpha$ in the form from the formulation of the lemma
\begin{equation}\label{eq5}
\alpha=\gamma x_sB_{k_s,s}x_{s+1}\dots x_n=\eta y_sB_{t_s,s}y_{s+1}\dots y_n
\end{equation}
and prove that $\gamma=\eta$, $k_s=t_s$ and $x_j=y_j$ for every $j=s,\dots,n$.

From the equality (\ref{eq5}) we have
\begin{equation}\label{eq6}
\eta^{-1}\gamma = y_sB_{t_s,s}y_{s+1}\dots y_n(x_sB_{k_s,s}x_{s+1}\dots x_n)^{-1}
\end{equation}
Without loss of generality consider that $t_s\leq k_s$ and look at the images of the braids from the right and left sides of this equality under the homomorphism $\iota$.
The permutation $\iota(y_sB_{t_s,s}y_{s+1}\dots y_n(x_sB_{k_s,s}x_{s+1}\dots x_n)^{-1})$ maps $s$ to $s$ if $t_s=k_s$ and maps $s$ to $t_s$ if $t_s<k_s$
 $$ \begin{CD}
s @>~\iota(y_sB_{t_s,s}y_{s+1}\dots y_n)~>> t_s @>~\iota((x_sB_{k_s,s}x_{s+1}\dots x_n)^{-1})~>> t_s.
\end{CD}$$
 At the same time since $\eta$ and $\gamma$ does not involve the strands $s,\dots,n$, then $s$ is fixed by the permutation $\iota(\eta^{-1}\gamma)$ and therefore $t_s=k_s$. Then from the equality (\ref{eq6}) we have $\iota(y_sB_{t_s,s}y_{s+1}\dots y_n(x_sB_{k_s,s}x_{s+1}\dots x_n)^{-1})=1$ and $\iota(\eta^{-1}\gamma)=1$. Therefore $\eta^{-1}\gamma\in UVP_n$ and since it does not involve the strands $s,s+1,\dots n$, we have $$\eta^{-1}\gamma\in\langle U_2,\dots,U_{s-1}\rangle.$$

On the other side the braid $y_sB_{t_s,s}y_{s+1}\dots y_n(x_sB_{k_s,s}x_{s+1}\dots x_n)^{-1}$ can be rewritten in the following form.
\begin{align}
\notag y_sB_{t_s,s}y_{s+1}\dots y_n(x_sB_{k_s,s}x_{s+1}\dots x_n)^{-1}&=y_sB_{t_s,s}y_{s+1}\dots y_nx_n^{-1}\dots x_{s+1}^{-1}B_{t_s,s}^{-1}x_s^{-1}\\
\notag&=y_sB_{t_s,s}B_{t_s,s}^{-1}(y_{s+1}\dots y_nx_n^{-1}\dots x_{s+1}^{-1}x_s^{-1})^{B_{t_s,s}^{-1}}x_s^{-1}\\
\notag&=y_s(y_{s+1}\dots y_nx_n^{-1}\dots x_{s+1}^{-1}x_s^{-1})^{B_{t_s,s}^{-1}}x_s^{-1}
\end{align}
By the theorem \ref{t2} for every $j=s+1,\dots,n$ we have $y_j^{B_{t_s,s}^{-1}}\in U_j$, $(x_j^{-1})^{B_{t_s,s}^{-1}}\in U_j$, therefore
$$y_sB_{t_s,s}y_{s+1}\dots y_n(x_sB_{k_s,s}x_{s+1}\dots x_n)^{-1}\in \langle U_s,\dots,U_n\rangle.$$
 From the equality (\ref{eq4}) we have $\langle U_2,\dots,U_{s-1}\rangle\cap\langle U_s,\dots,U_n\rangle=1$, therefore $\eta=\gamma$ and $x_sB_{k_s,s}x_{s+1}\dots x_n=y_sB_{t_s,s}y_{s+1}\dots y_n$.
The lemma is proved.

The following lemma has a bit different formulation in the paper \cite{BBD} and it completely proved there. However here we formulate the lemma in the other words and repeat the proof from \cite{BBD} since this lemma is extremely important in our paper.
\begin{lem}\label{lem3} There exists a map $\varrho:UVB_{\infty}\to UVB_{\infty}$ with the following conditions.
\begin{enumerate}
\item For any $\alpha$ the closures of $\alpha$ and $\varrho(\alpha)$ are equivalent as double welded links.
\item If $\alpha\in UVP_{\infty}$, then $\varrho(\alpha)=\alpha$.
\item If $\alpha\in UVB_n\setminus UVP_n$, then $\varrho(\alpha)\in UVB_{n-1}$.
\end{enumerate}
\end{lem}
\textbf{Proof.} We explicitly construct the image of the braid $\alpha\in UVB_n$ under the map $\varrho$. If $\alpha\in UVP_n$, then $\varrho(\alpha)=\alpha$ and there is nothing to construct.
So let $\alpha\in UVB_n\setminus UVP_n$, then the algorithm of finding $\varrho(\alpha)$ has the following steps.
\begin{itemize}
\item[\emph{Step 1.}] Find a maximal number  $s\in\{1,\dots,n\}$ with the condition $\iota(\alpha)(s)=k_s\neq s$.
\item[\emph{Step 2.}] Conjugate the braid $\alpha$ by  the element $B_{1,n}^{n-s}$.
$$\alpha_1=B_{1,n}^{s-n}\delta B_{1,n}^{n-s}$$
The fused link $\widehat{\alpha}_1$  is equivalent to $\widehat{\alpha}$, and the permutation induced by the braid $\alpha_1$ maps $n$ to $k_n=n-s+k_s< n$
$$ \begin{CD}
n @>~B_{1,n}^{s-n}~>> s @>~\alpha~>> k_s@>~B_{1,n}^{n-s}~>> n-s+k_s=k_n.
\end{CD}$$
\item[\emph{Step 3.}] By Lemma \ref{lem2} express the braid $\alpha_1$ in the form
$$\alpha_1=\gamma w_1\dots w_{n-1}\rho_{n-1}B_{k_n,n-1},$$
where the braid $\gamma$ does not involve the strand $n$ and $w_i\in\langle\lambda_{i,n},\lambda_{n,i}\rangle$ for $i=1,\dots,n$ (see Fig. \ref{al1}).\\
\begin{figure}[bh]
\noindent\centering{
\includegraphics[width=35mm]{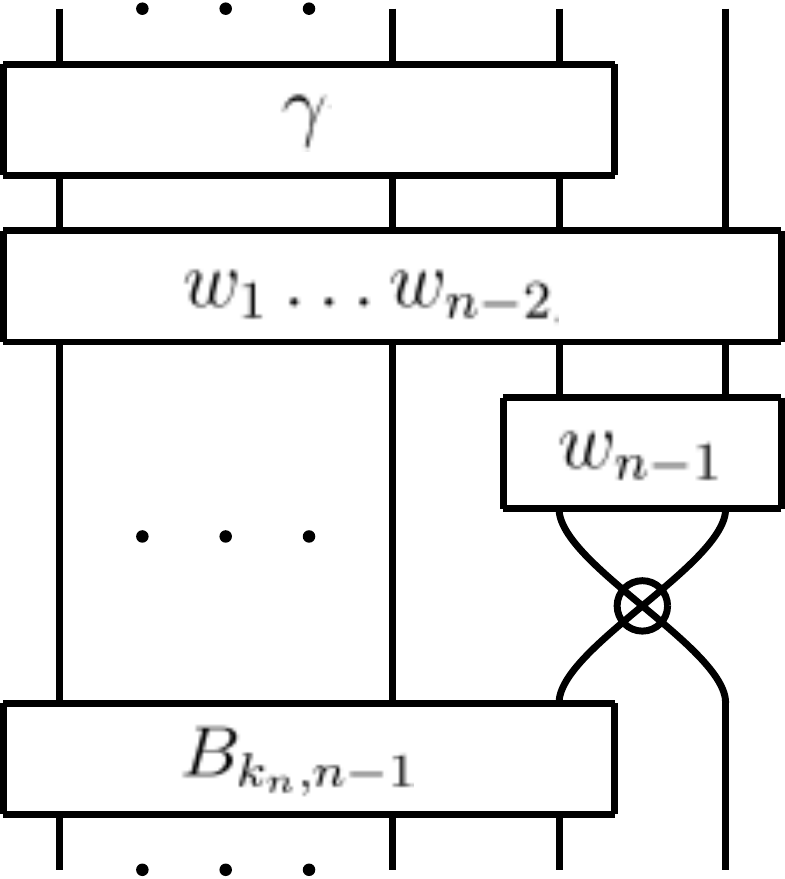}
}
\caption{The braid $\alpha_1$}
\label{al1}
\end{figure}

 By Lemma \ref{lem2} we have $\alpha_1=\gamma x_n B_{k_n,n}$, where the braid $\gamma$ does not involve an $n$-strand and $x_n\in U_n$. By the equality (\ref{eq2}) we have $x_n=w_1\dots w_{n-1}$ for $w_i\in\langle\lambda_{i,n},\lambda_{n,i}\rangle$. Sine $B_{k_n,n}\neq1$, then $B_{k_n,n}=\rho_{n-1}B_{k_n,n-1}$.
 \item[\emph{Step 4.}] Change the braid $\alpha_1$ by the braid
 $$\alpha_2=\gamma w_1\dots w_{n-2}\rho_{n-1}B_{k_n,n-1}.$$

 From the figure \ref{al1} it is obvious that all the crossings of $w_{n-1}$ belong to the same component of $\widehat{\alpha}_1$, therefore, we can virtualize all the classical crossings of $w_{n-1}\rho_{n-1}$ using Kanenobu's technique \cite[Proof of the theorem 1]{K}, i.~e. to change $\sigma_{n-1}$ to $\rho_{n-1}$ in the representation of $w_{n-1}$ and to obtain a braid $\alpha_2$ such that $\widehat{\alpha}_2$ and $\widehat{\alpha}_1$ are equivalent as fused links.

 Since $w_{n-1}\in \langle\lambda_{n-1,n},\lambda_{n,n-1}\rangle$ and $\lambda_{n-1,n}=\rho_{n-1}\sigma_{n-1}^{-1}$, $\lambda_{n,n-1}=\sigma_{n-1}^{-1}\rho_{n-1}$, then after virtualization of classical crossings the braid $w_{n-1}\rho_{n-1}$ involved an odd number of $\rho_{n-1}$ and we obtain $\rho_{n-1}$ instead of $w_{n-1}\rho_{n-1}$. Therefore the closure of the braid $\alpha_2=\gamma w_1\dots w_{n-2}\rho_{n-1}B_{k_n,n-1}$ is equivalent to the closure of the braid $\alpha_1$ as fused link.
 \item[\emph{Step 5.}] Construct $\varrho(\alpha)$.

 Note that by Lemma \ref{lem1} for all $\lambda_{j,n}$ and $\lambda_{n,j}$ with $j<n-1$ we have:
$$\lambda_{j,n}\rho_{n-1}=\rho_{n-1}\lambda_{j,n-1},$$
$$\lambda_{n,j}\rho_{n-1}=\rho_{n-1}\lambda_{n-1,j}.$$
Therefore  $\alpha_2=\gamma \rho_{n-1}w^{\prime}_1\dots w^{\prime}_{n-2}B_{k_n,n-1}$, where $w_i^{\prime}=w_i^{\rho_{n-1}}$ is a braid from $\langle\lambda_{i,n-1},\lambda_{n-1,i}\rangle$ and hence does not involve an $n$-strand (see Fig. \ref{al2}).

\begin{figure}[bh]
\noindent\centering{
\includegraphics[width=35mm]{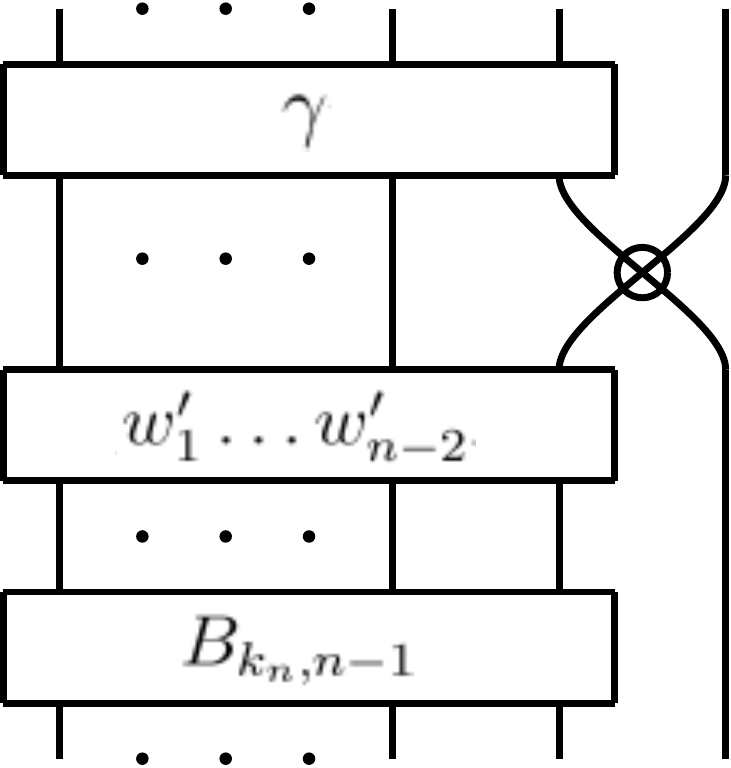}
}
\caption{The braid $\alpha_2$}
\label{al2}
\end{figure}

In the braid $\alpha_2$ there is only one (virtual) crossing on the $n$-strand, so, using Markov moves we obtain a new braid $\alpha_3=\gamma w^{\prime}_1\dots w^{\prime}_{n-2}B_{k_n,n-1}$ whose closure is again equivalent to $\widehat{\alpha}$ and has $n-1$ strands. We have constructed the braid $$\varrho(\alpha)=\alpha_3=\gamma w^{\prime}_1\dots w^{\prime}_{n-2}B_{k_n,n-1}.$$
\end{itemize}
The lemma is proved.
\begin{rmk} It is obvious that  the equality
 $\varrho(\alpha)=\varrho(\alpha_1)=\varrho(\alpha_2)$ holds in Lemma \ref{lem2}.
\end{rmk}
Let $\alpha\in UVB_n$, then the sequence $\varrho(\alpha),\varrho^2(\alpha),\dots$ is stabilized on some step. Let $\varrho^*$ be such a map which maps the braid $\alpha$ to the braid $\varrho^k(\alpha)$, where $k$ is a minimal integer such that $\varrho^k(\alpha)=\varrho^{k+1}(\alpha)$. Every step of finding $\varrho(\alpha)$ is clearly defined, therefore $\varrho^*$ is a correct function.

Using simple induction on the number $n$  it is easy to prove the following fact.
\begin{rmk}\label{rmk2} Let $\alpha\in UVB_n$ be the braid, such that $\varrho^*(\alpha)\in UVP_m$. If $\alpha=\beta\gamma$ for $\beta\in UVB_n$, $\gamma\in UVP_n^{\prime}$, then $\varrho^*(\alpha)=\varrho^*(\beta)\delta$, where $\delta\in UVP_m^{\prime}$.
\end{rmk}
From Lemma \ref{lem3} We have the following corollary which is formulated and proved in the paper \cite{BBD}.
\begin{prp} Let $\alpha\in UVB_n$ be an unrestricted virtual braid, such that its closure $\widehat{\alpha}$ has $m$ components. Then there exists a pure unrestricted virtual braid $\beta\in UVP_m$, such that $\widehat{\alpha}=\widehat{\beta}$.
\end{prp}
\textbf{Proof.} $\beta=\varrho^*(\alpha)$.

\section{Proof of the main results}

\begin{lem}\label{lem4} Let $\alpha\in UVP_n$ and $u,v\in \langle\lambda_{n-1,n},\lambda_{n,n-1}\rangle$. Then the closure of the braids $\alpha$ and $\alpha[u,v]$ are equivalent as fused links.
\end{lem}
\textbf{Proof.} Let $\gamma$ be the following braid from $UVB_{n+2}$
$$\gamma=\alpha u^{-1}\rho_{n-1}u^{\rho_{n}\rho_{n-1}}B_{n-1,n+1}B_{n-1,n+2}.$$
Here the braid $\alpha u^{-1}\rho_{n-1}$ does not involve the strands $n+1,n+2$ and the braid $u^{\rho_{n}\rho_{n-1}}$ belongs to $\langle\lambda_{n,n+1},\lambda_{n+1,n}\rangle$.
Let us find the braid $\varrho(\gamma)$.
\begin{itemize}
\item[Step 1.] Since $\iota(\gamma)(n+2)=n-1$, then the maximal number which is not fixed by $\gamma$ is equal to $n+2$.
\item[Step 2.]  $\gamma_1=\gamma$.
\item[Step 3.] $\gamma=\alpha u^{-1}\rho_{n-1}u^{\rho_{n}\rho_{n-1}}B_{n-1,n+1}B_{n-1,n+2}$, where the braid $\alpha u^{-1}\rho_{n-1}u^{\rho_{n}\rho_{n-1}}B_{n-1,n+1}$ does not involve the strand $n+2$.
\item[Step 4.]  $\gamma_2=\alpha u^{-1}\rho_{n-1}u^{\rho_{n}\rho_{n-1}}B_{n-1,n+1}\rho_{n+1}B_{n-1,n+1}$.
\item[Step 5.] $\varrho(\gamma)=\alpha u^{-1}\rho_{n-1}u^{\rho_{n}\rho_{n-1}}B_{n-1,n+1}B_{n-1,n+1}$.
\end{itemize}
Let us find the braid $\varrho^2(\gamma)$.
\begin{itemize}
\item[Step 1.] Since $\iota(\varrho(\gamma))(n+1)=n$, then the maximal number which is not fixed by $\varrho(\gamma)$ is equal to $n+1$.
\item[Step 2.]  $\varrho(\gamma)_1=\varrho(\gamma)$.
\item[Step 3.] The braid $\varrho(\gamma)$ can be rewritten
\begin{align}
\notag\varrho(\gamma)&=\alpha u^{-1}\rho_{n-1}u^{\rho_{n}\rho_{n-1}}B_{n-1,n+1}B_{n-1,n+1}\\
\notag&=\alpha u^{-1}\rho_{n-1}u^{\rho_{n}\rho_{n-1}}\rho_{n-1}\rho_n\\
\notag&=\alpha u^{-1}\rho_{n-1}u^{\rho_{n}\rho_{n-1}}\rho_{n-1}\rho_n=\alpha u^{-1}u^{\rho_{n}}\rho_n,
\end{align}
where the braid $\alpha u^{-1}$ does not involve the strand $n+1$ and $u^{\rho_{n}}$ belongs to $\langle\lambda_{n-1,n+1},\lambda_{n+1,n-1}\rangle$.
\item[Step 4.]  $\varrho(\gamma)_2=\alpha u^{-1}u^{\rho_{n}}\rho_n$.
\item[Step 5.] $\varrho^2(\gamma)=\alpha$.
\end{itemize}
Since $\varrho^2(\gamma)=\alpha$ is a pure braid, then $\varrho^*(\gamma)=\alpha$. Let $w=v^{\rho_{n-1}B_{n-1,n+1}B_{n-1,n+2}}$, then by Lemma \ref{lem1} the braid $w$ belongs to $\langle\lambda_{n+2,n+1},\lambda_{n+1,n+2}\rangle$. Denote by $\delta=w\gamma w^{-1}$ and find the braid $\varrho^*(\delta)$. At first find the braid $\varrho(\delta)$.
\begin{itemize}
\item[Step 1.] Since $w$ is a pure braid, then the maximal number, which is not fixed by $\delta$ is equal to the maximal number which is not fixed by $\gamma$ and is equal to $n+2$.
\item[Step 2.]  $\delta_1=\delta$.
\item[Step 3.] The braid $\delta$ can be rewritten in details
\begin{align}
\notag \delta&=w\gamma w^{-1}=w\alpha u^{-1}\rho_{n-1}u^{\rho_{n}\rho_{n-1}}B_{n-1,n+1}B_{n-1,n+2}w^{-1}\\
\notag&=\alpha u^{-1}\rho_{n-1}u^{\rho_{n}\rho_{n-1}}wB_{n-1,n+1}(w^{-1})^{B_{n-1,n+2}^{-1}}B_{n-1,n+2}\\
\notag&=\alpha u^{-1}\rho_{n-1}u^{\rho_{n}\rho_{n-1}}B_{n-1,n+1}w^{B_{n-1,n+1}}(w^{-1})^{B_{n-1,n+2}^{-1}}B_{n-1,n+2}\\
\notag&=\alpha u^{-1}\rho_{n-1}u^{\rho_{n}\rho_{n-1}}B_{n-1,n+1}(w^{-1})^{B_{n-1,n+2}^{-1}}w^{B_{n-1,n+1}}B_{n-1,n+2},
\end{align}
where the braid $\alpha u^{-1}\rho_{n-1}u^{\rho_{n}\rho_{n-1}}B_{n-1,n+1}(w^{-1})^{B_{n-1,n+2}^{-1}}$ does not involve the strand $n+2$ and the braid $w^{B_{n-1,n+1}}$ belongs to $\langle\lambda_{n+2,n},\lambda_{n,n+2}\rangle$.
\item[Step 4.]  $\delta_2=\alpha u^{-1}\rho_{n-1}u^{\rho_{n}\rho_{n-1}}B_{n-1,n+1}(w^{-1})^{B_{n-1,n+2}^{-1}}w^{B_{n-1,n+1}}\rho_{n+1}B_{n-1,n+1}$.
\item[Step 5.] $\varrho(\gamma)=\alpha u^{-1}\rho_{n-1}u^{\rho_{n}\rho_{n-1}}B_{n-1,n+1}(w^{-1})^{B_{n-1,n+2}^{-1}}w^{B_{n-1,n+1}\rho_{n+1}}B_{n-1,n+1}$.
\end{itemize}
Let us fin the braid $\varrho^2(\delta)$
\begin{itemize}
\item[Step 1.] Since $\iota(\varrho(\delta))(n+1)=n$, then $n+1$ is maximal number which is not fixed by $\varrho(\delta)$.
\item[Step 2.]  $\varrho(\delta)_1=\varrho(\delta)$.
\item[Step 3.] The braid $\varrho(\delta)$ can be rewritten in details
\begin{align}
\notag \varrho(\delta)&=\alpha u^{-1}\rho_{n-1}u^{\rho_{n}\rho_{n-1}}B_{n-1,n+1}(w^{-1})^{B_{n-1,n+2}^{-1}}w^{B_{n-1,n+1}\rho_{n+1}}B_{n-1,n+1}\\
\notag&=\alpha u^{-1}\rho_{n-1}u^{\rho_{n}\rho_{n-1}}(w^{-1})^{B_{n-1,n+2}^{-1}B_{n-1,n+1}^{-1}}w^{B_{n-1,n+1}\rho_{n+1}B_{n-1,n+1}^{-1}}B_{n-1,n+1}B_{n-1,n+1}\\
\notag&=\alpha u^{-1}\rho_{n-1}u^{\rho_{n}\rho_{n-1}}(w^{-1})^{B_{n-1,n+2}^{-1}B_{n-1,n+1}^{-1}}w^{B_{n-1,n+1}\rho_{n+1}B_{n-1,n+1}^{-1}}\rho_{n-1}\rho_n\\
\notag&=\alpha u^{-1}u^{\rho_{n}}(w^{-1})^{B_{n-1,n+2}^{-1}B_{n-1,n+1}^{-1}\rho_{n-1}}w^{B_{n-1,n+1}\rho_{n+1}B_{n-1,n+1}^{-1}\rho_{n-1}}\rho_n\\
\notag&=\alpha u^{-1}(w^{-1})^{B_{n-1,n+2}^{-1}B_{n-1,n+1}^{-1}\rho_{n-1}}u^{\rho_{n}}w^{B_{n-1,n+1}\rho_{n+1}B_{n-1,n+1}^{-1}\rho_{n-1}}\rho_n,
\end{align}
where the braid $\alpha u^{-1}(w^{-1})^{B_{n-1,n+2}^{-1}B_{n-1,n+1}^{-1}\rho_{n-1}}$ does not involve the strand $n+1$ and the braid $u^{\rho_n}ww^{B_{n-1,n+1}\rho_{n+1}B_{n-1,n+1}^{-1}\rho_{n-1}}$ belongs to $\langle\lambda_{n-1,n+1},\lambda_{n+1,n-1}\rangle$.
\item[Step 4.]  $\varrho(\delta)_2=\alpha u^{-1}(w^{-1})^{B_{n-1,n+2}^{-1}B_{n-1,n+1}^{-1}\rho_{n-1}}u^{\rho_{n}}w^{B_{n-1,n+1}\rho_{n+1}B_{n-1,n+1}^{-1}\rho_{n-1}}\rho_n$.
\item[Step 5.] The braid $\varrho^2(\delta)$ follows
$$\varrho^2(\delta)=\alpha u^{-1}(w^{-1})^{B_{n-1,n+2}^{-1}B_{n-1,n+1}^{-1}\rho_{n-1}}uw^{B_{n-1,n+1}\rho_{n+1}B_{n-1,n+1}^{-1}\rho_{n-1}\rho_n}.$$
Since $w=v^{\rho_{n-1}B_{n-1,n+1}B_{n-1,n+2}}$, then $(w^{-1})^{B_{n-1,n+2}^{-1}B_{n-1,n+1}^{-1}\rho_{n-1}}=v^{-1}$ and
\begin{align}
\notag w^{B_{n-1,n+1}\rho_{n+1}B_{n-1,n+1}^{-1}\rho_{n-1}\rho_n}&=v^{\rho_{n-1}B_{n-1,n+1}B_{n-1,n+2}B_{n-1,n+1}\rho_{n+1}\underline{B_{n-1,n+1}^{-1}}\rho_{n-1}\rho_n}\\
\notag&=v^{\rho_{n-1}B_{n-1,n+1}B_{n-1,n+2}B_{n-1,n+1}\rho_{n+1}\rho_{n-1}\underline{\rho_n\rho_{n-1}\rho_n}}\\
\notag&=v^{\rho_{n-1}\underline{B_{n-1,n+1}B_{n-1,n+2}}B_{n-1,n+1}\rho_{n+1}\rho_{n}\rho_{n-1}}\\
\notag&=v^{\rho_{n-1}B_{n-1,n+2}\underline{B_{n,n+2}B_{n-1,n+1}}\rho_{n+1}\rho_{n}\rho_{n-1}}\\
\notag&=v^{\rho_{n-1}B_{n-1,n+2}\rho_{n+1}\rho_{n-1}\rho_{n+1}\rho_{n}\rho_{n-1}}\\
\notag&=v^{\rho_{n-1}\underline{B_{n-1,n+2}\rho_{n-1}\rho_{n}}\rho_{n-1}}=v^{\rho_{n-1}\rho_{n+1}\rho_{n-1}}=v^{\rho_{n+1}}\\
\notag&=v~~~~~~~~~~~~~~~~(\text{since~}v\in\langle\lambda_{n-1,n},\lambda_{n,n-1}\rangle)
\end{align}
Therefore $\varrho^2(\delta)=\alpha u^{-1}v^{-1}uv=\alpha[u,v]$. Since $\varrho^2(\delta)$ is a pure braid, then $\varrho^*(\delta)=\varrho^2(\delta)=\alpha[u,v]$.
\end{itemize}
Since the closures of the braids $\gamma$ and $\delta=\gamma^{w^{-1}}$ define equivalent fused links, then the closures of the braids $\varrho^*(\gamma)=\alpha$ and $\varrho^*(\delta)=\alpha[u,v]$ define equivalent fused links. The lemms is proved.

Denote by $T$ the set of coset representatives of $UVP_n/UVP_n^{\prime}$ and for the element $\alpha\in UVP_n$ denote by $\overline{\alpha}\in T$ the unique coset representative $\alpha UVP_n^{\prime}=\overline{\alpha}UVP_n^{\prime}$. The following statement  almost immediately follows from Lemma \ref{lem4}.
\begin{cor}\label{cor1}Let $\alpha,\beta\in UVP_{n}$ be unrestricted virtual braids, such that $\overline{\alpha}$ and $\overline{\beta}$ are conjugated by the element from $S_n$. Then the fused links $\widehat{\alpha}$ and $\widehat{\beta}$ are equivalent.
\end{cor}
\textbf{Proof.} Since $\overline{\alpha}$ and $\overline{\beta}$ are conjugated by the element from $S_n$ then the closures of the braids $\overline{\alpha}$ and $\overline{\beta}$ are equivalent fused links. Therefore it is enough to prove that the closure of the braid $\alpha$ is equivalent to the closure of the braid $\overline{\alpha}$ and the closure of the braid $\beta$ is equivalent to the closure of the braid $\overline{\beta}$.

Since $\overline{\alpha}UVP_n^{\prime}=\alpha UVP_n^{\prime}$, then $\overline{\alpha}=\alpha[u_1,v_1][u_2,v_2]\dots[u_k,v_k]$ for certain elements $u_i,v_i\in UVP_n$, $i=1,\dots,k$. By Theorem \ref{t1} we can consider that $u_i,v_i\in\langle\lambda_{r_i,s_i},\lambda_{s_i,r_I}\rangle$.
Denote by $\alpha_0=\alpha$, $\alpha_1=\alpha_0[u_1,v_1]$, $\alpha_2=\alpha_1[u_2,v_2]$, $\dots$, $\alpha_k=\alpha_{k-1}[u_k,v_k]=\overline{\alpha}$. If we prove that the closure of $\alpha_i$ and the closure of $\alpha_{i-1}$ define the same fused link for every $i=1,\dots,k$, then we prove that the closures of the braids $\alpha$ and $\overline{\alpha}$ are equivalent as fused links.
Therefore we can consider that $\overline{\alpha}=\alpha[u,v]$ for some $u,v\in\langle\lambda_{r,s},\lambda_{s,r}\rangle$.

The closure of the braid $\overline{\alpha}$ is equivalent to the closure of the braid $\overline{\alpha}^{\mu}$ for any $\mu\in S_n<UVB_n$. If $\mu$ is a permutation which maps $s$ to $n-1$ and maps $r$ to $n$, then the closure of the braid $\overline{\alpha}$ is equivalent to the closure of the braid $\overline{\alpha}^{\mu}=\alpha^{\mu}[u^{\mu},v^{\mu}]$. By Theorem \ref{t2} the braids $u^{\mu},v^{\mu}$ belong to $\langle\lambda_{n-1,n},\lambda_{n,n-1}\rangle$, thus by Lemma \ref{lem4} the closure of the braid $\overline{\alpha}$ is equivalent to the closure of the braid $\alpha^{\mu}$ and is also equivalent to the closure of $\alpha$. The closure of $\beta$ is equivalent to the closure of $\overline{\beta}$ by the same reasons.
The corollary is proved.

\begin{prp}\label{pr1}Let $\alpha, \beta\in UVB_{\infty}$ be unrestricted virtual braids, such that $\widehat{\alpha}$ and $\widehat{\beta}$ are equivalent fused links with $m$ components. Then $\overline{\varrho^*(\alpha)}$ and $\overline{\varrho^*(\beta)}$ are conjugated by the element from $S_n$.
\end{prp}
\textbf{Proof.} Since the links $\widehat{\alpha}$ and $\widehat{\beta}$ are equivalent, then by Theorem \ref{tkam} the braids $\alpha$ and $\beta$ are related by the finite sequence of Markov's transformations. It means that there is a finite sequence of braids $$\alpha=\alpha_0,\alpha_1,\dots,\alpha_k=\beta$$
such that $\alpha_{j+1}$ is obtained from $\alpha_j$ by conjugation in $UVB_{\infty}$, by right stabilization or by inverse to right stabilization transformation.

If we prove that $\overline{\varrho^*(\alpha_j)}$ and $\overline{\varrho^*(\alpha_{j+1})}$ are conjugated by the element from $S_n$ for every $j=0,\dots,k-1$, than we prove that $\overline{\varrho^*(\alpha)}$ and $\overline{\varrho^*(\beta)}$ are conjugated by the element from $S_n$, therefore we can consider that $\beta$ is obtained from $\alpha$ using only one Markov's transformation.

\underline{Case 1.} The braid $\beta$ is obtained from the braid $\alpha$ by right stabilization of positive, negative or virtual type.
We consider only the case of right stabilization of positive type ($\beta=\alpha\sigma_n$), the cases of right stabilization of  negative and virtual type are similar.

Let us count the image of $\beta$ under the map $\varrho$:
\begin{itemize}
\item[Step 1.] Since $\beta=\alpha\sigma_n$, then $\iota(\beta)(n+1)=\iota(\alpha\sigma_n)(n+1)=n\neq n+1$
      $$ \begin{CD}
n+1 @>~\alpha~>> n+1 @>~\sigma_n~>> n,
\end{CD}$$
     therefore $n+1$ is a maximal number which is not fixed by $\beta$.

\item[Step 2.]  $\beta_1=B_{1,n+1}^{n+1-(n+1)}\beta B_{1,n+1}^{n+1-(n+1)}=\beta$.

\item[Step 3.] We can rewrite the braid $\beta_1$ in the following form $$\beta_1=\beta=\alpha\sigma_n=\alpha\sigma_n\rho_n\rho_n=\alpha\lambda_{n,n+1}^{-1}\rho_n,$$
where the braid $\alpha$ does not involve the strand $n+1$ and $\lambda_{n,n+1}^{-1}\in\langle\lambda_{n,n+1},\lambda_{n+1,n}\rangle$.

\item[Step 4.]  $\beta_2=\alpha\rho_n$.

\item[Step 5.] $\varrho(\beta)=\alpha$.
\end{itemize}
Since $\varrho(\beta)=\alpha$, then $\varrho^*(\beta)=\varrho^*(\alpha)$ (and they the braids $\overline{\varrho^*(\beta)}$ and $\overline{\varrho^*(\alpha)}$ are certainly conjugated).

\underline{Case 2.} The braid $\beta$ is obtained from the braid $\alpha$ by an opposite to a right stabilization of positive, negative or virtual type transformation.

In this case the braid $\alpha$ is obtained from the braid $\beta$ by right stabilization of positive (negative, virtual) type and by the case 1 we have $\varrho^*(\alpha)=\varrho^*(\beta)$.

\underline{Case 3.} The braid $\beta$ is obtained from the braid $\alpha$ by conjugation in $UVB_n$.

%In this case it is enough to prove that for conjugated braids $\alpha$ and $\beta$ the braids $\varrho^j(\alpha)$ and $\varrho^j(\beta)$ are conjugated for some $j$.
%Really if the braids $\varrho^j(\alpha)$ and $\varrho^j(\beta)$ are conjugated for some $j$, then denote $\delta_1=\varrho^j(\alpha)$, $\gamma_1=\varrho^j(\beta)$. Since $\delta_1$ and $\gamma_1$ are conjugated, then $\varrho^{j_1}(\delta_1)$ and $\varrho^{j_1}(\gamma_1)$ are conjugated for some $j_1$. But $\varrho^{j_1}(\delta_1)=\varrho^{j+j_1}(\alpha)$, $\varrho^{j_1}(\gamma_1)=\varrho^{j+j_1}(\beta)$. Denoting by $\delta_2=\varrho^{j_1}(\delta_1)$, $\gamma_2=\varrho^{j_1}(\gamma_1)$ and repeating this process enough number of times we conclude that $\varrho^*(\alpha)$ and $\varrho^*(\beta)$ are conjugated.

Since the braid $\alpha$ and $\beta$ are conjugated, then for some braid $\theta$ we have $\beta=\alpha^{\theta}$. Since $\sigma_i=\rho_{i}\lambda_{i,i+1}$, then $\theta=y_1y_2\dots y_r$, where $$y_i\in\{\rho_1,\dots,\rho_{n-1},\lambda_{1,2}^{\pm1},\lambda_{2,3}^{\pm1},\dots,\lambda_{n-1,n}^{\pm1}\}.$$
If we denote $\delta_1=\alpha$, $\delta_2=\delta_1^{y_1}$, $\delta_3=\delta_2^{y_2}$, $\dots$, $\delta_{r+1}=\delta_r^{y_r}=\beta$ and prove that $\overline{\varrho^*(\delta_j)}$ is conjugated with $\overline{\varrho^*(\delta_{j+1})}$ by the element from $S_n$ for every $j=1,
\dots,r+1$, then we prove that $\overline{\varrho^*(\alpha)}$ is conjugated with $\overline{\varrho^*(\beta)}$ by the element from $S_n$. Therefore we can consider that $\beta$ is obtained from $\alpha$ conjugating by $\rho_i$ or $\lambda_{i,i+1}^{\pm1}$ for some $i$.

We use the induction by the parameter $n-m$, i.~e. by the difference between the number of strands in the braid $\alpha$ and the number components in the link $\widehat{\alpha}$.

If $n-m=0$, then $\alpha$ and $\beta$ are pure braids and $\varrho^*(\alpha)=\alpha$, $\varrho^*(\beta)=\beta$.
Let $\beta=\alpha^{\mu}$ for some braid $\mu\in UVB_n$, then by Theorem \ref{t2} for certain elements $\beta\in UVP_n$ and $\pi\in S_n$ the braid $\mu$ can be presented as a product $\mu=\beta\pi$. Therefore
\begin{align}
\notag\varrho^*(\beta)&=\beta=\alpha^{\mu}=\mu^{-1}\alpha\mu=\pi^{-1}\beta^{-1}\alpha\beta\pi\\
\notag&=\pi^{-1}\alpha\beta^{-1}[\beta^{-1},\alpha]\beta\pi=\pi^{-1}\alpha[\beta^{-1},\alpha]^{\beta}\pi=\pi^{-1}\varrho^*(\alpha)[\beta^{-1},\alpha^{\beta}]\pi
\end{align}
Since $\alpha,\beta\in UVP_n$, then $[\beta^{-1},\alpha^{\beta}]\in UVP_n^{\prime}$, therefore $\overline{\varrho^*(\alpha)[\beta^{-1},\alpha^{\beta}]}=\overline{\varrho^*(\alpha)}$ and $\overline{\varrho^*{(\beta)}}=\overline{\varrho^*(\alpha)}^{\pi}$.

If $n-m>0$, then $\alpha,\beta\in UVB_n\setminus UVP_n$, $\varrho^*(\alpha)\neq\alpha$, $\varrho^*(\beta)\neq\beta$ and we have the following cases

\underline{Case 3.1.} The braid $\beta$ is obtained from $\alpha$ conjugating by $\rho_i$. Let us find $\varrho(\alpha)$.
\begin{itemize}
\item[Step 1.] Let $s$ be a maximal number which is not fixed by $\alpha$ and $\iota(\alpha)(s)=k_s$.

\item[Step 2.] Since $\alpha_1=B_{1,n}^{s-n}\alpha B_{1,n}^{n-s}$, then $\iota(\alpha_1)(n)=k_s+n-s=k_n\neq n$.
$$ \begin{CD}
n @>~B_{1,n}^{s-n}~>> s @>~\alpha~>> k_s@>~B_{1,n}^{n-s}~>>k_s+n-s=k_n.
\end{CD}$$
\item[Step 3.] Express the braid $\alpha_1$ in the form
\begin{equation}\label{a1}
\alpha_1=B_{1,n}^{s-n}\alpha B_{1,n}^{n-s}=\gamma w_1\dots w_{n-1}\rho_{n-1}B_{k_n,n-1},
\end{equation}
where the braid $\gamma$ does not involve the strand $n$ and $w_i\in\langle\lambda_{i,n},\lambda_{n,i}\rangle$.
\item[Step 4.] $\alpha_2=\gamma w_1\dots w_{n-2}\rho_{n-1}B_{k_n,n-1}.$

\item[Step 5.] The braid $\varrho(\alpha)$ has the following form
\begin{equation}\label{a2}
\varrho(\alpha)=\gamma w_1^{\rho_{n-1}}\dots w_{n-2}^{\rho_{n-1}}B_{k_n,n-1},
\end{equation}
where $w_{j}^{\rho_{n-1}}\in U_{n-1}$  does not involve the strand $n$.
\end{itemize}

\underline{Case 3.1.1.} $i\leq k_s-2$. Let us count the braid $\varrho(\beta)=\varrho(\rho_i\alpha\rho_i)$.

\begin{itemize}
\item[Step 1.] Since $i\leq k_s-2$, then $\rho_i$ fixes $s$ and $k_s$, therefore $s$ is a maximal number which is not fixed by $\beta$ and $\iota(\beta)(s)=k_s$.
        $$ \begin{CD}
s @>~\rho_{i}~>> s @>~\alpha~>> k_s@>~\rho_{i}~>> k_s
\end{CD}$$
\item[Step 2.] We have $\beta_1=B_{1,n}^{s-n}\beta B_{1,n}^{n-s}$, then the permutation $\iota(\beta)$ maps $n$ to $k_n\neq n$.
    $$ \begin{CD}
n @>~B_{1,n}^{s-n}~>> s @>~\beta~>> k_s@>~B_{1,n}^{n-s}~>> k_s+n-s=k_n
\end{CD}$$
\item[Step 3.] For $i\leq k_s-2$ we have
\begin{align}
\notag\rho_iB_{1,n}&=\rho_i\rho_{n-1}\dots\rho_1=\rho_i\rho_{n-1}\dots\rho_{i+2}\rho_{i+1}\rho_i\rho_{i-1}\dots\rho_1\\
\notag&=\rho_{n-1}\dots\rho_{i+2}\underline{\rho_i\rho_{i+1}\rho_i}\rho_{i-1}\dots\rho_1\\
\notag&=\rho_{n-1}\dots\rho_{i+2}\rho_{i+1}\rho_{i}\rho_{i+1}\rho_{i-1}\dots\rho_1\\
\notag&=\rho_{n-1}\dots\rho_{i+2}\rho_{i+1}\rho_{i}\rho_{i-1}\dots\rho_1\rho_{i+1}=B_{1,n}\rho_{i+1},
\end{align}
therefore $\rho_iB_{1,n}^{n-s}=B_{1,n}^{n-s}\rho_{i+n-s}$ and hence we have
\begin{align}
\notag\beta_1&=B_{1,n}^{s-n}\beta B_{1,n}^{n-s}=B_{1,n}^{s-n}\rho_i\alpha\rho_i B_{1,n}^{n-s}=\rho_{i+n-s}B_{1,n}^{s-n}\alpha B_{1,n}^{n-s}\rho_{i+n-s}\\
\notag&=\rho_{i+n-s}\gamma w_1\dots w_{n-1}\rho_{n-1}B_{k_n,n-1}\rho_{i+n-s}\\
\notag&=\gamma^{\rho_{i+n-s}} (w_1\dots w_{n-2})^{\rho_{i+n-s}}w_{n-1}^{\rho_{i+n-s}}\rho_{n-1}^{\rho_{i+n-s}}B_{k_n,n-1}^{\rho_{i+n-s}}
\end{align}

Since $i\leq k_s-2$, then $i+n-s\leq k_s-2+n-s=k_n-2<n-2$, therefore $B_{k_n,n-1}^{\rho_{i+n-s}}=B_{k_n,n-1}$, $w_{n-1}^{\rho_{i+n-s}}=w_{n-1}$, $\rho_{n-1}^{\rho_{i+n-s}}=\rho_{n-1}$ and the braid
$\gamma^{\rho_{i+n-s}}$ does not involve the strand $n$. Since $\rho_{i+n-s}$ fixes $n$, then
$(w_1\dots w_{n-2})^{\rho_{i+n-s}}$ belongs to $\langle\lambda_{1,n},\lambda_{n,1}\rangle\times\dots\times\langle\lambda_{n-2,n},\lambda_{n,n-2}\rangle$. Then we have
$$\beta_1=\gamma^{\rho_{i+n-s}} (w_1\dots w_{n-2})^{\rho_{i+n-s}}w_{n-1}\rho_{n-1}B_{k_n,n-1}.$$
\item[Step 4.] $\beta_2=\gamma^{\rho_{i+n-s}} (w_1\dots w_{n-2})^{\rho_{i+n-s}}\rho_{n-1}B_{k_s,n-1}$

\item[Step 5.] Since $i+n-s<n-2$, then $\rho_{n-1}$ and $\rho_{i+n-s}$ commute, therefore the braid $\varrho(\beta)$ has the following form
\begin{align}
\notag\varrho(\beta)&=\gamma^{\rho_{i+n-s}} (w_1\dots w_{n-2})^{\rho_{i+n-s}\rho_{n-1}}B_{k_n,n-1}\\
\notag&=\gamma^{\rho_{i+n-s}} (w_1\dots w_{n-2})^{\rho_{n-1}\rho_{i+n-s}}B_{k_n,n-1}\\
\notag &=\gamma^{\rho_{i+n-s}} ((w_1\dots w_{n-2})^{\rho_{n-1}})^{\rho_{i+n-s}}B_{k_n,n-1}^{\rho_{i+n-s}}\\
\notag &=(\gamma(w_1\dots w_{n-2})^{\rho_{n-1}}B_{k_n,n-1})^{\rho_{i+n-s}}=\varrho(\alpha)^{\rho_{i+n-s}}
\end{align}
\end{itemize}
Therefore $\varrho(\beta)=\varrho(\alpha)^{\rho_{i+n-s}}$ and by
the induction hypothesis the braids $\overline{\varrho^*(\alpha)}$ and $\overline{\varrho^*(\beta)}$ are conjugated by the element from $S_n$.

\underline{Case 3.1.2.} $i= k_s-1$. Let us find the braid $\varrho(\beta)$.
\begin{itemize}
\item[Step 1.]

Since $k_s\leq s-1$, then $i=k_s-1\leq s-2$ and $\rho_i$ fixes the strands $s$. Then the image of $s$ under the permutation
 $\beta=\rho_{k_s-1}\alpha\rho_{k_s-1}$ is equal to $k_s-1$.
$$ \begin{CD}
s @>~\rho_{k_s-1}~>> s @>~\alpha~>> k_s@>~\rho_{k_s-1}~>> k_s-1.
\end{CD}$$
Therefore $s$ is a maximal number  which is not foxed by $\beta$.

\item[Step 2.]  We have $\beta_1=B_{1,n}^{s-n}\beta B_{1,n}^{n-s}$ and $\iota(\beta_1)(n)=k_n-1$.
$$ \begin{CD}
n @>~B_{1,n}^{s-n}~>> s @>~\beta~>> k_s-1@>~ B_{1,n}^{n-s}~>> n-s+k_s-1=k_n-1
\end{CD}$$
\item[Step 3.] Using the same arguments as in the case 3.1.1 we have
$$\rho_{k_s-1}B_{1,n}^{n-s}=B_{1,n}^{n-s}\rho_{n-s+k_s-1}=B_{1,n}^{n-s}\rho_{k_n-1}.$$
Therefore the braid $\beta_1$ has the following form
\begin{align}
\notag\beta_1&=B_{1,n}^{s-n}\beta B_{1,n}^{n-s}=B_{1,n}^{s-n}\rho_{k_s-1}\alpha\rho_{k_s-1} B_{1,n}^{n-s}=\rho_{k_n-1}B_{1,n}^{s-n}\alpha B_{1,n}^{n-s}\rho_{k_n-1}\\
\notag&=\rho_{k_n-1}\gamma w_1\dots w_{n-1}\rho_{n-1}B_{k_n,n-1}\rho_{k_n-1}\\
\notag&=\gamma^{\rho_{k_n-1}} (w_1\dots w_{n-2})^{\rho_{k_n-1}}w_{n-1}^{\rho_{k_n-1}}\rho_{n-1}^{\rho_{k_n-1}}B_{k_n,n-1}^{\rho_{k_n-1}}
\end{align}
Since $\rho_{n-1}^{\rho_{k_n-1}}B_{k_n,n-1}^{\rho_{k_n-1}}=\rho_{k_n-1}B_{k_n-1,n}$, then we have
\begin{align}
\notag\beta_1&=\gamma^{\rho_{k_n-1}} (w_1\dots w_{n-2})^{\rho_{k_n-1}}w_{n-1}^{\rho_{k_n-1}}\rho_{k_n-1}B_{n,k_n-1}\\
\notag&=\gamma^{\rho_{k_n-1}}\rho_{k_n-1} (w_1\dots w_{n-2})w_{n-1}B_{n,k_n-1}\\
\notag&=\gamma^{\rho_{k_n-1}}\rho_{k_n-1} (w_1\dots w_{n-2})w_{n-1}\rho_{n-1}B_{n-1,k_n-1}
\end{align}

Since $k_n\leq n-1$, then $k_n-1\leq n-2$, therefore the braid
$\gamma^{\rho_{k_n-1}}\rho_{k_n-1}$ does not involve the strand $n$.
\item[Step 4.] $\beta_2=\gamma^{\rho_{k_n-1}}\rho_{k_n-1} (w_1\dots w_{n-2})\rho_{n-1}B_{n-1,k_n-1}$
\item[Step 5.] The braid $\varrho(\beta)$ has the following form.
\begin{align}
\notag\varrho(\beta)&=\gamma^{\rho_{k_n-1}}\rho_{k_n-1} (w_1\dots w_{n-2})^{\rho_{n-1}}B_{n-1,k_n-1}\\
\notag&=\gamma^{\rho_{k_n-1}} (w_1\dots w_{n-2})^{\rho_{n-1}\rho_{k_n-1}}\rho_{k_n-1}B_{n-1,k_n-1}\\
\notag&=\gamma^{\rho_{k_n-1}} (w_1\dots w_{n-2})^{\rho_{n-1}\rho_{k_n-1}}B_{n-1,k_n}^{\rho_{k_n-1}}=\varrho(\alpha)^{\rho_{k_n-1}}
\end{align}
\end{itemize}
Therefore $\varrho(\beta)=\varrho(\alpha)^{\rho_{k_n-1}}$ and by
the induction hypothesis the braids $\overline{\varrho^*(\alpha)}$ and $\overline{\varrho^*(\beta)}$ are conjugated by the element from $S_n$.

\underline{Case 3.1.3.}  $i=k_s$.

\underline{Case 3.1.3.1.} $k_s\leq s-2$. Let us find the image of the braid $\beta=\rho_{k_s}\alpha\rho_{k_s}$ under the map $\varrho$.
\begin{itemize}
\item[Step 1.] Since $k_s\leq s-2$, then and $\rho_{k_s}$ fixes $s$ and the image of $s$ under the permutation $\iota(\beta)$ is equal to $k_s+1$.
     $$ \begin{CD}
s @>~\rho_{k_s}~>> s @>~\alpha~>> k_s@>~\rho_{k_s}~>> k_s+1
\end{CD}$$
Therefore $s$ is a maximal number which is not fixed by $\beta$.
\item[Step 2.]  We have $\beta_1=B_{1,n}^{s-n}\beta B_{1,n}^{n-s}$ and $\iota(\beta_1)(n)=k_n+1$.
 $$ \begin{CD}
n @>~B_{1,n}^{s-n}~>> s @>~\beta~>> k_s+1@>~B_{1,n}^{n-s}~>> k_n+1,
\end{CD}$$
\item[Step 3.] Analogically to the case 3.1.1, since $k_s\leq s-2$, then
$\rho_{k_s}B_{1,n}^{n-s}=B_{1,n}^{n-s}\rho_{k_n}$. Therefore the braid $\beta_1$ has the following form.
\begin{align}
\notag\beta_1&=B_{1,n}^{s-n}\beta B_{1,n}^{n-s}=B_{1,n}^{s-n}\rho_{k_s}\alpha\rho_{k_s} B_{1,n}^{n-s}=\rho_{k_n}B_{1,n}^{s-n}\alpha B_{1,n}^{n-s}\rho_{k_n}\\
\notag&=\rho_{k_n}\gamma w_1\dots w_{n-1}\rho_{n-1}B_{k_n,n-1}\rho_{k_n}\\
\notag&=\gamma^{\rho_{k_n}} (w_1\dots w_{n-2})^{\rho_{k_n}}w_{n-1}^{\rho_{k_n}}\rho_{n-1}^{\rho_{k_n}}B_{k_n,n-1}^{\rho_{k_n}}
\end{align}
Since $\rho_{n-1}^{\rho_{k_n}}B_{k_n,n-1}^{\rho_{k_n}}\rho_{k_n}B_{k_n+1,n}$, then we have
\begin{align}
\notag\beta_1&=\gamma^{\rho_{k_n}} (w_1\dots w_{n-2})^{\rho_{k_n}}w_{n-1}^{\rho_{k_n}}\rho_{k_n}B_{k_n+1,n}\\
\notag&=\gamma^{\rho_{k_n}}\rho_{k_n} (w_1\dots w_{n-2})w_{n-1} B_{k_n+1,n}.
\end{align}
Since $k_s\leq s-2$, then $k_n\leq n-2$ and the braid
$\gamma^{\rho_{k_n}}\rho_{k_n}$ does not involve the strand $n$.
\item[Step 4.] $\beta_2=\gamma^{\rho_{k_n}}\rho_{k_n} (w_1\dots w_{n-2})\rho_{n-1} B_{k_n+1,n-1}$

\item[Step 5.] The braid $\varrho(\beta)$ has the form
\begin{align}
\notag\varrho(\beta)&=\gamma^{\rho_{k_n}}\rho_{k_n} (w_1\dots w_{n-2})^{\rho_{n-1}} B_{k_n+1,n-1}\\
\notag&=\gamma^{\rho_{k_n}} (w_1\dots w_{n-2})^{\rho_{n-1}\rho_{k_n}} \rho_{k_n}B_{k_n+1,n-1}\\
\notag&=\gamma^{\rho_{k_n}} (w_1\dots w_{n-2})^{\rho_{n-1}\rho_{k_n}} B_{k_n,n-1}^{\rho_{k_n}}=\varrho(\alpha)^{\rho_{k_n}}
\end{align}
\end{itemize}
    Therefore $\varrho(\beta)=\varrho(\alpha)^{\rho_{k_n}}$ and by
the induction hypothesis the braids $\overline{\varrho^*(\alpha)}$ and $\overline{\varrho^*(\beta)}$ are conjugated by the element from $S_n$.

\underline{Case 3.1.3.2.} $k_s=s-1$.

\underline{Case 3.1.3.2.1.} $\iota(\alpha)(s-1)=s$. From the formulas (\ref{a1}) and (\ref{a2}) we have
$$\alpha_1=\gamma w_1\dots w_{n-1}\rho_{n-1}B_{k_n,n-1},$$
$$\varrho(\alpha)=\gamma (w_1\dots w_{n-2})^{\rho_{n-1}}B_{k_n,n-1}.$$

Let us rewrite this equalities in more details. Since $\iota(\alpha)(s)=s-1$, then $\iota(\alpha_1)(n)=n-1$ and $B_{k_n,n}=\rho_{n-1}$.
 $$ \begin{CD}
n @>~B_{s,n}^{s-n}~>> s @>~\alpha~>> s-1@>~B_{s,n}^{n-s}~>> n-1
\end{CD}$$
 Since $\iota(\alpha)(s-1)=s$ then $\iota(\alpha_1)$ maps $(n-1)$ to $n$.
$$ \begin{CD}
n-1 @>~B_{s,n}^{s-n}~>> s-1 @>~\alpha~>> s@>~B_{s,n}^{n-s}~>> n
\end{CD}$$
Therefore, since $\gamma=\alpha_1\left(B_{k_n,n}w_1\dots w_{n-1}\right)^{-1}$, then  $\iota(\gamma)(n-1)=n-1$
$$ \begin{CD}
n-1 @>~\alpha_1~>> n @>~\left(B_{k_n,n}w_1\dots w_{n-1}\right)^{-1}~>> n-1
\end{CD}$$
and by Lemma \ref{lem2} we have
$$\gamma=\eta v_1\dots v_{n-2},$$
where $\eta$ does not involve the strand $n-1$ and $v_j\in\langle\lambda_{j,n-1},\lambda_{n-1,j}\rangle$ for $j=1,\dots,n-2$.
Therefore
$$\alpha_1=\eta v_1\dots v_{n-2} w_1\dots w_{n-1}\rho_{n-1},$$
$$\varrho(\alpha)=\eta v_1\dots v_{n-2} (w_1\dots w_{n-2})^{\rho_{n-1}}.$$

Let us find the braid $\varrho(\beta)$.
\begin{itemize}
\item[Step 1.] Since $k_s=s-1$ and $\iota(\alpha)(s-1)=s$, then the image of $s$ under the permutation $\iota(\beta)=\iota(\rho_{k_s}\alpha\rho_{k_s})=\iota(\rho_{s-1}\alpha\rho_{s-1})$ is equal to $s-1$.
     $$ \begin{CD}
s @>~\rho_{s-1}~>> s-1 @>~\alpha~>> s@>~\rho_{s-1}a~>> s-1,
\end{CD}$$
then $s$ is a maximal number  which is not fixed by $\beta$.
\item[Step 2.] We have $\beta_1=B_{1,n}^{s-n}\beta B_{1,n}^{n-s}$ and $\iota(\beta_1)(n)=n-1$.
    $$ \begin{CD}
n @>~B_{1,n}^{s-n}~>> s @>~\beta~>> s-1@>~B_{1,n}^{n-s}~>> n-1,
\end{CD}$$
\item[Step 3.] Since $\rho_{s-1}B_{1,n}^{n-s}=B_{1,n}^{n-s}\rho_{n-1}$, then we have
\begin{align}
\notag\beta_1&=B_{1,n}^{s-n}\beta B_{1,n}^{n-s}=B_{1,n}^{s-n}\rho_{s-1}\alpha\rho_{s-1} B_{1,n}^{n-s}=\rho_{n-1}B_{1,n}^{s-n}\alpha B_{1,n}^{n-s}\rho_{n-1}\\
\notag&=\rho_{n-1}\eta v_1\dots v_{n-2} w_1\dots w_{n-1}\rho_{n-1}\rho_{n-1}\\
\notag&=\eta^{\rho_{n-1}} v_1^{\rho_{n-1}}\dots v_{n-2}^{\rho_{n-1}} w_1^{\rho_{n-1}}\dots w_{n-2}^{\rho_{n-1}}w_{n-1}^{\rho_{n-1}}\rho_{n-1}
\end{align}
Since the braid $\eta$ does not involve the strands $n-1,n$, then $\eta^{\rho_{n-1}}=\eta$. Also by Lemma \ref{lem1} for $j=1,\dots, n-2$ we have
$$v_j^{\rho_{n-1}}\in\langle\lambda_{j,n},\lambda_{n,j}\rangle,~~~~~
w_j^{\rho_{n-1}}\in\langle\lambda_{j,n-1},\lambda_{n-1,j}\rangle$$
and $w_{n-1}^{\rho_{n-1}}\in\langle\lambda_{n,n-1},\lambda_{n-1,n}\rangle$. Therefore
$$\beta_1=\eta w_1^{\rho_{n-1}}\dots w_{n-2}^{\rho_{n-1}} v_1^{\rho_{n-1}}\dots v_{n-2}^{\rho_{n-1}} w_{n-1}^{\rho_{n-1}}\rho_{n-1},$$
where the braid
$\eta w_1^{\rho_{n-1}}\dots w_{n-2}^{\rho_{n-1}}$ does not involve the strand $n$, the braid
$v_j^{\rho_{n-1}}$ belongs to $\langle\lambda_{j,n},\lambda_{n,j}\rangle$ for $j=1,\dots,n-2$ and
$w_{n-1}^{\rho_{n-1}}\in\langle\lambda_{n,n-1},\lambda_{n-1,n}\rangle$.

\item[Step 4.] $\beta_2=\eta w_1^{\rho_{n-1}}\dots w_{n-2}^{\rho_{n-1}} v_1^{\rho_{n-1}}\dots v_{n-2}^{\rho_{n-1}} \rho_{n-1}$

\item[Step 5.] The braid $\varrho(\beta)$ has the following form
\begin{align}
\notag\varrho(\beta)&=\eta w_1^{\rho_{n-1}}\dots w_{n-2}^{\rho_{n-1}} v_1\dots v_{n-2}\\
\notag&=\eta v_1\dots v_{n-2}w_1^{\rho_{n-1}}\dots w_{n-2}^{\rho_{n-1}}[w_1^{\rho_{n-1}}\dots w_{n-2}^{\rho_{n-1}},v_1\dots v_{n-2}]
\end{align}
\end{itemize}
Therefore by the remark \ref{rmk2} the braids $\varrho^*(\beta)$ and $\varrho^*(\alpha)$ are equal modulo $UVP_m^{\prime}$, i.~e. $\overline{\varrho^*(\alpha)}=\overline{\varrho^*(\beta)}$ and these braids are certainly conjugated.

\underline{Case 3.1.3.2.2.} $\iota(\alpha)(s-1)=k_{s-1}\leq s-2$. By the equalities (\ref{a1}) and (\ref{a2}) we have
$$\alpha_1=\gamma w_1\dots w_{n-1}\rho_{n-1}B_{k_n,n-1},$$
$$\varrho(\alpha)=\gamma (w_1\dots w_{n-2})^{\rho_{n-1}}B_{k_n,n-1}.$$
Let us rewrite these equalities in more details.
Since $\iota(\alpha_1)(n)=k_n=n-1$, then $B_{k_n,n-1}=1$. Also we have $\iota(\alpha_1)(n-1)=n-s+k_{s-1}=k_{n-1}$
  $$ \begin{CD}
n-1 @>~B_{1,n}^{s-n}~>> s-1 @>~\alpha~>> k_{s-1}@>~B_{1,n}^{n-s}~>> n-s+k_{s-1}=k_{n-1},
\end{CD}$$
 Since $\gamma=\alpha_1\left(B_{k_n,n}w_1\dots w_{n-1}\right)^{-1}$, then $\iota(\gamma)$ maps $n-1$ to $k_{n-1}$
$$ \begin{CD}
n-1 @>~\alpha_1~>> k_{n-1} @>~\left(B_{k_n,n}w_1\dots w_{n-1}\right)^{-1}~>> k_{n-1},
\end{CD}$$
  therefore  by Lemma \ref{lem2} we have
$$\gamma=\eta v_1\dots v_{n-2}B_{k_{n-1},n-1},$$
where $\eta$ does not involve the strand $n-1,n$ and $v_j$ belongs to $\langle\lambda_{j,n-1},\lambda_{n-1,j}\rangle$ for $j=1,\dots,n-2$.
Therefore we have
$$\alpha_1=\eta v_1\dots v_{n-2}B_{k_{n-1},n-1} w_1\dots w_{n-1}\rho_{n-1},$$
$$\varrho(\alpha)=\eta v_1\dots v_{n-2}B_{k_{n-1},n-1} (w_1\dots w_{n-2})^{\rho_{n-1}}.$$
Let us count $\varrho^2 (\alpha)$.
\begin{itemize}
\item[Step 1.] Since $\eta$ does not involve the strand $n-1$, then $\iota(\varrho(\alpha))(n-1)=k_{n-1}$,
then $n-1$ is a maximal number  which is not fixed by $\varrho(\alpha)$.
\item[Step 2.]  $\varrho(\alpha)_1=\varrho(\alpha)=\eta v_1\dots v_{n-2}B_{k_{n-1},n-1} (w_1\dots w_{n-2})^{\rho_{n-1}}$.\\
\item[Step 3.] Let us rewrite the braid $\varrho(\alpha)_1$ in the required form.
\begin{align}
\notag\varrho(\alpha)_1&=\eta v_1\dots v_{n-2}B_{k_{n-1},n-1} (w_1\dots w_{n-2})^{\rho_{n-1}}\\
\notag&=\eta v_1\dots v_{n-2} (w_1\dots w_{n-2})^{\rho_{n-1}B_{k_{n-1},n-1}^{-1}}B_{k_{n-1},n-1}\\
\notag&=\eta v_1\dots v_{n-2} (w_1\dots w_{k_{n-1}-1}w_{k_{n-1}}w_{k_{n-1}+1}\dots w_{n-2})^{\rho_{n-1}B_{k_{n-1},n-1}^{-1}}B_{k_{n-1},n-1}\\
\notag&=\eta (w_1\dots w_{k_{n-1}-1}w_{k_{n-1}+1}\dots w_{n-2})^{\rho_{n-1}B_{k_{n-1},n-1}^{-1}}\\
\notag &~\cdot~ v_1\dots v_{n-3} v_{n-2}w_{k_{n-1}}^{\rho_{n-1}B_{k_{n-1},n-1}^{-1}}B_{k_{n-1},n-1}
\end{align}
By Theorem \ref{t2} the braid
$\eta (w_1\dots w_{k_{n-1}-1}w_{k_{n-1}+1}\dots w_{n-2})^{\rho_{n-1}B_{k_{n-1},n-1}^{-1}}$ does not involve the strand $n-1$, the braid $v_j\in\langle\lambda_{j,n-1},\lambda_{n-1,j}\rangle$ for $j=1,\dots,n-3$ and $v_{n-2}w_{k_{n-1}}^{\rho_{n-1}B_{k_{n-1},n-1}^{-1}}\in\langle\lambda_{n-1,n-2},\lambda_{n-2,n-1}\rangle$.
\item[Step 4.] $\varrho(\alpha)_2=\eta (w_1\dots w_{k_{n-1}-1}w_{k_{n-1}+1}\dots w_{n-2})^{\rho_{n-1}B_{k_{n-1},n-1}^{-1}} v_1\dots v_{n-3} \rho_{n-2}B_{k_{n-1},n-2}$

\item[Step 5.] The braid $\varrho^2(\alpha)$ follows
\begin{align}\notag\varrho(\alpha)_2&=\eta (w_1\dots w_{k_{n-1}-1}w_{k_{n-1}+1}\dots w_{n-2})^{\rho_{n-1}B_{k_{n-1},n-1}^{-1}} (v_1\dots v_{n-3})^{\rho_{n-2}} B_{k_{n-1},n-2}
\end{align}
\end{itemize}

Let us count $\varrho(\beta)$
\begin{itemize}
\item[Step 1.] Since $k_s=s-1$ and $\iota(\alpha)(s-1)=k_{s-1}\leq s-2$, then the image of $s$ under the permutation $\beta=\rho_{k_s}\alpha\rho_{k_s}=\rho_{s-1}\alpha\rho_{s-1}$ is equal to $k_{s-1}$.
     $$ \begin{CD}
s @>~\rho_{s-1}~>> s-1 @>~\alpha~>> k_{s-1}@>~\rho_{s-1}~>> k_{s-1},
\end{CD}$$
then $s$ is a maximal number  which is not fixed by $\beta$.
The image of $s-1$ under the permutation $\beta=\rho_{k_s}\alpha\rho_{k_s}=\rho_{s-1}\alpha\rho_{s-1}$ is equal to $s$.
     $$ \begin{CD}
s-1 @>~\rho_{s-1}~>> s @>~\alpha~>> s-1@>~\rho_{s-1}~>> s,
\end{CD}$$
\item[Step 2.] We have $\beta_1=B_{n,1}^{s-n}\beta B_{n,1}^{n-s}$ and
 $\iota(\beta_1)(n)=k_{s-1}+n-s=k_{n-1}$, $\iota(\beta_1)(n-1)=n$
 $$ \begin{CD}
n @>~B_{n,1}^{s-n}~>> s @>~\beta~>> k_{s-1}@>~B_{1,n}^{n-s}>> k_{s-1}+n-s~=k_{n-1},
\end{CD}$$
$$ \begin{CD}
n-1 @>~B_{n,1}^{s-n}~>> s-1 @>~\beta~>> s@>~B_{n,1}^{s-n}~>> n,
\end{CD}$$
\item[Step 3.]Since $k_s=s-1< s$ then $\rho_{s-1}B_{1,n}^{n-s}=B_{1,n}^{n-s}\rho_{n-1}$, therefore we have
\begin{align}
\notag\beta_1&=B_{1,n}^{s-n}\beta B_{1,n}^{n-s}=B_{1,n}^{s-n}\rho_{s-1}\alpha\rho_{s-1}B_{1,n}^{n-s}= \rho_{n-1}B_{1,n}^{s-n}\alpha B_{1,n}^{n-s}\rho_{n-1}\\
\notag&=\rho_{n-1}\eta v_1\dots v_{n-2}B_{k_{n-1},n-1} w_1\dots w_{n-1}\rho_{n-1}\rho_{n-1}\\
\notag&=\eta^{\rho_{n-1}} (v_1\dots v_{n-2})^{\rho_{n-1}}B_{k_{n-1},n-1}^{\rho_{n-1}} (w_1\dots w_{n-1})^{\rho_{n-1}}\rho_{n-1}\\
\notag&=\eta^{\rho_{n-1}} (v_1\dots v_{n-2})^{\rho_{n-1}} (w_1\dots w_{n-1})^{\rho_{n-1}\left(B_{k_{n-1},n-1}^{\rho_{n-1}}\right)^{-1}}B_{k_{n-1},n-1}^{\rho_{n-1}}\rho_{n-1}\\
\notag&=\eta^{\rho_{n-1}} (v_1\dots v_{n-2})^{\rho_{n-1}} (w_1\dots w_{n-1})^{B_{k_{n-1},n-1}^{-1}\rho_{n-1}}B_{k_{n-1},n}\\
\notag&=\eta^{\rho_{n-1}} (v_1\dots v_{n-2})^{\rho_{n-1}} (w_1\dots w_{n-1})^{B_{k_{n-1},n}^{-1}}B_{k_{n-1},n}\\
\notag&=\eta^{\rho_{n-1}} (v_1\dots v_{n-2})^{\rho_{n-1}} (w_1\dots w_{k_{n-1}-1} w_{k_{n-1}}w_{k_{n-1}+1}\dots w_{n-1})^{B_{k_{n-1},n}^{-1}}B_{k_{n-1},n}\\
\notag&=\eta^{\rho_{n-1}} (w_1\dots w_{k_{n-1}-1} w_{k_{n-1}+1}\dots w_{n-1})^{B_{k_{n-1},n}^{-1}}(v_1\dots v_{n-2})^{\rho_{n-1}}w_{k_{n-1}}^{B_{k_{n-1},n}^{-1}}B_{k_{n-1},n}
\end{align}
By Theorem \ref{t2} the braid
$\eta^{\rho_{n-1}} (w_1\dots w_{k_{n-1}-1} w_{k_{n-1}+1}\dots w_{n-1})^{B_{k_{n-1},n}^{-1}}$ does not involve the strand $n$,
$v_j^{\rho_{n-1}}\in\langle\lambda_{j,n},\lambda_{n,j}\rangle$ for $j=1,\dots n-2$ and
$w_{k_{n-1}}^{B_{k_{n-1},n}^{-1}}$ belongs to $\langle\lambda_{n-1,n},\lambda_{n,n-1}\rangle$.

\item[Step 4.] $\beta_2=\eta^{\rho_{n-1}} (w_1\dots w_{k_{n-1}-1} w_{k_{n-1}+1}\dots w_{n-1})^{B_{k_{n-1},n}^{-1}}(v_1\dots v_{n-2})^{\rho_{n-1}}\rho_{n-1}B_{k_{n-1},n-1}$

\item[Step 5.] The braid $\varrho(\beta)$ has the following form
\begin{equation}
\label{eq1112}
\varrho(\beta)=\eta^{\rho_{n-1}} (w_1\dots w_{k_{n-1}-1} w_{k_{n-1}+1}\dots w_{n-1})^{B_{k_{n-1},n}^{-1}}(v_1\dots v_{n-2})B_{k_{n-1},n-1}
\end{equation}
\end{itemize}
Let us count $\varrho^2(\beta)$
\begin{itemize}
\item[Step 1.] From the equality (\ref{eq1112}) the image of $n-1$ under the permutation $\iota(\varrho(\beta))$ is equal to $k_{n-1}$,
therefore $n-1$ is a maximal number  which is not fixed by $\varrho(\beta)$.
\item[Step 2.]  $\varrho(\beta)_1=\varrho(\beta$).
\item[Step 3.] Rewrite the braid $\varrho(\beta)_1$ in the required form
\begin{align}
\notag\varrho(\beta)_1&=\eta^{\rho_{n-1}} (w_1\dots w_{k_{n-1}-1} w_{k_{n-1}+1}\dots w_{n-1})^{B_{k_{n-1},n}^{-1}}(v_1\dots v_{n-2})B_{k_{n-1},n-1}\\
\notag&=\eta  (w_1\dots w_{k_{n-1}-1} w_{k_{n-1}+1}\dots w_{n-2})^{B_{k_{n-1},n}^{-1}}v_1\dots v_{n-3}w_{n-1}^{B_{k_{n-1},n}^{-1}}v_{n-2}B_{k_{n-1},n-1}
\end{align}
Since $\eta$ does not involve the strands $n-1,n$, then $\eta^{\rho_{n-1}}=\eta$ does not involve the strand $n-1$. By Lemma \ref{lem1} the braid $v_{n-2}w_{n-1}^{B_{k_{n-1},n}^{-1}}$ belongs to $\langle\lambda_{n-1,n-2},\lambda_{n-2,n-1}\rangle$, the braid
$(w_1\dots w_{k_{n-1}-1} w_{k_{n-1}+1}\dots w_{n-2})^{B_{k_{n-1},n}^{-1}}v_1\dots v_{n-3}$ belongs to $\langle\lambda_{1,n-1},\lambda_{n-1,1}\rangle\times\dots\times \langle\lambda_{n-1,n-2},\lambda_{n-2,n-1}\rangle$.

\item[Step 4.] $\varrho(\beta)_2=\eta  (w_1\dots w_{k_{n-1}-1} w_{k_{n-1}+1}\dots w_{n-2})^{B_{k_{n-1},n}^{-1}}v_1\dots v_{n-3}\rho_{n-2}B_{k_{n-1},n-2}$

\item[Step 5.]
$\varrho^2({\beta})=\eta  (w_1\dots w_{k_{n-1}-1} w_{k_{n-1}+1}\dots w_{n-2})^{B_{k_{n-1},n}^{-1}\rho_{n-2}}(v_1\dots v_{n-3})^{\rho_{n-2}}B_{k_{n-1},n-2}$
\end{itemize}
Using simple calculations in symmetric group it is easy to show that $B_{k_{n-1},n}^{-1}\rho_{n-2}=\rho_{n-1}B_{k_{n-1},n-1}^{-1}\rho_{n-1}$, therefore
\begin{multline*}(w_1\dots w_{k_{n-1}-1} w_{k_{n-1}+1}\dots w_{n-2})^{B_{k_{n-1},n}^{-1}\rho_{n-2}}=\\
=\left((w_1\dots w_{k_{n-1}-1} w_{k_{n-1}+1}\dots w_{n-2})^{\rho_{n-1}B_{k_{n-1},n-1}^{-1}}\right)^{\rho_{n-1}}
\end{multline*}
and since $(w_1\dots w_{k_{n-1}-1} w_{k_{n-1}+1}\dots w_{n-2})^{\rho_{n-1}B_{k_{n-1},n-1}^{-1}}$ does not involve the strands $n-1$, $n$, then
\begin{multline*}\left((w_1\dots w_{k_{n-1}-1} w_{k_{n-1}+1}\dots w_{n-2})^{\rho_{n-1}B_{k_{n-1},n-1}^{-1}}\right)^{\rho_{n-1}}=\\
=(w_1\dots w_{k_{n-1}-1} w_{k_{n-1}+1}\dots w_{n-2})^{\rho_{n-1}B_{k_{n-1},n-1}^{-1}}
\end{multline*}

Therefore $\varrho^2(\alpha)=\varrho^2(\beta)$ and by the induction hypothesis the braids $\overline{\varrho^*(\alpha)}$ and $\overline{\varrho^*(\beta)}$ are conjugated by the element from $S_n$.

\underline{Case 3.1.4.}  $k_s+1\leq i \leq s-2$. Let us count the braid $\varrho(\beta)=\varrho(\rho_{i}\alpha\rho_{i})$.
\begin{itemize}
\item[Step 1.] Since $k_s\leq s-3$ and $k_s+1\leq i \leq s-2$, then the braid $\rho_i$ fixes $s$ and $k_s$, and therefore the image of $s$ under the permutation $\iota(\beta)=\iota(\rho_{i}\alpha\rho_{i})$ is equal to $k_s$.
$$ \begin{CD}
s @>~\rho_{i}~>> s @>~\alpha~>> k_s@>~\rho_{i}~>> k_s
\end{CD}$$
Therefore  $s$ is a maximal number  which is not fixed by $\beta$.
\item[Step 2.]  We have $\beta_1=B_{1,n}^{s-n}\beta B_{1,n}^{n-s}$ and $\iota(\beta_1)(n)=k_n$.
 $$ \begin{CD}
n @>~B_{1,n}^{s-n}~>> s @>~\beta~>> k_s@>~B_{1,n}^{n-s}~>> k_n,
\end{CD}$$
\item[Step 3.] Since $i\leq s-2$, then $\rho_{i}B_{1,n}^{n-s}=B_{1,n}^{n-s}\rho_{n-s+i},$ and therefore the braid $\beta_1$ follows.
\begin{align}
\notag\beta_1&=B_{1,n}^{s-n}\beta B_{1,n}^{n-s}=B_{1,n}^{s-n}\rho_{i}\alpha\rho_{i} B_{1,n}^{n-s}=\rho_{i+n-s}B_{1,n}^{s-n}\alpha B_{1,n}^{n-s}\rho_{i+n-s}\\
\notag&=\rho_{i+n-s}\gamma w_1\dots w_{n-1}\rho_{n-1}B_{k_n,n-1}\rho_{i+n-s}\\
\notag&=\gamma^{\rho_{i+n-s}} (w_1\dots w_{n-2})^{\rho_{i+n-s}}w_{n-1}^{\rho_{i+n-s}}\rho_{n-1}^{\rho_{i+n-s}}B_{k_n,n-1}^{\rho_{i+n-s}}
\end{align}
Using simple calculation in the permutation group it is easy to show that
$\rho_{n-1}^{\rho_{i+n-s}}B_{k_n,n-1}^{\rho_{i+n-s}}=B_{k_n,n}^{\rho_{i+n-s}}=\rho_{i+n-s}\rho_{i+n-s-1}B_{k_n,n}$, hence
\begin{align}
\notag\beta_1&=\gamma^{\rho_{i+n-s}} (w_1\dots w_{n-2})^{\rho_{i+n-s}}w_{n-1}^{\rho_{i+n-s}}\rho_{i+n-s}\rho_{i+n-s-1}B_{k_n,n}\\
\notag&=\gamma^{\rho_{i+n-s}}\rho_{i+n-s}\rho_{i+n-s-1} (w_1\dots w_{n-2})^{\rho_{i+n-s-1}}w_{n-1}^{\rho_{i+n-s-1}}B_{k_n,n}
\end{align}
Since $i\leq s-2$, then $i+n-s-1\leq n-3$ and the braid
$\gamma^{\rho_{i+n-s}}\rho_{i+n-s}\rho_{i+n-s-1}$ does not involve the strand $n$. Also by Lemma \ref{lem1}
$(w_1\dots w_{n-2})^{\rho_{i+n-s-1}}$ belongs to $\langle\lambda_{1,n},\lambda_{n,1}\rangle\times\dots\times\langle\lambda_{n-2,n},\lambda_{n,n-2}\rangle$ and
$w_{n-1}^{\rho_{i+n-s-1}}=w_{n-1}\in\langle\lambda_{n-1,n},\lambda_{n,n-1}\rangle$.

\item[Step 4.] $\beta_2=\gamma^{\rho_{i+n-s}}\rho_{i+n-s}\rho_{i+n-s-1} (w_1\dots w_{n-2})^{\rho_{i+n-s-1}}\rho_{n-1}B_{k_n,n-1}$

\item[Step 5.] The braid $\beta$ has the following form
\begin{align}
\notag\varrho(\beta)&=\gamma^{\rho_{i+n-s}}\rho_{i+n-s}\rho_{i+n-s-1} (w_1\dots w_{n-2})^{\rho_{i+n-s-1}\rho_{n-1}}B_{k_n,n-1}\\
\notag&=\gamma^{\rho_{i+n-s}} (w_1\dots w_{n-2})^{\rho_{i+n-s-1}\rho_{n-1}\rho_{i+n-s-1}\rho_{i+n-s}}\rho_{i+n-s}\rho_{i+n-s-1}B_{k_n,n-1}\\
\notag&=\gamma^{\rho_{i+n-s}} (w_1\dots w_{n-2})^{\rho_{n-1}\rho_{i+n-s}}B_{k_n,n-1}^{\rho_{i+n-s}}=\varrho(\alpha)^{\rho_{i+n-s}}
\end{align}
\end{itemize}
    Therefore $\varrho(\beta)=\varrho(\alpha)^{\rho_{i+n-s}}$ and by the induction hypothesis the braids $\overline{\varrho^*(\alpha)}$ and $\overline{\varrho^*(\beta)}$ are conjugated by the element from $S_n$.

\underline{Case 3.1.5.}  $i=s-1$. In this case we can consider that $k_s< s-1$ since the case when $k_s=s-1=i$ is already solved in the case 3.1.3.2.

\underline{Case 3.1.5.1.}  $\iota(\alpha)(s-1)=s$. By the equalities (\ref{a1}) and (\ref{a2}) the braids $\alpha_1$ and $\varrho(\alpha)$ have the following forms
$$\alpha_1=\gamma w_1\dots w_{n-1}\rho_{n-1}B_{k_n,n-1},$$
$$\varrho(\alpha)=\gamma (w_1\dots w_{n-2})^{\rho_{n-1}}B_{k_n,n-1}.$$
Let us rewrite this equalities in more details. Since $\iota(\alpha)(s-1)=s$, then the permutation $\iota(\alpha_1)$ maps $n-1$ to $n$.
  $$ \begin{CD}
n-1 @>~B_{1,n}^{s-n}~>> s-1 @>~\alpha~>> s@>~B_{1,n^{n-s}}~>> n
\end{CD}$$
 Therefore the permutation $\iota(\gamma)=\iota(\alpha_1\left(w_1\dots w_{n-1}\rho_{n-1}B_{k_n,n-1}\right)^{-1})$ fixes $n-1$
   $$ \begin{CD}
n-1 @>~\alpha_1~>> n @>~\left(w_1\dots w_{n-1}\rho_{n-1}B_{k_n,n-1}\right)^{-1}~>> n-1,
\end{CD}$$
 and by Lemma \ref{lem2} we have
$$\gamma=\eta v_1\dots v_{n-2},$$
where $\eta$ does not involve the strands $n-1,n$ and $v_j\in\langle\lambda_{j,n-1},\lambda_{n-1,j}\rangle$.
Therefore the braids $\alpha_1$, $\varrho(\alpha)$ can be rewritten.
$$\alpha_1=\eta v_1\dots v_{n-2} w_1\dots w_{n-1}\rho_{n-1}B_{k_{n},n-1}$$
$$\varrho(\alpha)=\eta v_1\dots v_{n-2} (w_1\dots w_{n-2})^{\rho_{n-1}}B_{k_{n},n-1}$$
Let us count $\varrho^2(\alpha)$
\begin{itemize}
\item[Step 1.] Since $k_s<s-1$, then $k_n<n-1$, therefore $B_{k_1,n-1}\neq1$,
then $n-1$ is a maximal number  which is not fixed by $\varrho(\alpha)$.
\item[Step 2.]  $\varrho(\alpha)_1=\varrho(\alpha)$.
\item[Step 3.]The braid $\varrho(\alpha)_1$ has the following form.
\begin{align}
\notag\varrho(\alpha)_1&=\varrho(\alpha)=\eta v_1\dots v_{n-2} (w_1\dots w_{n-2})^{\rho_{n-1}}B_{k_{n},n-1}\\
\notag&=\eta v_1\dots v_{n-3} (w_1\dots w_{n-3})^{\rho_{n-1}}v_{n-2}w_{n-2}^{\rho_{n-1}}B_{k_{n},n-1}
\end{align}
Here the braid $\eta$ does not involve the strand $n-1$, by Lemma \ref{lem1} the braid
$v_1\dots v_{n-3} (w_1\dots w_{n-3})^{\rho_{n-1}}$ belongs to $\langle\lambda_{1,n-1},\lambda_{n-1,1}\rangle\times\dots\times\langle\lambda_{1,n-3},\lambda_{n-3,1}\rangle$ and
$v_{n-2}w_{n-2}^{\rho_{n-1}}$ belongs to $\langle\lambda_{n-1,n-2},\lambda_{n-2,n-1}\rangle$.
\item[Step 4.] $\varrho(\alpha)_2=\eta v_1\dots v_{n-3} (w_1\dots w_{n-3})^{\rho_{n-1}}{\rho_{n-2}}B_{k_{n},n-2}$
\item[Step 5.] The braid $\varrho^2(\alpha)$ follows.
\begin{equation}
\label{eq12}
\varrho^{2}({\alpha})=\eta (v_1\dots v_{n-3})^{\rho_{n-2}} (w_1\dots w_{n-3})^{\rho_{n-1}\rho_{n-2}}B_{k_{n},n-2}
\end{equation}
\end{itemize}
Let us find the braid $\varrho(\beta)=\varrho(\rho_{i}\alpha\rho_{i})=\varrho(\rho_{s-1}\alpha\rho_{s-1})$.
\begin{itemize}
\item[Step 1.] Since $i=s-1$ and $\iota(\alpha)(s-1)=s$, then the image of $s$ under the permutation $\beta=\rho_{s-1}\alpha\rho_{s-1}$ is equal to $s-1$
$$ \begin{CD}
s @>~\rho_{s-1}~>> s-1 @>~\alpha~>> s@>~\rho_{s-1}~>> s-1,
\end{CD}$$
then $s$ is a maximal number  which is not fixed by $\beta$.
\item[Step 2.]  We have $\beta_1=B_{1,n}^{s-n}\beta B_{1,n}^{n-s}$ and $\iota(\beta_1)(n)=n-1$.
 $$ \begin{CD}
n @>~B_{1,n}^{s-n}~>> s @>~\beta~>> s-1@>~B_{1,n}^{n-s}~>> n-1
\end{CD}$$
\item[Step 3.] Since $\rho_{s-1}B_{1,n}^{n-s}=B_{1,n}^{n-s}\rho_{n-1}$, then then we have
\begin{align}
\notag\beta_1&=B_{1,n}^{s-n}\beta B_{1,n}^{n-s}=B_{1,n}^{s-n}\rho_{s-1}\alpha\rho_{s-1} B_{1,n}^{n-s}=\rho_{n-1}B_{1,n}^{s-n}\alpha B_{1,n}^{n-s}\rho_{n-1}\\
\notag&=\rho_{n-1}\eta v_1\dots v_{n-2} w_1\dots w_{n-1}\rho_{n-1}B_{k_{n},n-1}\rho_{n-1}\\
\notag&=\eta^{\rho_{n-1}} (v_1\dots v_{n-2})^{\rho_{n-1}} (w_1\dots w_{n-1})^{\rho_{n-1}}\rho_{n-1}B_{k_{n},n-1}^{\rho_{n-1}}\\
\notag&=\eta^{\rho_{n-1}} (v_1\dots v_{n-2})^{\rho_{n-1}} (w_1\dots w_{n-1})^{\rho_{n-1}}B_{k_{n},n-1}{\rho_{n-1}}\\
\notag&=\eta^{\rho_{n-1}}B_{k_{n},n-1} (v_1\dots v_{n-2})^{\rho_{n-1}B_{k_{n},n-1}} (w_1\dots w_{n-1})^{\rho_{n-1}B_{k_{n},n-1}}{\rho_{n-1}}\\
\notag&=\eta^{\rho_{n-1}}B_{k_{n},n-1} (v_1\dots v_{n-2})^{B_{k_{n},n}} (w_1\dots w_{n-1})^{B_{k_{n},n}}{\rho_{n-1}}\\
\notag&=\eta^{\rho_{n-1}}B_{k_{n},n-1} (v_1\dots v_{n-3})^{B_{k_{n},n}} (w_1\dots w_{n-2})^{B_{k_{n},n}}w_{n-1}^{B_{k_{n},n}}v_{n-2}^{B_{k_{n},n}}{\rho_{n-1}}\\
\notag&=\eta^{\rho_{n-1}}B_{k_{n},n-1}(w_1\dots w_{n-2})^{B_{k_{n},n}} (v_1\dots v_{n-3})^{B_{k_{n},n}} w_{n-1}^{B_{k_{n},n}}v_{n-2}^{B_{k_{n},n}}\rho_{n-1}
\end{align}
Since $\eta$ does not involve the strands $n-1,n$, then $\eta^{\rho_{n-1}}=\eta$ and the braid
$\eta^{\rho_{n-1}}B_{k_{n},n-1}(w_1\dots w_{n-2})^{B_{k_{n},n}}$ does not involve the strand $n$. By Lemma \ref{lem2} the braid
$(v_1\dots v_{n-3})^{B_{k_{n},n}} w_{n-1}^{B_{k_{n},n}}$ belongs to $\langle\lambda_{1,n},\lambda_{n,1}\rangle\times\dots\times\langle\lambda_{n-2,n},\lambda_{n,n-2}\rangle$ and
$v_{n-2}^{B_{k_{n},n}}$ belongs to $\langle\lambda_{n-1,n},\lambda_{n,n-1}\rangle$.
\item[Step 4.] $\beta_2=\eta^{\rho_{n-1}}B_{k_{n},n-1}(w_1\dots w_{n-2})^{B_{k_{n},n}} (v_1\dots v_{n-3})^{B_{k_{n},n}} w_{n-1}^{B_{k_{n},n}}\rho_{n-1}$

\item[Step 5.]
$\varrho(\beta)=\eta^{\rho_{n-1}}B_{k_{n},n-1}(w_1\dots w_{n-2})^{B_{k_{n},n}} (v_1\dots v_{n-3})^{B_{k_{n},n}{\rho_{n-1}}} w_{n-1}^{B_{k_{n},n}{\rho_{n-1}}}$
\end{itemize}
    Let us count  the braid $\varrho^2(\beta)=\varrho(\varrho(\beta))$.
  \begin{itemize}
\item[Step 1.] Since $\eta$ does not involve the strand $n-1,n$ therefore $\eta^{\rho_{n-1}}=\eta$ and $\iota(\eta)$ fixes $n-1$. Moreover since $(w_1\dots w_{n-2})^{B_{k_{n},n}} (v_1\dots v_{n-3})^{B_{k_{n},n}{\rho_{n-1}}} w_{n-1}^{B_{k_{n},n}{\rho_{n-1}}}$ is a pure braid, that the image of $n-1$ under the permutation $\iota(\varrho(\beta))$ is equal to $\iota(B_{k_n,n-1})(n-1)=k_n$. Therefore
$n-1$ is a maximal number  which is not fixed by $\varrho(\beta)$.
\item[Step 2.]  $\varrho(\beta)_1=\varrho(\beta)$.
\item[Step 3.] The braid $\varrho(\beta)_1$ has the following form
\begin{align}
\notag \varrho(\beta)_1&=\eta^{\rho_{n-1}}B_{k_{n},n-1}(w_1\dots w_{n-2})^{B_{k_{n},n}} (v_1\dots v_{n-3})^{B_{k_{n},n}{\rho_{n-1}}} w_{n-1}^{B_{k_{n},n}{\rho_{n-1}}}\\
\notag&=\eta^{\rho_{n-1}}(w_1\dots w_{n-2})^{B_{k_{n},n}B_{k_{n},n-1}^{-1}} (v_1\dots v_{n-3})^{B_{k_{n},n}{\rho_{n-1}}B_{k_{n},n-1}^{-1}}\\
\notag&~\cdot~w_{n-1}^{B_{k_{n},n}{\rho_{n-1}}B_{k_{n},n-1}^{-1}}B_{k_{n},n-1}
\end{align}
 Since $B_{k_{n},n}B_{k_{n},n-1}^{-1}=\rho_{n-1}$, $B_{k_{n},n}{\rho_{n-1}}B_{k_{n},n-1}^{-1}=\rho_{n-2}\rho_{n-1}$, then we have
 \begin{align}
\notag \varrho(\beta)_1&=\eta^{\rho_{n-1}}(w_1\dots w_{n-2})^{\rho_{n-1}} (v_1\dots v_{n-3})^{\rho_{n-2}\rho_{n-1}} w_{n-1}^{\rho_{n-2}\rho_{n-1}}B_{k_{n},n-1}\\
\notag &=\eta^{\rho_{n-1}}(w_1\dots w_{n-3})^{\rho_{n-1}} (v_1\dots v_{n-3})^{\rho_{n-2}\rho_{n-1}} w_{n-2}^{\rho_{n-1}}w_{n-1}^{\rho_{n-2}\rho_{n-1}}B_{k_{n},n-1}
\end{align}
Here the braid $\eta^{\rho_{n-1}}=\eta$ does not involve the strand $n-1$. Also by Lemma \ref{lem2} the braid
$(w_1\dots w_{n-3})^{\rho_{n-1}} (v_1\dots v_{n-3})^{\rho_{n-2}\rho_{n-1}}$ belongs to the group $\langle\lambda_{1,n-1},\lambda_{n-1,1}\rangle\times\dots\times\langle\lambda_{n-3,n-1},\lambda_{n-1,n-3}\rangle$ and
$w_{n-1}^{\rho_{n-2}\rho_{n-1}}\in\langle\lambda_{n-1,n-2},\lambda_{n-2,n-1}\rangle$.
\item[Step 4.] $\varrho(\beta)_2=\eta^{\rho_{n-1}}(w_1\dots w_{n-3})^{\rho_{n-1}} (v_1\dots v_{n-3})^{\rho_{n-2}\rho_{n-1}} \rho_{n-2}B_{k_{n},n-2}$

\item[Step 5.] The braid $\varrho^2(\beta)$ has the following form
\begin{align}\notag\varrho^2(\beta)&=\eta^{\rho_{n-1}}(w_1\dots w_{n-3})^{\rho_{n-1}\rho_{n-2}} (v_1\dots v_{n-3})^{\rho_{n-2}\rho_{n-1}\rho_{n-2}} B_{k_{n},n-2}\\
\notag\varrho^2(\beta)&=\eta^{\rho_{n-1}}(v_1\dots v_{n-3})^{\rho_{n-2}\rho_{n-1}\rho_{n-2}}(w_1\dots w_{n-3})^{\rho_{n-1}\rho_{n-2}}\\
\notag&~\cdot~[(w_1\dots w_{n-3})^{\rho_{n-1}\rho_{n-2}},(v_1\dots v_{n-3})^{\rho_{n-2}\rho_{n-1}\rho_{n-2}}] B_{k_{n},n-2}
\end{align}

By Lemma \ref{lem2} it is obvious that $(v_1\dots v_{n-3})^{\rho_{n-2}\rho_{n-1}\rho_{n-2}}=(v_1\dots v_{n-3})^{\rho_{n-2}}$ and therefore by the remark \ref{rmk2} the braids $\varrho^*(\beta)$ and $\varrho^*(\alpha)$ are equal modulo $UVP_m^{\prime}$, i.~e. $\overline{\varrho^*(\alpha)}=\overline{\varrho^*(\beta)}$.
\end{itemize}

\underline{Case 3.1.5.2.}  $\iota(\alpha)(s-1)= s-1$.
Since $\beta=\rho_{s-1}\alpha\rho_{s-1}$ then $\iota(\beta)$ fixes $s$ and maps $s-1$ to $k_s$.
$$ \begin{CD}
s @>~\rho_{s-1}~>> s-1 @>~\alpha~>> s-1@>~\rho_{s-1}~>> s
\end{CD}$$
$$ \begin{CD}
s-1 @>~\rho_{s-1}~>> s @>~\alpha~>> k_s (\leq s-2)@>~\rho_{s-1}~>> k_s
\end{CD}$$
Therefore $s-1$ is maximal number which is not fixed by $\beta$, and $\beta_1$ has the following form.
\begin{align}\notag \beta_1&=B_{1,n}^{s-1-n}\beta B_{1,n}^{n-s+1}=B_{1,n}^{-1}B_{1,n}^{s-n}\rho_{s-1}\alpha\rho_{s-1}B_{1,n}^{n-s}B_{1,n}\\
\notag&=B_{1,n}^{-1}\rho_{n-1}B_{1,n}^{s-n}\alpha B_{1,n}^{n-s}\rho_{n-1}B_{1,n}=B_{1,n-1}^{-1}\alpha_1B_{1,n-1}
\end{align}
If we denote by $\delta_{1}=\beta_1$, $\delta_{2}=\delta_{1}^{\rho_{1}}$, $\delta_{3}=\delta_{2}^{\rho_{2}}$, $\dots$, $\delta_{n-1}=\delta_{n-2}^{\rho_{n-2}}=\alpha_1$, then it is obvious that $n$ is not fixed by $\iota(\delta_j)$ for every $j$. Since $\delta_{j+1}$ is obtained from $\delta_j$ conjugating by $\rho_j~(j<n-1)$, then by the cases 3.1.1--3.1.4 the braids $\overline{\varrho^*(\delta_{j+1})}$ and $\overline{\varrho^*(\delta_{j})}$ are conjugated by the element from $S_n$, i.~e.
 $\overline{\varrho^*(\delta_j)}=\overline{\varrho^*(\delta_{j+1})}^{\mu_j}$.
 Therefore we have
 $$\overline{\varrho^*(\beta_1)}=\overline{\varrho^*(\delta_1)}=\overline{\varrho^*(\delta_2)}^{\mu_1}=\overline{\varrho^*(\delta_3)}^{\mu_2\mu_1}=\dots=\overline{\varrho^*(\delta_{n-1})}^{\mu_{n-2}\dots\mu_{1}}=\overline{\varrho^*(\alpha_1)}^{\mu_{n-2}\dots\mu_{1}},$$
 so the braids $\overline{\varrho^*(\beta)}=\overline{\varrho^*(\beta_1)}$ and $\overline{\varrho^*({\alpha})}=\overline{\varrho^*(\alpha_1)}$ are conjugated by the element from $S_n$.

\underline{Case 3.1.5.3.}  $\iota(\alpha)(s-1)\leq s-2$.

\underline{Case 3.1.5.3.1.}  $\iota(\alpha)(s-1)\geq \iota(\alpha)(s)+1$. By the equalities (\ref{a1}) and (\ref{a2}) the braids $\alpha_1$ and $\varrho(\alpha)$ have the following forms.
$$\alpha_1=\gamma w_1\dots w_{n-1}\rho_{n-1}B_{k_n,n-1}$$
$$\varrho(\alpha)=\gamma (w_1\dots w_{n-2})^{\rho_{n-1}}B_{k_n,n-1}$$
Let us rewrite this equalities in more details. Since $\iota(\alpha)(s-1)\leq s-2$, then $\iota(\alpha_1)(n-1)=n-s+\iota(\alpha)(s-1)$.
  $$ \begin{CD}
n-1 @>~B_{1,n}^{s-n}~>> s-1 @>~\alpha~>> \iota(\alpha)(s-1)@>~B_{1,n}^{n-s}~>> n-s+\iota(\alpha)(s-1)
\end{CD}$$
 Therefore since $\iota(\alpha)(s-1)\geq \iota(\alpha)(s)+1$, then the braid $\gamma=\alpha_1\left(w_1\dots w_{n-1}\rho_{n-1}B_{k_n,n-1}\right)^{-1}$ maps $n-1$ to $k_{n-1}=n-s+\iota(\alpha)(s-1)-1\geq n-s+\iota(\alpha)(s)=k_n$
$$ \begin{CD}
n-1 @>~\alpha_1~>> n-s+\iota(\alpha)(s-1) @>~\left(w_1\dots w_{n-1}\rho_{n-1}B_{k_n,n-1}\right)^{-1}~>> n-s+\iota(\alpha)(s-1)-1,
\end{CD}$$
 and by Lemma \ref{lem2} we have
$$\gamma=\eta v_1\dots v_{n-2}B_{k_{n-1},n-1},$$
where $\eta$ does not involve the strands $n-1,n$ and $v_j$ belongs to $\langle\lambda_{j,n-1},\lambda_{n-1,j}\rangle$ for $j=1,\dots,n-1$.
Therefore
$$\alpha_1=\eta v_1\dots v_{n-2}B_{k_{n-1},n-1} w_1\dots w_{n-1}\rho_{n-1}B_{k_{n},n-1}$$
$$\varrho(\alpha)=\eta v_1\dots v_{n-2} B_{k_{n-1},n-1}(w_1\dots w_{n-2})^{\rho_{n-1}}B_{k_{n},n-1}$$
and $k_{n-1}\geq k_n$.

Let us count $\varrho^2(\alpha)$

\begin{itemize}
\item[Step 1] The permutation  $\iota(\varrho(\alpha))$ maps $n-1$ to $k_{n-1}+1=n-s+\iota(\alpha)(s-1)-1+1=n-s+\iota(\alpha)(s-1)\neq n$, therefore $n-1$ is a maximal number which is not foxed by $\varrho(\alpha)$.
\item[Step 2] $\varrho(\alpha)_1=\varrho(\alpha)$.
\item[Step 3] The braid $\varrho(\alpha)_1=\varrho(\alpha)$ can be rewritten.
\begin{align}\notag\varrho(\alpha_1)&=\eta v_1\dots v_{n-2} B_{k_{n-1},n-1}(w_1\dots w_{n-2})^{\rho_{n-1}}B_{k_{n},n-1}\\
\notag&=\eta v_1\dots v_{n-2} (w_1\dots w_{n-2})^{\rho_{n-1}B_{k_{n-1},n-1}^{-1}}B_{k_{n-1},n-1}B_{k_{n},n-1}
\end{align}
Since $k_{n-1}\geq k_n$, then $B_{k_{n-1},n-1}B_{k_{n},n-1}=B_{k_n,n-2}B_{k_{n-1}+1,n-1}$ and
\begin{align}\notag\varrho(\alpha_1)&=\eta v_1\dots v_{n-2} (w_1\dots w_{n-2})^{\rho_{n-1}B_{k_{n-1},n-1}^{-1}}B_{k_n,n-2}B_{k_{n-1}+1,n-1}\\
\notag&=\eta(w_1\dots w_{k_{n-1}-1}w_{k_{n-1}+1}\dots w_{n-2})^{\rho_{n-1}B_{k_{n-1},n-1}^{-1}} v_1\dots v_{n-2}w_{k_{n-1}}^{\rho_{n-1}B_{k_{n-1},n-1}^{-1}} \\
\notag &~\cdot~B_{k_n,n-2}B_{k_{n-1}+1,n-1}\\
\notag&=\eta(w_1\dots w_{k_{n-1}-1}w_{k_{n-1}+1}\dots w_{n-2})^{\rho_{n-1}B_{k_{n-1},n-1}^{-1}}B_{k_n,n-2}\\
 \notag&~\cdot~ (v_1\dots v_{n-2})^{B_{k_n,n-2}}w_{k_{n-1}}^{\rho_{n-1}B_{k_{n-1},n-1}^{-1}B_{k_n,n-2}} B_{k_{n-1}+1,n-1}\\
 \notag&=\eta(w_1\dots w_{k_{n-1}-1}w_{k_{n-1}+1}\dots w_{k_{n-1}})^{\rho_{n-1}B_{k_{n-1},n-1}^{-1}}B_{k_n,n-2}\\
 \notag&~\cdot~ (v_1\dots v_{n-4}v_{n-2})^{B_{k_n,n-2}}w_{k_{n-1}}^{\rho_{n-1}B_{k_{n-1},n-1}^{-1}B_{k_n,n-2}} v_{n-3}^{B_{k_n,n-2}} B_{k_{n-1}+1,n-1}
\end{align}

By Lemma \ref{lem1} the braid
$\eta(w_1\dots w_{k_{n-1}-1}w_{k_{n-1}+1}\dots w_{n-2})^{\rho_{n-1}B_{k_{n-1},n-1}^{-1}}B_{k_n,n-2}$ does not involve the strand $n-1$, the braid
$(v_1\dots v_{n-4}v_{n-2})^{B_{k_n,n-2}}w_{k_{n-1}}^{\rho_{n-1}B_{k_{n-1},n-1}^{-1}B_{k_n,n-2}}$ belongs to $\langle\lambda_{1,n-1},\lambda_{n-1,1}\rangle\times\dots\times\langle\lambda_{n-3,n-1},\lambda_{n-1,n-3}\rangle$ and
$v_{n-3}^{B_{k_n,n-2}}\in\lambda_{n-2,n-1},\lambda_{n-1,n-2}\rangle$.
\item[Step 4] The braid $\varrho(\alpha)_2$ has the following form
\begin{align}
 \notag&\varrho(\alpha)_2=\eta(w_1\dots w_{k_{n-1}-1}w_{k_{n-1}+1}\dots w_{n-2})^{\rho_{n-1}B_{k_{n-1},n-1}^{-1}}B_{k_n,n-2}\\
 \notag&~\cdot~ (v_1\dots v_{n-4}v_{n-2})^{B_{k_n,n-2}}w_{k_{n-1}}^{\rho_{n-1}B_{k_{n-1},n-1}^{-1}B_{k_n,n-2}} \rho_{n-2} B_{k_{n-1}+1,n-2}
\end{align}
\item[Step 5]Finally, the braid $\varrho^2(\alpha)$ follows
\begin{align}
 \notag\varrho^2(\alpha)&=\eta(w_1\dots w_{k_{n-1}-1}w_{k_{n-1}+1}\dots w_{n-2})^{\rho_{n-1}B_{k_{n-1},n-1}^{-1}}B_{k_n,n-2}\\
 \notag& ~\cdot~(v_1\dots v_{n-4}v_{n-2})^{B_{k_n,n-2}\rho_{n-2}}w_{k_{n-1}}^{\rho_{n-1}B_{k_{n-1},n-1}^{-1}B_{k_n,n-2}\rho_{n-2}}  B_{k_{n-1}+1,n-2}\\
 \notag&=\eta(w_1\dots w_{k_{n-1}-1}w_{k_{n-1}+1}\dots w_{n-2})^{\rho_{n-1}B_{k_{n-1},n-1}^{-1}}\\
 \notag& ~\cdot~(v_1\dots v_{n-4}v_{n-2})^{B_{k_n,n-2}\rho_{n-2}B_{k_n,n-2}^{-1}}w_{k_{n-1}}^{\rho_{n-1}B_{k_{n-1},n-1}^{-1}B_{k_n,n-2}\rho_{n-2}B_{k_n,n-2}^{-1}}\\
 \notag&~\cdot~B_{k_n,n-2}B_{k_{n-1}+1,n-2}
\end{align}
\end{itemize}
Let us count $\varrho(\beta)=\varrho(\rho_{s-1}\alpha\rho_{s-1})$.
\begin{itemize}
\item[Step 1] An image of $s$ under the braid $\beta=\rho_{s-1}\alpha\rho_{s-1}$ is equal to $\iota(\alpha)(s-1)$.
     $$ \begin{CD}
s @>~\rho_{s-1}~>> s-1 @>~\alpha~>> \alpha(s-1)@>~\rho_{s-1}~>> \alpha(s-1)
\end{CD}$$
Therefore $s$ is a maximal number which is not fixed by $\beta$.
\item[Step 2] $\beta_1=B_{1,n}^{s-n}\beta B_{1,n}^{n-s}.$
\item[Step 3]
Since $\rho_{s-1}B_{1,n}^{n-s}=B_{1,n}^{n-s}\rho_{n-1}$, then we have
\begin{align}
\notag\beta_1&=B_{1,n}^{s-n}\beta B_{1,n}^{n-s}=B_{1,n}^{s-n}\rho_{s-1}\alpha\rho_{s-1} B_{1,n}^{n-s}=\rho_{n-1}B_{1,n}^{s-n}\alpha B_{1,n}^{n-s}\rho_{n-1}=\rho_{n-1}\alpha_1\rho_{n-1}\\
\notag&=\rho_{n-1}\eta v_1\dots v_{n-2}B_{k_{n-1},n-1} w_1\dots w_{n-1}\rho_{n-1}B_{k_{n},n-1}\rho_{n-1}\\
\notag&=\eta^{\rho_{n-1}} (v_1\dots v_{n-2})^{\rho_{n-1}}B_{k_{n-1},n-1}^{\rho_{n-1}} (w_1\dots w_{n-1})^{\rho_{n-1}}B_{k_{n},n}^{\rho_{n-1}}\\
\notag&=\eta^{\rho_{n-1}} (v_1\dots v_{n-2})^{\rho_{n-1}} (w_1\dots w_{n-1})^{B_{k_{n-1},n-1}^{-1}\rho_{n-1}}B_{k_{n-1},n-1}^{\rho_{n-1}}B_{k_{n},n}^{\rho_{n-1}}
\end{align}
Using simple calculations in the symmetric group, it is easy to see that
\begin{equation}\label{big1}
B_{k_{n-1},n-1}^{\rho_{n-1}}B_{k_{n},n}^{\rho_{n-1}}=B_{k_n,n-1}\rho_{n-2}B_{k_{n-1}+1,n},
\end{equation}
therefore
\begin{align}
\notag\beta_1&=\eta^{\rho_{n-1}} (v_1\dots v_{n-2})^{\rho_{n-1}} (w_1\dots w_{n-1})^{B_{k_{n-1},n-1}^{-1}\rho_{n-1}}B_{k_n,n-1}\rho_{n-2}B_{k_{n-1}+1,n}\\
\notag&=\eta^{\rho_{n-1}}(w_1\dots w_{k_{n-1}-1}w_{k_{n-1}+1}\dots w_{n-1})^{B_{k_{n-1},n-1}^{-1}\rho_{n-1}}\\
\notag&~\cdot~(v_1\dots v_{n-2})^{\rho_{n-1}}w_{k_{n-1}}^{B_{k_{n-1},n-1}^{-1}\rho_{n-1}}B_{k_n,n-1}\rho_{n-2}B_{k_{n-1}+1,n}\\
\notag&=\eta^{\rho_{n-1}}(w_1\dots w_{k_{n-1}-1}w_{k_{n-1}+1}\dots w_{n-1})^{B_{k_{n-1},n-1}^{-1}\rho_{n-1}}B_{k_n,n-1}\rho_{n-2}\\
\notag&~\cdot~(v_1\dots v_{n-2})^{\rho_{n-1}B_{k_n,n-1}\rho_{n-2}}w_{k_{n-1}}^{B_{k_{n-1},n-1}^{-1}\rho_{n-1}B_{k_n,n-1}\rho_{n-2}}B_{k_{n-1}+1,n}\\
\notag&=\eta^{\rho_{n-1}}(w_1\dots w_{k_{n-1}-1}w_{k_{n-1}+1}\dots w_{n-1})^{B_{k_{n-1},n-1}^{-1}\rho_{n-1}}B_{k_n,n-1}\rho_{n-2}\\
\notag&~\cdot~(v_1\dots v_{n-4}v_{n-2})^{\rho_{n-1}B_{k_n,n-1}\rho_{n-2}}w_{k_{n-1}}^{B_{k_{n-1},n-1}^{-1}\rho_{n-1}B_{k_n,n-1}\rho_{n-2}}v_{n-3}^{\rho_{n-1}B_{k_n,n-1}\rho_{n-2}}B_{k_{n-1}+1,n}
\end{align}
By Theorem \ref{t2} the braid
$$\eta^{\rho_{n-1}}(w_1\dots w_{k_{n-1}-1}w_{k_{n-1}+1}\dots w_{n-1})^{B_{k_{n-1},n-1}^{-1}\rho_{n-1}}B_{k_n,n-1}\rho_{n-2}$$ does not involve the strand $n$, the braid
$$(v_1\dots v_{n-4}v_{n-2})^{\rho_{n-1}B_{k_n,n-1}\rho_{n-2}}w_{k_{n-1}}^{B_{k_{n-1},n-1}^{-1}\rho_{n-1}B_{k_n,n-1}\rho_{n-2}}$$ belongs to $\langle\lambda_{1,n},\lambda_{n,1}\rangle\times\dots\times\langle\lambda_{n-2,n},\lambda_{n,n-2}\rangle$ and the braid
$v_{n-3}^{\rho_{n-1}B_{k_n,n-1}\rho_{n-2}}$ belongs to $\langle\lambda_{n-1,n},\lambda_{n,n-1}\rangle$.
\item[Step 4] The braid $\beta_2$ follows
\begin{align}
\notag\beta_2&=\eta^{\rho_{n-1}}(w_1\dots w_{k_{n-1}-1}w_{k_{n-1}+1}\dots w_{n-1})^{B_{k_{n-1},n-1}^{-1}\rho_{n-1}}B_{k_n,n-1}\rho_{n-2}\\
\notag&~\cdot~(v_1\dots v_{n-4}v_{n-2})^{\rho_{n-1}B_{k_n,n-1}\rho_{n-2}}w_{k_{n-1}}^{B_{k_{n-1},n-1}^{-1}\rho_{n-1}B_{k_n,n-1}\rho_{n-2}}\rho_{n-1}B_{k_{n-1}+1,n-1}
\end{align}
\item[Step 5] The braid $\varrho(\beta)$ is the following.
\begin{align}
\notag\varrho(\beta)&=\eta^{\rho_{n-1}}(w_1\dots w_{k_{n-1}-1}w_{k_{n-1}+1}\dots w_{n-1})^{B_{k_{n-1},n-1}^{-1}\rho_{n-1}}B_{k_n,n-1}\rho_{n-2}\\
\notag&~\cdot~(v_1\dots v_{n-4}v_{n-2})^{\rho_{n-1}B_{k_n,n-1}\rho_{n-2}\rho_{n-1}}w_{k_{n-1}}^{B_{k_{n-1},n-1}^{-1}\rho_{n-1}B_{k_n,n-1}\rho_{n-2}\rho_{n-1}}B_{k_{n-1}+1,n-1}
\end{align}
\end{itemize}
Let us find the braid $\varrho^2(\beta)$.
\begin{itemize}
\item[Step 1] Since $\eta^{\rho_{n-1}}$ does not involve the strands $n-1,n$, then the image of $n-1$ under the permutation $\iota(\varrho(\beta))$ is equal to the image of $n-1$ under the permutation $\iota(B_{k_n,n-1}\rho_{n-2}B_{k_{n-1}+1,n-1})$ and is equal to $k_n$. Then $n-1$ is a maximal number which is not fixed by $\varrho(\beta)$.
\item[Step 2]$\varrho(\beta)_1=\varrho(\beta)$.
\item[Step 3] The braid $\varrho(\beta)_1$ can be rewritten
\begin{align}
\notag\varrho(\beta)_1&=\eta^{\rho_{n-1}}(w_1\dots w_{k_{n-1}-1}w_{k_{n-1}+1}\dots w_{n-1})^{B_{k_{n-1},n-1}^{-1}\rho_{n-1}}\\
\notag&~\cdot~(v_1\dots v_{n-4}v_{n-2})^{\rho_{n-1}B_{k_n,n-1}\rho_{n-2}\rho_{n-1}\rho_{n-2}B_{k_n,n-1}^{-1}}\\
\notag& ~\cdot~ w_{k_{n-1}}^{B_{k_{n-1},n-1}^{-1}\rho_{n-1}B_{k_n,n-1}\rho_{n-2}\rho_{n-1}\rho_{n-2}B_{k_n,n-1}^{-1}}\\
\notag&~\cdot~B_{k_n,n-1}\rho_{n-2}B_{k_{n-1}+1,n-1}
\end{align}
The following equality is faithful in the symmetric group
\begin{equation}\label{big2}
B_{k_n,n-1}\rho_{n-2}B_{k_{n-1}+1,n-1}=\rho_{n-3}B_{k_{n-1},n-2}B_{k_n,n-1}=B_{k_{n-1},n-3}B_{k_n,n-1},
\end{equation}
therefore
\begin{align}
\notag\varrho(\beta)_1&=\eta^{\rho_{n-1}}(w_1\dots w_{k_{n-1}-1}w_{k_{n-1}+1}\dots w_{n-1})^{B_{k_{n-1},n-1}^{-1}\rho_{n-1}}\\
\notag&~\cdot~(v_1\dots v_{n-4}v_{n-2})^{\rho_{n-1}B_{k_n,n-1}\rho_{n-2}\rho_{n-1}\rho_{n-2}B_{k_n,n-1}^{-1}}\\
\notag&~\cdot~w_{k_{n-1}}^{B_{k_{n-1},n-1}^{-1}\rho_{n-1}B_{k_n,n-1}\rho_{n-2}\rho_{n-1}\rho_{n-2}B_{k_n,n-1}^{-1}}\\
\notag&~\cdot~B_{k_{n-1},n-3}B_{k_n,n-1}\\
\notag&=\eta^{\rho_{n-1}}(v_1\dots v_{n-4}v_{n-2})^{\rho_{n-1}B_{k_n,n-1}\rho_{n-2}\rho_{n-1}\rho_{n-2}B_{k_n,n-1}^{-1}}\\
\notag&~\cdot~(w_1\dots w_{k_{n-1}-1}w_{k_{n-1}+1}\dots w_{n-2})^{B_{k_{n-1},n-1}^{-1}\rho_{n-1}}\\
\notag&~\cdot~w_{k_{n-1}}^{B_{k_{n-1},n-1}^{-1}\rho_{n-1}B_{k_n,n-1}\rho_{n-2}\rho_{n-1}\rho_{n-2}B_{k_n,n-1}^{-1}}\\
\notag&~\cdot~w_{n-1}^{B_{k_{n-1},n-1}^{-1}\rho_{n-1}}B_{k_{n-1},n-3}B_{k_n,n-1}
\end{align}
\begin{align}
\notag~~~~~~~~~~~~~~~&=\eta^{\rho_{n-1}}(v_1\dots v_{n-4}v_{n-2})^{\rho_{n-1}B_{k_n,n-1}\rho_{n-2}\rho_{n-1}\rho_{n-2}B_{k_n,n-1}^{-1}}B_{k_{n-1},n-3}\\
\notag&~\cdot~(w_1\dots w_{k_{n-1}-1}w_{k_{n-1}+1}\dots w_{n-2})^{B_{k_{n-1},n-1}^{-1}\rho_{n-1}B_{k_{n-1},n-3}}\\
\notag&~\cdot~w_{k_{n-1}}^{B_{k_{n-1},n-1}^{-1}\rho_{n-1}B_{k_n,n-1}\rho_{n-2}\rho_{n-1}\rho_{n-2}B_{k_n,n-1}^{-1}B_{k_{n-1},n-3}}\\
\notag&~\cdot~w_{n-1}^{B_{k_{n-1},n-1}^{-1}\rho_{n-1}B_{k_{n-1},n-3}}B_{k_n,n-1}
\end{align}

Here the braid
$\eta^{\rho_{n-1}}(v_1\dots v_{n-4}v_{n-2})^{\rho_{n-1}B_{k_n,n-1}\rho_{n-2}\rho_{n-1}\rho_{n-2}B_{k_n,n-1}^{-1}}B_{k_{n-1},n-3}$ does not involve the strand $n-1$, the braid
\begin{multline*}
(w_1\dots w_{k_{n-1}-1}w_{k_{n-1}+1}\dots w_{n-2})^{B_{k_{n-1},n-1}^{-1}\rho_{n-1}B_{k_{n-1},n-3}}\cdot\\
\cdot w_{k_{n-1}}^{B_{k_{n-1},n-1}^{-1}\rho_{n-1}B_{k_n,n-1}\rho_{n-2}\rho_{n-1}\rho_{n-2}B_{k_n,n-1}^{-1}B_{k_{n-1},n-3}}
\end{multline*}
belongs to $\langle\lambda_{n-1,1},\lambda_{1,n-1}\rangle\times\dots\times\langle\lambda_{n-3,n-1},\lambda_{n-1,n-3}\rangle$ and
$w_{n-1}^{B_{k_{n-1},n-1}^{-1}\rho_{n-1}B_{k_{n-1},n-3}}$ belongs to $\langle\lambda_{n-2,n-1},\lambda_{n-1,n-2}\rangle$.
\item[Step 4] The braid $\varrho(\beta)_2$ follows.
\begin{align}
\notag\varrho(\beta)_2&=\eta^{\rho_{n-1}}(v_1\dots v_{n-4}v_{n-2})^{\rho_{n-1}B_{k_n,n-1}\rho_{n-2}\rho_{n-1}\rho_{n-2}B_{k_n,n-1}^{-1}}B_{k_{n-1},n-3}\\
\notag&~\cdot~(w_1\dots w_{k_{n-1}-1}w_{k_{n-1}+1}\dots w_{n-2})^{B_{k_{n-1},n-1}^{-1}\rho_{n-1}B_{k_{n-1},n-3}}\\
\notag&~\cdot~w_{k_{n-1}}^{B_{k_{n-1},n-1}^{-1}\rho_{n-1}B_{k_n,n-1}\rho_{n-2}\rho_{n-1}\rho_{n-2}B_{k_n,n-1}^{-1}B_{k_{n-1},n-3}}\rho_{n-2}B_{k_n,n-2}
\end{align}
\item[Step 5] The braid $\varrho^2(\beta)$ has the following form.
\begin{align}
\notag\varrho^2(\beta)&=\eta^{\rho_{n-1}}(v_1\dots v_{n-4}v_{n-2})^{\rho_{n-1}B_{k_n,n-1}\rho_{n-2}\rho_{n-1}\rho_{n-2}B_{k_n,n-1}^{-1}}B_{k_{n-1},n-3}\\
\notag&~\cdot~(w_1\dots w_{k_{n-1}-1}w_{k_{n-1}+1}\dots w_{n-2})^{B_{k_{n-1},n-1}^{-1}\rho_{n-1}B_{k_{n-1},n-3}\rho_{n-2}}\\
\notag&~\cdot~w_{k_{n-1}}^{B_{k_{n-1},n-1}^{-1}\rho_{n-1}B_{k_n,n-1}\rho_{n-2}\rho_{n-1}\rho_{n-2}B_{k_n,n-1}^{-1}B_{k_{n-1},n-3}\rho_{n-2}}B_{k_n,n-2}\\
\notag&=\eta^{\rho_{n-1}}(v_1\dots v_{n-4}v_{n-2})^{\rho_{n-1}B_{k_n,n-1}\rho_{n-2}\rho_{n-1}\rho_{n-2}B_{k_n,n-1}^{-1}}\\
\notag&~\cdot~(w_1\dots w_{k_{n-1}-1}w_{k_{n-1}+1}\dots w_{n-2})^{B_{k_{n-1},n-1}^{-1}\rho_{n-1}B_{k_{n-1},n-3}\rho_{n-2}B_{k_{n-1},n-3}^{-1}}\\
\notag&~\cdot~w_{k_{n-1}}^{B_{k_{n-1},n-1}^{-1}\rho_{n-1}B_{k_n,n-1}\rho_{n-2}\rho_{n-1}\rho_{n-2}B_{k_n,n-1}^{-1}B_{k_{n-1},n-3}\rho_{n-2}B_{k_{n-1},n-3}^{-1}}\\
\notag&~\cdot~B_{k_{n-1},n-3}B_{k_n,n-2}
\end{align}
\begin{align}
\notag~~~~~~~~~~~&=\eta^{\rho_{n-1}}(w_1\dots w_{k_{n-1}-1}w_{k_{n-1}+1}\dots w_{n-2})^{B_{k_{n-1},n-1}^{-1}\rho_{n-1}B_{k_{n-1},n-3}\rho_{n-2}B_{k_{n-1},n-3}^{-1}}\\
\notag&~\cdot~(v_1\dots v_{n-4}v_{n-2})^{\rho_{n-1}B_{k_n,n-1}\rho_{n-2}\rho_{n-1}\rho_{n-2}B_{k_n,n-1}^{-1}}\\
\notag&~\cdot~[(v_1\dots v_{n-4}v_{n-2})^{\rho_{n-1}B_{k_n,n-1}\rho_{n-2}\rho_{n-1}\rho_{n-2}B_{k_n,n-1}^{-1}},\\
\notag&~~~~(w_1\dots w_{k_{n-1}-1}w_{k_{n-1}+1}\dots w_{n-2})^{B_{k_{n-1},n-1}^{-1}\rho_{n-1}B_{k_{n-1},n-3}\rho_{n-2}B_{k_{n-1},n-3}^{-1}}]\\
\notag&~\cdot~w_{k_{n-1}}^{B_{k_{n-1},n-1}^{-1}\rho_{n-1}B_{k_n,n-1}\rho_{n-2}\rho_{n-1}\rho_{n-2}B_{k_n,n-1}^{-1}B_{k_{n-1},n-3}\rho_{n-2}B_{k_{n-1},n-3}^{-1}}\\
\notag&~\cdot~B_{k_{n-1},n-3}B_{k_n,n-2}
\end{align}
\end{itemize}
Now we have
\begin{align}
\notag\varrho^2(\beta)&=\eta^{\rho_{n-1}}(w_1\dots w_{k_{n-1}-1}w_{k_{n-1}+1}\dots w_{n-2})^{B_{k_{n-1},n-1}^{-1}\rho_{n-1}B_{k_{n-1},n-3}\rho_{n-2}B_{k_{n-1},n-3}^{-1}}\\
\notag&~\cdot~(v_1\dots v_{n-4}v_{n-2})^{\rho_{n-1}B_{k_n,n-1}\rho_{n-2}\rho_{n-1}\rho_{n-2}B_{k_n,n-1}^{-1}}\\
\notag&~\cdot~[(v_1\dots v_{n-4}v_{n-2})^{\rho_{n-1}B_{k_n,n-1}\rho_{n-2}\rho_{n-1}\rho_{n-2}B_{k_n,n-1}^{-1}},\\
\notag&~~~~(w_1\dots w_{k_{n-1}-1}w_{k_{n-1}+1}\dots w_{n-2})^{B_{k_{n-1},n-1}^{-1}\rho_{n-1}B_{k_{n-1},n-3}\rho_{n-2}B_{k_{n-1},n-3}^{-1}}]\\
\notag&~\cdot~w_{k_{n-1}}^{B_{k_{n-1},n-1}^{-1}\rho_{n-1}B_{k_n,n-1}\rho_{n-2}\rho_{n-1}\rho_{n-2}B_{k_n,n-1}^{-1}B_{k_{n-1},n-3}\rho_{n-2}B_{k_{n-1},n-3}^{-1}}\\
\notag&~\cdot~B_{k_{n-1},n-3}B_{k_n,n-2}\\
 \notag\varrho^2(\alpha)&=\eta(w_1\dots w_{k_{n-1}-1}w_{k_{n-1}+1}\dots w_{n-2})^{\rho_{n-1}B_{k_{n-1},n-1}^{-1}}\\
 \notag& (v_1\dots v_{n-4}v_{n-2})^{B_{k_n,n-2}\rho_{n-2}B_{k_n,n-2}^{-1}}w_{k_{n-1}}^{\rho_{n-1}B_{k_{n-1},n-1}^{-1}B_{k_n,n-2}\rho_{n-2}B_{k_n,n-2}^{-1}}\\
 \notag&B_{k_n,n-2}B_{k_{n-1}+1,n-2}
\end{align}
Note that the following equalities are faithful in the symmetric group
\begin{align}
\notag B_{k_{n-1},n-3}B_{k_n,n-2}&=B_{k_n,n-2}B_{k_{n-1}+1,n-2}\\
\notag\rho_{n-1}B_{k_n,n-1}\rho_{n-2}\rho_{n-1}\rho_{n-2}B_{k_n,n-1}^{-1}&=\rho_{n-1}\underline{\rho_{n-2}\rho_{n-3}\rho_{n-2}}\rho_{n-1}\rho_{n-2}\rho_{n-3}\rho_{n-2}\\
\notag&=\rho_{n-1}\rho_{n-3}\rho_{n-2}\underline{\rho_{n-3}\rho_{n-1}}\rho_{n-2}\rho_{n-3}\rho_{n-2}\\
\notag&=\rho_{n-1}\rho_{n-3}\rho_{n-2}\rho_{n-1}\rho_{n-3}\underline{\rho_{n-2}\rho_{n-3}\rho_{n-2}}\\
\notag&=\rho_{n-1}\rho_{n-3}\rho_{n-2}\rho_{n-1}\underline{\rho_{n-3}\rho_{n-3}}\rho_{n-2}\rho_{n-3}\\
\notag&=\underline{\rho_{n-1}\rho_{n-3}}\rho_{n-2}\rho_{n-1}\rho_{n-2}\rho_{n-3}\\
\notag&=\rho_{n-3}\rho_{n-1}\underline{\rho_{n-2}\rho_{n-1}\rho_{n-2}}\rho_{n-3}\\
\notag&=\rho_{n-3}\underline{\rho_{n-1}\rho_{n-1}}\rho_{n-2}\rho_{n-1}\rho_{n-3}\\
\notag&=\rho_{n-3}\rho_{n-2}\underline{\rho_{n-1}\rho_{n-3}}\\
\notag&=\rho_{n-3}\rho_{n-2}\rho_{n-3}\rho_{n-1}\\
\notag B_{k_n,n-2}\rho_{n-2}B_{k_n,n-2}^{-1}&=\rho_{n-3}\rho_{n-2}\rho_{n-3}
\end{align}
Therefore we have
$$(v_1\dots v_{n-4}v_{n-2})^{\rho_{n-1}B_{k_n,n-1}\rho_{n-2}\rho_{n-1}\rho_{n-2}B_{k_n,n-1}^{-1}}=((v_1\dots v_{n-4}v_{n-2})^{\rho_{n-3}\rho_{n-2}\rho_{n-3}})^{\rho_{n-1}},$$
and since $(v_1\dots v_{n-4}v_{n-2})^{\rho_{n-3}\rho_{n-2}\rho_{n-3}}$ does not involve the strands $n-1,n$, then
\begin{multline*}((v_1\dots v_{n-4}v_{n-2})^{\rho_{n-3}\rho_{n-2}\rho_{n-3}})^{\rho_{n-1}}=\\
=(v_1\dots v_{n-4}v_{n-2})^{\rho_{n-3}\rho_{n-2}\rho_{n-3}}=\\
=(v_1\dots v_{n-4}v_{n-2})^{B_{k_n,n-2}\rho_{n-2}B_{k_n,n-2}^{-1}}
\end{multline*}
and we have
$$(v_1\dots v_{n-4}v_{n-2})^{\rho_{n-1}B_{k_n,n-1}\rho_{n-2}\rho_{n-1}\rho_{n-2}B_{k_n,n-1}^{-1}}=(v_1\dots v_{n-4}v_{n-2})^{B_{k_n,n-2}\rho_{n-2}B_{k_n,n-2}^{-1}}.$$
Also, since
\begin{align}
\notag B_{k_{n-1},n-1}^{-1}\rho_{n-1}B_{k_{n-1},n-3}\rho_{n-2}B_{k_{n-1},n-3}^{-1}&=B_{k_{n-1},n-1}^{-1}\rho_{n-1}\rho_{n-2}\\
\notag&=B_{k_{n-1},n-1}^{-1}\rho_{n-1}\rho_{n-2}\\
\notag&=\rho_{k_{n-1}}\dots\rho_{n-3}\rho_{n-2}\rho_{n-1}\rho_{n-2}\\
\notag&=\rho_{k_{n-1}}\dots\rho_{n-3}\rho_{n-1}\rho_{n-2}\rho_{n-1}\\
\notag&=\rho_{n-1}\rho_{k_{n-1}}\dots\rho_{n-3}\rho_{n-2}\rho_{n-1}=\rho_{n-1}B_{k_{n-1},n-1}^{-1}\rho_{n-1},
\end{align}
then
\begin{multline*}(w_1\dots w_{k_{n-1}-1}w_{k_{n-1}+1}\dots w_{n-2})^{B_{k_{n-1},n-1}^{-1}\rho_{n-1}B_{k_{n-1},n-3}\rho_{n-2}B_{k_{n-1},n-3}^{-1}}=\\
=\left((w_1\dots w_{k_{n-1}-1}w_{k_{n-1}+1}\dots w_{n-2})^{\rho_{n-1}B_{k_{n-1},n-1}^{-1}}\right)^{\rho_{n-1}},
\end{multline*}
and since the braid $(w_1\dots w_{k_{n-1}-1}w_{k_{n-1}+1}\dots w_{n-2})^{\rho_{n-1}B_{k_{n-1},n-1}^{-1}}$ does not involve the strands $n-1,n$, then
\begin{multline*}(w_1\dots w_{k_{n-1}-1}w_{k_{n-1}+1}\dots w_{n-2})^{B_{k_{n-1},n-1}^{-1}\rho_{n-1}B_{k_{n-1},n-3}\rho_{n-2}B_{k_{n-1},n-3}^{-1}}=\\
=(w_1\dots w_{k_{n-1}-1}w_{k_{n-1}+1}\dots w_{n-2})^{\rho_{n-1}B_{k_{n-1},n-1}^{-1}}
\end{multline*}
\begin{multline*}
w_{k_{n-1}}^{B_{k_{n-1},n-1}^{-1}\rho_{n-1}B_{k_n,n-1}\rho_{n-2}\rho_{n-1}\rho_{n-2}B_{k_n,n-1}^{-1}B_{k_{n-1},n-3}\rho_{n-2}B_{k_{n-1},n-3}^{-1}}=\\
=w_{k_{n-1}}^{\rho_{n-1}B_{k_{n-1},n-1}^{-1}B_{k_n,n-2}\rho_{n-2}B_{k_n,n-2}^{-1}}
\end{multline*}
Finally since $\eta$ does not involve the strands $n-1,n$, then $\eta^{\rho_{n-1}}=\eta$. Therefore the braids $\varrho^2(\alpha)$ and $\varrho^2(\beta)$ are equal modulo $UVP_{n-2}^{\prime}$ and by the remark \ref{rmk2} the braids $\varrho^*(\beta)$ and $\varrho^*(\alpha)$ are equal modulo $UVP_m^{\prime}$, i.~e. $\overline{\varrho^*(\alpha)}=\overline{\varrho^*(\beta)}$.

\underline{Case 3.1.5.3.2.}  $\alpha(s-1)< \alpha(s)+1$. This case literally repeats the case 3.1.5.3.1 using the fact that in this case $k_{n-1}<k_n$ and then using the equalities
$$B_{k_{n-1},n-1}B_{k_{n},n-1}=B_{k_n-1,n-2}B_{k_{n-1},n-1}$$
$$B_{k_{n-1},n-1}^{\rho_{n-1}}B_{k_{n},n}^{\rho_{n-1}}=B_{k_n-1,n-1}\rho_{n-2}B_{k_{n-1},n}$$
instead of equalities (\ref{big1}) and (\ref{big2}).

\underline{Case 3.1.6.}  $i=s$. In this case the braid $\beta=\rho_{s}\alpha\rho_{s}$ maps $s+1$ to $k_s$.
$$\begin{CD}
s+1 @>~\rho_{s}~>> s @>~\alpha~>> k_s@>~\rho_{s}~>> k_s
\end{CD}$$
Therefore $s+1$ is a maximal number, which is not fixed by $\beta$ and we have
\begin{align}\notag \beta_1&=B_{1,n}^{s+1-n}\beta B_{1,n}^{n-s-1}=B_{1,n}B_{1,n}^{s-n}\rho_{s-1}\alpha\rho_{s-1}B_{1,n}^{n-s}B_{1,n}^{-1}\\
\notag&=B_{1,n}\rho_{n-1}B_{1,n}^{s-n}\alpha B_{1,n}^{n-s}\rho_{n-1}B_{1,n}^{-1}\\
\notag&=B_{1,n}\rho_{n-1}\alpha_1\rho_{n-1}B_{1,n}^{-1}=\rho_{n-2}B_{1,n}\alpha_1B_{1,n}^{-1}\rho_{n-2}
\end{align}
If we denote by  $\delta_{1}=\alpha_1$,
  $\delta_{2}=\delta_{1}^{\rho_{1}}$,
 $\delta_{3}=\delta_{2}^{\rho_{2}}$,
 $\dots$,
  $\delta_{n-1}=\delta_{n-2}^{\rho_{n-2}}$,
$\delta_{n}=\delta_{n-1}^{\rho_{n-1}}=\alpha_1^{B_{1,n}^{-1}}$,
$\delta_{n+1}=\delta_n^{\rho_{n-2}}=\beta_1$, then it is obvious, that   $\iota(\delta_i)(n)\neq n$ for $i=1\dots n-1$.
  Therefore by the previous cases $\overline{\varrho^*(\delta_i)}$ and $\overline{\varrho^*(\delta_{i+1})}$ are conjugated by the element from $S_n$ for every $i=1\dots n-1$.
  $$\varrho(\delta_i)=\varrho(\delta_{i+1})^{\mu_i}$$
 Also since $\iota(\delta_n)(n-1)\neq n-1$
 $$\begin{CD}
n-1 @>~B_{1,n}~>> n @>~\alpha_1~>> k_n\neq n@>~B_{1,n}~>> k_n-1~(\text{mod}~n)\neq n-1
\end{CD}$$
then by the previous cases $\overline{\varrho^*(\delta_n)}=\overline{\varrho^*(\delta_{n+1})}^{\mu_n}$ and we have
 $$\overline{\varrho^*(\beta_1)}=\overline{\varrho^*(\delta_1)}=\overline{\varrho^*(\delta_2)}^{\mu_1}=\overline{\varrho^*(\delta_3)}^{\mu_2\mu_1}=\dots=\overline{\varrho^*(\delta_{n-1})}^{\mu_{n-2}\dots\mu_{1}}=\overline{\varrho^*(\alpha_1)}^{\mu_{n-2}\dots\mu_{1}}$$
 Therefore the braids $\overline{\varrho^*(\beta)}=\overline{\varrho^*(\beta_1)}$ and $\overline{\varrho^*({\alpha})}=\overline{\varrho^*(\alpha_1)}$ are conjugated by the element from $S_n$.

 \underline{Case 3.1.7.}  $i\geq s+1$. In this case the maximal number which is not fixed by the braid $\beta=\rho_{i}\alpha\rho_{i}$ is $s$ and we have
$$ \begin{CD}
s @>~\rho_{i}~>> s @>~\alpha~>> k_s@>~\rho_{i}~>> k_s
\end{CD}$$
Since $i\geq s+2$, then we have
\begin{align}
\notag\rho_iB_{1,n}^{n-s}&=\rho_{i}B_{1,n}^{n-i-1}B_{1,n}B_{1,n}^{i-s}=B_{1,n}^{n-i-1}\rho_{n-1}B_{1,n}B_{1,n}^{i-s}\\
\notag&=B_{1,n}^{n-i-1}B_{1,n-1}B_{1,n}^{i-s}=B_{1,n}^{n-i-1}B_{1,n-1}B_{1,n}B_{1,n}^{i-s-1}\\
\notag&=B_{1,n}^{n-i-1}B_{1,n}B_{2,n}B_{1,n}^{i-s-1}=B_{1,n}^{n-i-1}B_{1,n}B_{2,n}\rho_1\rho_1B_{1,n}^{i-s-1}\\
\notag&=B_{1,n}^{n-i+1}\rho_1B_{1,n}^{i-s-1}=B_{1,n}^{n-i+1}B_{1,n}^{i-s-1}\rho_{i-s}=B_{1,n}^{n-s}\rho_{i-s}
\end{align}
Therefore the braid $\beta_1$ can be rewritten
\begin{align}\notag\beta_1&=B_{1,n}^{s-n}\beta B_{1,n}^{n-s}=B_{1,n}^{s-n}\rho_i\alpha\rho_i B_{1,n}^{n-s}=\rho_{i-s}B_{1,n}^{s-n}\alpha B_{1,n}^{n-s}\rho_{i-s}=\rho_{i-s}\alpha_1\rho_{i-s}
\end{align}
Since $\iota(\alpha_1)(n)\neq n$ and $i-s<n$, then $\overline{\varrho^*(\alpha_1)}$ and $\overline{\varrho^*(\beta_1)}$ are conjugated by the previous cases.

\underline{Case 3.2.} The braid $\beta$ is obtained from $\alpha$ conjugating by $\lambda_{i,i+1}^{\pm1}$. If $\beta$ is obtained from $\alpha$ conjugating by $\lambda_{i,i+1}^{-1}$, then $\alpha$ is obtained from $\beta$ conjugating by $\lambda_{i,i+1}$ and we can consider that $\beta$ is obtained from $\alpha$ conjugating by $\lambda_{i,i+1}$.

As we already found in the case 3.1 (equalities (\ref{a2}) and (\ref{a1})) the braids $\varrho(\alpha)$ and $\alpha_1$  has the following forms
$$\alpha_1=\gamma w_1\dots w_{n-1}\rho_{n-1}B_{k_n,n-2},$$
$$\varrho(\alpha)=\gamma (w_1\dots w_{n-2})^{\rho_{n-1}}B_{k_n,n-1},$$
where $k_n=n-s+k_s$ and $k_s$ is a maximal number which is not fixed by $\alpha$.

\underline{Case 3.2.1.} $i\leq k_s-2$. Let us count the braid $\varrho(\beta)=\varrho(\lambda_{i,i+1}^{-1}\alpha\lambda_{i,i+1})$.
\begin{itemize}
\item[Step 1.] Since $\lambda_{i,i+1}$ is a pure braid, then the maximal number which is not fixed by $\beta$ is equal to the maximal number which is not fixed by $\beta$ and is equal to $s$.
\item[Step 2.] We have $\beta_1=B_{1,n}^{s-n}\beta B_{1,n}^{n-s}$ and $\iota(\beta_1)(n)=k_s+n-s=k_n\neq n$.
\item[Step 3.] By Lemma \ref{lem1} we have $\lambda_{i,i+1}B_{1,n}^{n-s}=B_{1,n}^{n-s}\lambda_{i+n-s,i+n-s+1}$. If we denote by $j=i+n-s$, then the braid $\beta_1$ can be rewritten in more details
\begin{align}
\notag\beta_1&=B_{1,n}^{s-n}\beta B_{1,n}^{n-s}=B_{1,n}^{s-n}\lambda_{i,i+1}^{-1}\alpha\lambda_{i,i+1} B_{1,n}^{n-s}\\
\notag&=\lambda_{j,j+1}^{-1}B_{1,n}^{s-n}\alpha B_{1,n}^{n-s}\lambda_{j,j+1}=\lambda_{j,j+1}^{-1}\alpha_1\lambda_{j,j+1}\\
\notag&=\lambda_{j,j+1}^{-1}\gamma w_1\dots w_{n-1}B_{k_n,n}\lambda_{j,j+1}\\
\notag&=\lambda_{j,j+1}^{-1}\gamma w_1\dots w_{n-1}\lambda_{j,j+1}^{B_{k_n,n}^{-1}}B_{k_n,n}
\end{align}
Since $i<k_s-1$, then $j=i+n-s<k_n-1$ and therefore $\lambda_{j,j+1}^{B_{k_n,n}^{-1}}=\lambda_{j,j+1}$. Therefore we have
$$\beta_1=\lambda_{j,j+1}^{-1}\gamma\lambda_{j,j+1} w_1\dots w_{n-1}B_{k_n,n},$$
where the braid $\lambda_{j,j+1}^{-1}\gamma\lambda_{j,j+1}$ does not involve the strand $n$ and $w_r$ belongs to $\langle\lambda_{r,n},\lambda_{n,r}\rangle$ for $r=1,\dots,w_{n-1}$.
\item[Step 4.] $\beta_2=\lambda_{j,j+1}^{-1}\gamma\lambda_{j,j+1} w_1\dots w_{n-2}\rho_{n-1}B_{k_n,n-1}.$
\item[Step 5.] The braid $\varrho(\beta)$ has the following form.
\begin{align}
\notag\varrho(\beta)&=\lambda_{j,j+1}^{-1}\gamma\lambda_{j,j+1} (w_1\dots w_{n-2})^{\rho_{n-1}}B_{k_n,n-1}\\
\notag&=\lambda_{j,j+1}^{-1}\gamma (w_1\dots w_{n-2})^{\rho_{n-1}}B_{k_n,n-1}\lambda_{j,j+1}=\varrho(\alpha)^{\lambda_{j,j+1}}
\end{align}
\end{itemize}
Therefore $\varrho(\alpha)$ and $\varrho(\beta)$ are conjugated and by the induction hypothesis the braids $\varrho^*(\alpha)$ and $\varrho^*(\beta)$ are conjugated by the element from $S_n$.

\underline{Case 3.2.2.} $i=k_s-1$. Let us find the braid $\varrho(\beta)$.
\begin{itemize}
\item[Step 1.] The maximal number which is not fixed by $\alpha$ is equal to $s$.
\item[Step 2.] We have $\beta_1=B_{1,n}^{s-n}\beta B_{1,n}^{n-s}$ and $\iota(\beta_1)(n)=k_s+n-s=k_n\neq n$.
\item[Step 3.] Similarly to the case 3.2.1 we conclude that $\beta_1=\lambda_{k_n-1,k_n}^{-1}\alpha_1\lambda_{k_n-1,k_n}$ and therefore we have
\begin{align}
\notag\beta_1&=\lambda_{k_n-1,k_n}^{-1}\gamma w_1\dots w_{n-1}B_{k_n,n}\lambda_{k_n-1,k_n}\\
\notag&=\lambda_{k_n-1,k_n}^{-1}\gamma w_1\dots w_{n-1}\lambda_{k_n-1,k_n}^{B_{k_n,n}^{-1}}B_{k_n,n}\\
\notag&=\lambda_{k_n-1,k_n}^{-1}\gamma w_1\dots w_{n-1}\lambda_{k_n-1,n}B_{k_n,n}\\
\notag&=\lambda_{k_n-1,k_n}^{-1}\gamma w_1\dots w_{n-2}\lambda_{k_n-1,n}w_{n-1}B_{k_n,n}
\end{align}
Since $k_n<n$, then the braid $\lambda_{k_n-1,k_n}^{-1}\gamma$ does not involve the strand $n$.
\item[Step 4.] $\beta_2=\lambda_{k_n-1,k_n}^{-1}\gamma w_1\dots w_{n-2}\lambda_{k_n-1,n}\rho_{n-1}B_{k_n,n-1}.$
\item[Step 5.] The braid $\varrho(\beta)$ has the following form.
\begin{align}
\notag\varrho(\beta)&=\lambda_{k_n-1,k_n}^{-1}\gamma (w_1\dots w_{n-2}\lambda_{k_n-1,n})^{\rho_{n-1}}B_{k_n,n-1}\\
\notag&=\lambda_{k_n-1,k_n}^{-1}\gamma (w_1\dots w_{n-2})^{\rho_{n-1}}\lambda_{k_n-1,n-1}B_{k_n,n-1}\\
\notag&=\lambda_{k_n-1,k_n}^{-1}\gamma (w_1\dots w_{n-2})^{\rho_{n-1}}B_{k_n,n-1}\lambda_{k_n-1,n-1}^{B_{k_n,n-1}}\\
\notag&=\lambda_{k_n-1,k_n}^{-1}\gamma (w_1\dots w_{n-2})^{\rho_{n-1}}B_{k_n,n-1}\lambda_{k_n-1,k_n}=\varrho(\alpha)^{\lambda_{k_n-1,k_n}}.
\end{align}
\end{itemize}
Therefore the braids $\varrho(\alpha)$ and $\varrho(\beta)$ are conjugated and by the induction hypothesis the braids $\varrho^*(\alpha)$ and $\varrho^*(\beta)$ are conjugated by the element from $S_n$.

\underline{Case 3.2.3.} $i=k_s$.

\underline{Case 3.2.3.1.} $k_s\leq s-2$. Let us count the braid $\varrho(\beta)$.
\begin{itemize}
\item[Step 1.] The maximal number which is not fixed by $\alpha$ is equal to $s$.
\item[Step 2.] We have $\beta_1=B_{1,n}^{s-n}\beta B_{1,n}^{n-s}$ and $\iota(\beta_1)(n)=k_s+n-s=k_n\neq n$.
\item[Step 3.] Similarly to the case 3.2.1 we have $\beta_1=\lambda_{k_n,k_n+1}^{-1}\alpha_1\lambda_{k_n,k_n+1}$ and therefore
\begin{align}
\notag\beta_1&=\lambda_{k_n,k_n+1}^{-1}\gamma w_1\dots w_{n-1}B_{k_n,n}\lambda_{k_n,k_n+1}\\
\notag&=\lambda_{k_n,k_n+1}^{-1}\gamma w_1\dots w_{n-1}\lambda_{k_n,k_n+1}^{B_{k_n,n}^{-1}}B_{k_n,n}\\
\notag&=\lambda_{k_n,k_n+1}^{-1}\gamma w_1\dots w_{n-1}\lambda_{n,k_n}B_{k_n,n}\\
\notag&=\lambda_{k_n,k_n+1}^{-1}\gamma w_1\dots w_{n-2}\lambda_{n,k_n}w_{n-1}B_{k_n,n}
\end{align}
Since $k_s<s-1$, then $k_n,k_n+1<n$ and the braid $\lambda_{k_n,k_n+1}^{-1}\gamma$ does not involve the strand $n$.
\item[Step 4.] $\beta_2=\lambda_{k_n,k_n+1}^{-1}\gamma w_1\dots w_{n-2}\lambda_{n,k_n}\rho_{n-1}B_{k_n,n-1}.$
\item[Step 5.] The braid $\varrho(\beta)$ has the following form.
\begin{align}
\notag\varrho(\beta)&=\lambda_{k_n,k_n+1}^{-1}\gamma (w_1\dots w_{n-2}\lambda_{n,k_n})^{\rho_{n-1}}B_{k_n,n-1}\\
\notag&=\lambda_{k_n,k_n+1}^{-1}\gamma (w_1\dots w_{n-2})^{\rho_{n-1}}\lambda_{n-1,k_n}B_{k_n,n-1}\\
\notag&=\lambda_{k_n,k_n+1}^{-1}\gamma (w_1\dots w_{n-2})^{\rho_{n-1}}B_{k_n,n-1}\lambda_{n-1,k_n}^{B_{k_n,n-1}}\\
\notag&=\lambda_{k_n,k_n+1}^{-1}\gamma (w_1\dots w_{n-2})^{\rho_{n-1}}B_{k_n,n-1}\lambda_{k_n,k_n+1}=\varrho(\alpha)^{\lambda_{k_n,k_n+1}}
\end{align}
\end{itemize}
Therefore the braids $\varrho(\alpha)$ and $\varrho(\beta)$ are conjugated and by the induction hypothesis the braids $\varrho^*(\alpha)$ and $\varrho^*(\beta)$ are conjugated by the element from $S_n$.

\underline{Case 3.2.3.2} $k_s=s-1.$

\underline{Case 3.2.3.2.1.} $\iota(\alpha)(s-1)=s$. Let us find the braid $\varrho(\beta)=\varrho(\lambda_{k_s,k_s+1}^{-1}\alpha\lambda_{k_s,k_s+1})$.
\begin{itemize}
\item[Step 1.] The maximal number which is not fixed by $\alpha$ is equal to $s$.
\item[Step 2.] We have $\beta_1=B_{1,n}^{s-n}\beta B_{1,n}^{n-s}$ and $\iota(\beta_1)(n)=k_s+n-s=k_n\neq n$.
\item[Step 3.] Similarly to the case 3.2.1 we have $\beta_1=\lambda_{k_n,k_n+1}^{-1}\alpha_1\lambda_{k_n,k_n+1}$. Since $k_s=s-1$, then $k_n=n-1$, $B_{k_n,n}=\rho_{n-1}$ and   $\beta_1=\lambda_{n-1,n}^{-1}\alpha_1\lambda_{n-1,n}$.  Then we have
$$
\beta_1=\lambda_{n-1,n}^{-1}\gamma w_1\dots w_{n-1}\rho_{n-1}\lambda_{n-1,n}.
$$
Since $\iota(\alpha)(s-1)=s$, then $\iota(\alpha_1)(n-1)=n$, therefore the braid $\gamma=\alpha_1\left(w_1\dots w_{n-1}\rho_{n-1}\right)^{-1}$ maps $n-1$ to $n-1$. Hence the braids $\gamma$ and $\lambda_{n-1,n}^{-1}$ commute and we have
\begin{align}
\notag\beta_1&=\gamma\lambda_{n-1,n}^{-1} w_1\dots w_{n-1}\rho_{n-1}\lambda_{n-1,n}\\
\notag&=\gamma w_1\dots w_{n-2}\lambda_{n-1,n}^{-1}w_{n-1}\rho_{n-1}\lambda_{n-1,n}\\
\notag&=\gamma w_1\dots w_{n-2}\lambda_{n-1,n}^{-1}w_{n-1}\lambda_{n,n-1}\rho_{n-1}
\end{align}
where the braid $\gamma$ does not involve the strand $n$, the braid $w_r$ belongs to $\langle\lambda_{n,r},\lambda_{r,n}\rangle$ for $r=1,\dots,n-2$ and $\lambda_{n-1,n}^{-1}w_{n-1}\lambda_{n,n-1}$ belongs to $\langle\lambda_{n,n-1},\lambda_{n-1,n}\rangle$.
\item[Step 4.] $\beta_2=\gamma w_1\dots w_{n-2}\rho_{n-1}.$
\item[Step 5.] $\varrho(\alpha)=\varrho(\beta)$.
\end{itemize}
Therefore the braids $\varrho(\alpha)$ and $\varrho(\beta)$ are conjugated and by the induction hypothesis the braids $\varrho^*(\alpha)$ and $\varrho^*(\beta)$ are conjugated by the element from $S_n$.

\underline{Case 3.2.3.2.2.} $\iota(\alpha)(s-1)=k_{s-1}\leq s-2$. As already founded in the case 3.1.3.2.2. for $k_s=s-1$ and $\iota(\alpha)(s-1)=k_{s-1}\leq s-2$ we have
$$\alpha_1=\eta v_1\dots v_{n-2}B_{k_{n-1},n-1} w_1\dots w_{n-1}\rho_{n-1},$$
$$\varrho(\alpha)=\eta v_1\dots v_{n-2}B_{k_{n-1},n-1} (w_1\dots w_{n-2})^{\rho_{n-1}},$$
$$\varrho^2(\alpha)=\eta (w_1\dots w_{k_{n-1}-1}w_{k_{n-1}+1}\dots w_{n-2})^{\rho_{n-1}B_{k_{n-1},n-1}^{-1}} (v_1\dots v_{n-3})^{\rho_{n-2}} B_{k_{n-1},n-2},
$$
where the braid $\eta$ does not involve the strands $n-1,n$ the braid $v_r$ belongs to $\langle\lambda_{r,n-1},\lambda_{n-1,r}\rangle$ for $r=1,\dots,n-2$ and the braid $w_r$ belongs to $\langle\lambda_{r,n},\lambda_{n,r}\rangle$ for $r=1,\dots,n-1$.
Let us count the braid $\varrho(\beta)=\varrho(\lambda_{s-1,s}\alpha\lambda_{s-1,s})$.
\begin{itemize}
\item[Step 1.] The maximal number which is not fixed by $\beta$ is equal to $s$.
\item[Step 2.] We have $\beta_1=B_{1,n}^{s-n}\beta B_{1,n}^{n-s}$ and $\iota(\beta_1)(n)=k_s+n-s=k_n\neq n$.
\item[Step 3.] Similarly to the case 3.2.3.2.1. $\beta_1=\lambda_{n-1,n}^{-1}\alpha_1\lambda_{n-1,n}$ and we have
\begin{align}
\notag\beta_1&=\lambda_{n-1,n}^{-1}\eta v_1\dots v_{n-2}B_{k_{n-1},n-1} w_1\dots w_{n-1}\rho_{n-1}\lambda_{n-1,n}\\
\notag&=\eta^{\lambda_{n-1,n}}\lambda_{n-1,n}^{-1} v_1\dots v_{n-2}B_{k_{n-1},n-1} w_1\dots w_{n-1}\lambda_{n,n-1}\rho_{n-1}\\
\notag&=\eta^{\lambda_{n-1,n}} v_1\dots v_{n-2}B_{k_{n-1},n-1}(\lambda_{n-1,n}^{-1})^{B_{k_{n-1},n-1}} w_1\dots w_{n-1}\lambda_{n,n-1}\rho_{n-1}\\
\notag&=\eta^{\lambda_{n-1,n}} v_1\dots v_{n-2}B_{k_{n-1},n-1}\lambda_{k_{n-1},n}^{-1} w_1\dots w_{n-1}\lambda_{n,n-1}\rho_{n-1}
\end{align}
Since the braid $\eta$ does not involve the strands $n-1,n$, then $\eta^{\lambda_{n-1,n}}=\eta$.
\item[Step 4.] $\beta_2=\eta v_1\dots v_{n-2}B_{k_{n-1},n-1}\lambda_{k_{n-1},n}^{-1} w_1\dots w_{n-2}\rho_{n-1}.$
\item[Step 5.] $\varrho(\beta)=\eta v_1\dots v_{n-2}B_{k_{n-1},n-1}(\lambda_{k_{n-1},n}^{-1} w_1\dots w_{n-2})^{\rho_{n-1}}$.
\end{itemize}
Let us count the braid $\varrho^2(\beta)$.
\begin{itemize}
\item[Step 1.] Since $k_{s-1}\leq s-2$, then $k_{n-1}\leq n-2$, therefore $B_{k_{n-1},n-1}\neq 1$ and the maximal number which is not fixed by $\varrho(\beta)$ is equal to $n-1$.
\item[Step 2.] $\varrho(\beta)_1=\varrho(\beta)$.
\item[Step 3.] The braid $\varrho(\beta)$ can be rewritten
\begin{align}
\notag\varrho(\beta)&=\eta v_1\dots v_{n-2}B_{k_{n-1},n-1}(\lambda_{k_{n-1},n}^{-1} w_1\dots w_{n-2})^{\rho_{n-1}}\\
\notag&=\eta v_1\dots v_{n-2}(\lambda_{k_{n-1},n}^{-1} w_1\dots w_{n-2})^{\rho_{n-1}B_{k_{n-1},n-1}^{-1}}B_{k_{n-1},n-1}\\
\notag&=\eta (w_1\dots w_{k_{n-1}-1}w_{k_{n-1}+1}\dots w_{n-2})^{\rho_{n-1}B_{k_{n-1},n-1}^{-1}} v_1\dots v_{n-3}\\
\notag&~\cdot~v_{n-2}\lambda_{{n-1},n-2}^{-1}w_{k_{n-1}}^{\rho_{n-1}B_{k_{n-1},n-1}^{-1}}B_{k_{n-1},n-1}
\end{align}
\item[Step 4.] $\varrho(\beta)_2=\eta (w_1\dots w_{k_{n-1}-1}w_{k_{n-1}+1}\dots w_{n-2})^{\rho_{n-1}B_{k_{n-1},n-1}^{-1}} v_1\dots v_{n-3}
\rho_{n-2}B_{k_{n-1},n-2}$.
\item[Step 5.] $\varrho^2(\beta)=\eta (w_1\dots w_{k_{n-1}-1}w_{k_{n-1}+1}\dots w_{n-2})^{\rho_{n-1}B_{k_{n-1},n-1}^{-1}} (v_1\dots v_{n-3})^{\rho_{n-2}}B_{k_{n-1},n-2}$.
\end{itemize}
Therefore $\varrho^2(\alpha)=\varrho^2(\beta)$ and by the induction hypothesis the braids $\varrho^*(\alpha)$ and $\varrho^*(\beta)$ are conjugated by the element from $S_n$.

\underline{Case 3.2.4.} $k_s+1\leq i \leq s-2$. Let us find the braid $\varrho(\beta)=\varrho(\lambda_{i,i+1}^{-1}\alpha\lambda_{i,i+1})$.
\begin{itemize}
\item[Step 1.] The maximal number which is not fixed by $\beta$ is equal to $s$.
\item[Step 2.] We have $\beta_1=B_{1,n}^{s-n}\beta B_{1,n}^{n-s}$ and $\iota(\beta_1)(n)=k_s+n-s=k_n\neq n$.
\item[Step 3.] By Lemma \ref{lem1} we have $\lambda_{i,i+1}B_{1,n}^{n-s}=B_{1,n}^{n-s}\lambda_{j,j+1}$, where $j=i+n-s$. Since $k_s+1\leq i \leq s-2$ then $k_n+1\leq j \leq n-2$ and the braid $\beta$ can be rewritten in more details.
\begin{align}
\notag\beta_1&=B_{1,n}^{s-n}\beta B_{1,n}^{n-s}=B_{1,n}^{s-n}\lambda_{i,i+1}^{-1}\alpha\lambda_{i,i+1} B_{1,n}^{n-s}=\lambda_{j,j+1}^{-1}\alpha_1\lambda_{j,j+1}\\
\notag&=\lambda_{j,j+1}^{-1}\gamma w_1\dots w_{n-1}B_{k_n,n}\lambda_{j,j+1}\\
\notag&=\lambda_{j,j+1}^{-1}\gamma w_1\dots w_{n-1}\lambda_{j,j+1}^{B_{k_n,n}^{-1}}B_{k_n,n}\\
\notag&=\lambda_{j,j+1}^{-1}\gamma w_1\dots w_{n-1}\lambda_{j-1,j}B_{k_n,n}\\
\notag&=\lambda_{j,j+1}^{-1}\gamma\lambda_{j-1,j} w_1\dots w_{n-1}B_{k_n,n},
\end{align}
Where the braid $\lambda_{j,j+1}^{-1}\gamma\lambda_{j-1,j}$ does not involve the strand $n$.
\item[Step 4.] $\beta_2=\lambda_{j,j+1}^{-1}\gamma\lambda_{j-1,j} w_1\dots w_{n-2}\rho_{n-1}B_{k_n,n-1}$.
\item[Step 5.] The braid $\varrho(\beta)$ follows
\begin{align}
\notag\varrho(\beta)&=\lambda_{j,j+1}^{-1}\gamma\lambda_{j-1,j} (w_1\dots w_{n-2})^{\rho_{n-1}}B_{k_n,n-1}\\
\notag&=\lambda_{j,j+1}^{-1}\gamma (w_1\dots w_{n-2})^{\rho_{n-1}}\lambda_{j-1,j}B_{k_n,n-1}\\
\notag&=\lambda_{j,j+1}^{-1}\gamma (w_1\dots w_{n-2})^{\rho_{n-1}}B_{k_n,n-1}\lambda_{j-1,j}^{B_{k_n,n-1}}\\
\notag&=\lambda_{j,j+1}^{-1}\gamma (w_1\dots w_{n-2})^{\rho_{n-1}}B_{k_n,n-1}\lambda_{j,j+1}=\varrho(\alpha)^{\lambda_{j,j+1}}
\end{align}
\end{itemize}
Therefore the braids $\varrho(\alpha)$ and $\varrho(\beta)$ are conjugated and by the induction hypothesis the braids $\varrho^*(\alpha)$ and $\varrho^*(\beta)$ are conjugated by the element from $S_n$.

\underline{Case 3.2.5.} $i=s-1$. In this case we can consider that $k_s< s-1$ since the case when $k_s=s-1=i$ is already solved in the case 3.2.3.2.

\underline{Case 3.2.5.1.}  $\iota(\alpha)(s-1)=s$. From the case 3.1.5.1. we have
$$\alpha_1=\eta v_1\dots v_{n-2} w_1\dots w_{n-1}\rho_{n-1}B_{k_{n},n-1},$$
$$\varrho(\alpha)=\eta v_1\dots v_{n-2} (w_1\dots w_{n-2})^{\rho_{n-1}}B_{k_{n},n-1},$$
$$\varrho^{2}({\alpha})=\eta (v_1\dots v_{n-3})^{\rho_{n-2}} (w_1\dots w_{n-3})^{\rho_{n-1}\rho_{n-2}}B_{k_{n},n-2},$$
where the braid $\eta$ does not involve the strands $n-1,n$.
Let us count $\varrho(\beta)$.
\begin{itemize}
\item[Step 1.] The maximal number which is not fixed by $\beta$ is equal to $s$.
\item[Step 2.] We have $\beta_1=B_{1,n}^{s-n}\beta B_{1,n}^{n-s}$ and $\iota(\beta_1)(n)=k_s+n-s=k_n\neq n$.
\item[Step 3.] By Lemma \ref{lem1} we have $\lambda_{s-1,s}B_{1,n}^{n-s}=B_{1,n}^{n-s}\lambda_{n-1,n}$, therefore
\begin{align}
\notag\beta_1&=B_{1,n}^{s-n}\beta B_{1,n}^{n-s}=B_{1,n}^{s-n}\lambda_{s-1,s}^{-1}\alpha\lambda_{s-1,s} B_{1,n}^{n-s}=\lambda_{n-1,n}^{-1}\alpha_1\lambda_{n-1,n}\\
\notag&=\lambda_{n-1,n}^{-1}\eta v_1\dots v_{n-2} w_1\dots w_{n-1}B_{k_{n},n}\lambda_{n-1,n}\\
\notag&=\eta v_1\dots v_{n-2} \lambda_{n-1,n}^{-1}w_1\dots w_{n-1}\lambda_{n-1,n}^{B_{k_{n},n}^{-1}}B_{k_{n},n}\\
\notag&=\eta v_1\dots v_{n-2} \lambda_{n-1,n}^{-1}w_1\dots w_{n-1}\lambda_{n-2,n-1}B_{k_{n},n}\\
\notag&=\eta v_1\dots v_{n-2}\lambda_{n-2,n-1} w_1\dots w_{n-2} \lambda_{n-1,n}^{-1}w_{n-1}B_{k_{n},n},
\end{align}
where the braid $\eta v_1\dots v_{n-2}\lambda_{n-2,n-1}$ does not involve the strand $n$.
\item[Step 4.] $\beta_2=\eta v_1\dots v_{n-2}\lambda_{n-2,n-1} w_1\dots w_{n-2} \rho_{n-1}B_{k_{n},n-1}$.
\item[Step 5.] $\varrho(\beta)=\eta v_1\dots v_{n-2}\lambda_{n-2,n-1} (w_1\dots w_{n-2})^{\rho_{n-1}}B_{k_{n},n-1}$.
\end{itemize}
Let us count $\varrho^2(\beta)$.
\begin{itemize}
\item[Step 1.] Since the braid $\eta$ does not involve the strand $n-1$, then the image of $n-1$ under the permutation $\iota(\varrho(\beta))$ is equal to $\iota(B_{k_n,n-1})(n-1)=k_n\neq n-1$.
\item[Step 2.]  $\varrho(\beta)_1=\varrho(\beta)$.
\item[Step 3.]The braid $\varrho(\alpha)_1$ has the following form.
\begin{align}
\notag\varrho(\beta)_1&=\varrho(\beta)=\eta v_1\dots v_{n-2}\lambda_{n-2,n-1} (w_1\dots w_{n-2})^{\rho_{n-1}}B_{k_{n},n-1}\\
\notag&=\eta v_1\dots v_{n-3} (w_1\dots w_{n-3})^{\rho_{n-1}}v_{n-2}\lambda_{n-2,n-1}w_{n-2}^{\rho_{n-1}}B_{k_{n},n-1}
\end{align}
\item[Step 4.] $\varrho(\beta)_2=\eta v_1\dots v_{n-3} (w_1\dots w_{n-3})^{\rho_{n-1}}{\rho_{n-2}}B_{k_{n},n-2}$
\item[Step 5.] $\varrho^2(\beta)=\eta v_1\dots v_{n-3} (w_1\dots w_{n-3})^{\rho_{n-1}\rho_{n-2}}B_{k_{n},n-2}=\varrho^2(\alpha)$.
\end{itemize}
By the induction hypothesis the braids $\varrho^*(\alpha)$ and $\varrho^*(\beta)$ are conjugated by the element from $S_n$.

\underline{Case 3.2.5.2.}  $\iota(\alpha)(s-1)=s-1$. In this case we have
\begin{align}
\notag\alpha_1&=\eta v_1\dots v_{n-2}B_{k_{n-1},n-1} w_1\dots w_{n-1}\rho_{n-1}B_{k_{n},n-1}\\
\notag\varrho(\alpha)&=\eta v_1\dots v_{n-2} B_{k_{n-1},n-1}(w_1\dots w_{n-2})^{\rho_{n-1}}B_{k_{n},n-1},
\end{align}
where the braid $\eta$ does not involve the strands $n-1,n$ and $k_{n-1}\geq k_n$.
Since $\iota(\alpha)(s-1)=s-1$, then $\iota(\alpha_1)(n-1)=n-1$, therefore $k_{n-1}=n-2$ and $B_{k_{n-1},n-1}=\rho_{n-2}$.
\begin{align}
\notag\alpha_1&=\eta v_1\dots v_{n-2}\rho_{n-2} w_1\dots w_{n-1}\rho_{n-1}B_{k_{n},n-1}\\
\notag\varrho(\alpha)&=\eta v_1\dots v_{n-2} \rho_{n-2}(w_1\dots w_{n-2})^{\rho_{n-1}}B_{k_{n},n-1},
\end{align}
Let us count the braid $\varrho(\beta)$.
\begin{itemize}
\item[Step 1.] The maximal number which is not fixed by $\beta$ is equal to $s$.
\item[Step 2.] We have $\beta_1=B_{1,n}^{s-n}\beta B_{1,n}^{n-s}$ and $\iota(\beta_1)(n)=k_s+n-s=k_n\neq n$.
\item[Step 3.] Similarly to the case 3.2.5.1. we have $\beta_1=\lambda_{n-1,n}^{-1}\alpha_1\lambda_{n-1,n}$ and therefore
\begin{align}
\notag\beta_1&=\lambda_{n-1,n}^{-1}\eta v_1\dots v_{n-2}\rho_{n-2} w_1\dots w_{n-1}B_{k_{n},n}\lambda_{n-1,n}\\
\notag&=\eta v_1\dots v_{n-2}\lambda_{n-1,n}^{-1}\rho_{n-2} w_1\dots w_{n-1}\rho_{n-1}\lambda_{n-1,n}^{B_{k_{n},n}^{-1}}B_{k_{n},n}\\
\notag&=\eta v_1\dots v_{n-2}\rho_{n-2}(\lambda_{n-1,n}^{-1})^{\rho_{n-2}} w_1\dots w_{n-1}\lambda_{n-2,n-1}B_{k_{n},n}\\
\notag&=\eta v_1\dots v_{n-2}\rho_{n-2}\lambda_{n-2,n-1}\lambda_{n-2,n}^{-1} w_1\dots w_{n-1}B_{k_{n},n},
\end{align}
where the braid $=\eta v_1\dots v_{n-2}\rho_{n-2}\lambda_{n-2,n-1}$ does not involve the strand $n$.
\item[Step 4.] $\beta_1=\eta v_1\dots v_{n-2}\rho_{n-2}\lambda_{n-2,n-1}\lambda_{n-2,n}^{-1} w_1\dots w_{n-2}\rho_{n-1}B_{k_{n},n-1}$.
\item[Step 5.] The braid $\varrho(\beta)$ follows.
\begin{align}\notag\varrho(\beta)&=\eta v_1\dots v_{n-2}\rho_{n-2}\lambda_{n-2,n-1}(\lambda_{n-2,n}^{-1} w_1\dots w_{n-2})^{\rho_{n-1}}B_{k_{n},n-1}\\
\notag\varrho(\beta)&=\eta v_1\dots v_{n-2}\rho_{n-2}\underline{\lambda_{n-2,n-1}\lambda_{n-2,n-1}^{-1}} (w_1\dots w_{n-2})^{\rho_{n-1}}B_{k_{n},n-1}=\varrho(\alpha).
\end{align}
\end{itemize}
By the induction hypothesis the braids $\varrho^*(\alpha)$ and $\varrho^*(\beta)$ are conjugated by the element from $S_n$.

\underline{Case 3.2.5.3.}  $\iota(\alpha)(s-1)\leq s-2$.

\underline{Case 3.2.5.3.1.}  $\iota(\alpha)(s-1)\geq \iota(\alpha)(s)+1$. From the case 3.1.5.3.1. we have
\begin{align}
\notag\alpha_1&=\eta v_1\dots v_{n-2}B_{k_{n-1},n-1} w_1\dots w_{n-1}\rho_{n-1}B_{k_{n},n-1}\\
\notag\varrho(\alpha)&=\eta v_1\dots v_{n-2} B_{k_{n-1},n-1}(w_1\dots w_{n-2})^{\rho_{n-1}}B_{k_{n},n-1}\\
 \notag\varrho^2(\alpha)&=\eta(w_1\dots w_{k_{n-1}-1}w_{k_{n-1}+1}\dots
 w_{n-2})^{\rho_{n-1}B_{k_{n-1},n-1}^{-1}}B_{k_n,n-2}\\
 \notag& ~\cdot~(v_1\dots v_{n-4}v_{n-2})^{B_{k_n,n-2}\rho_{n-2}}w_{k_{n-1}}^{\rho_{n-1}B_{k_{n-1},n-1}^{-1}B_{k_n,n-2}\rho_{n-2}}  B_{k_{n-1}+1,n-2}
\end{align}
Where the braid $\eta$ does not involve the strands $n-1,n$ and $k_{n-1}\geq k_n$.
Let us count the braid $\varrho(\beta)=\varrho(\lambda_{n-1,n}^{-1}\beta\lambda_{n-1,n})$.
\begin{itemize}
\item[Step 1.] The maximal number which is not fixed by $\beta$ is equal to $s$.
\item[Step 2.] We have $\beta_1=B_{1,n}^{s-n}\beta B_{1,n}^{n-s}$ and $\iota(\beta_1)(n)=k_s+n-s=k_n\neq n$.
\item[Step 3.] Similarly to the case 3.2.5.1. we have $\beta_1=\lambda_{n-1,n}^{-1}\alpha_1\lambda_{n-1,n}$ and therefore
\begin{align}
\notag\beta_1&=\lambda_{n-1,n}^{-1}\eta v_1\dots v_{n-2}B_{k_{n-1},n-1} w_1\dots w_{n-1}B_{k_{n},n}\lambda_{n-1,n}\\
\notag&=\eta v_1\dots v_{n-2}\lambda_{n-1,n}^{-1}B_{k_{n-1},n-1} w_1\dots w_{n-1}\lambda_{n-1,n}^{B_{k_{n},n}^{-1}}B_{k_{n},n}\\
\notag&=\eta v_1\dots v_{n-2}B_{k_{n-1},n-1}(\lambda_{n-1,n}^{-1})^{B_{k_{n-1},n-1}} w_1\dots w_{n-1}\lambda_{n-2,n-1}B_{k_{n},n}\\
\notag&=\eta v_1\dots v_{n-2}B_{k_{n-1},n-1}\lambda_{k_{n-1},n}^{-1} w_1\dots w_{n-1}\lambda_{n-2,n-1}B_{k_{n},n}\\
\notag&=\eta v_1\dots v_{n-2}B_{k_{n-1},n-1}\lambda_{n-2,n-1}\lambda_{k_{n-1},n}^{-1} w_1\dots w_{n-1}B_{k_{n},n}
\end{align}
where the braid $\eta v_1\dots v_{n-2}B_{k_{n-1},n-1}\lambda_{n-2,n-1}$ does not involve the strand $n$.
\item[Step 4.] $\beta_2=\eta v_1\dots v_{n-2}B_{k_{n-1},n-1}\lambda_{n-2,n-1}\lambda_{k_{n-1},n}^{-1} w_1\dots w_{n-2}\rho_{n-1}B_{k_{n},n-1}$.
\item[Step 5.] $\varrho(\beta)=\eta v_1\dots v_{n-2}B_{k_{n-1},n-1}\lambda_{n-2,n-1}(\lambda_{k_{n-1},n}^{-1} w_1\dots w_{n-2})^{\rho_{n-1}}B_{k_{n},n-1}$.
\end{itemize}
Let us count $\varrho^2(\beta)$.
\begin{itemize}
\item[Step 1] The maximal number which is not fixed by the braid $\varrho(\beta)$ is $n-1$ and $\iota(\varrho(\beta))(n-1)=k_{n-1}+1$.
\item[Step 2] $\varrho(\beta)_1=\varrho(\beta)$.
\item[Step 3] The braid $\varrho(\beta)_1$ can be rewritten
    \begin{align}\notag\varrho(\alpha_1)&=\eta v_1\dots v_{n-2}B_{k_{n-1},n-1}\lambda_{n-2,n-1}(\lambda_{k_{n-1},n}^{-1} w_1\dots w_{n-2})^{\rho_{n-1}}B_{k_{n},n-1}\\
\notag&=\eta v_1\dots v_{n-2}\lambda_{n-2,n-1}^{B_{k_{n-1},n-1}^{-1}} (\lambda_{k_{n-1},n}^{-1}w_1\dots w_{n-2})^{\rho_{n-1}B_{k_{n-1},n-1}^{-1}}B_{k_{n-1},n-1}B_{k_{n},n-1}\\
\notag&=\eta v_1\dots v_{n-2}\lambda_{n-3,n-2} (\lambda_{k_{n-1},n}^{-1}w_1\dots w_{n-2})^{\rho_{n-1}B_{k_{n-1},n-1}^{-1}}B_{k_{n-1},n-1}B_{k_{n},n-1}
\end{align}
Since $k_{n-1}\geq k_n$, then
\begin{equation}\label{a3}
B_{k_{n-1},n-1}B_{k_{n},n-1}=B_{k_n,n-2}B_{k_{n-1}+1,n-1}
\end{equation}
and and therefore
\begin{align}\notag\varrho(\alpha_1)&=\eta v_1\dots v_{n-2}\lambda_{n-3,n-2} (\lambda_{k_{n-1},n}^{-1}w_1\dots w_{n-2})^{\rho_{n-1}B_{k_{n-1},n-1}^{-1}}B_{k_n,n-2}B_{k_{n-1}+1,n-1}\\
\notag&=\eta\lambda_{n-3,n-2}(w_1\dots w_{k_{n-1}-1}w_{k_{n-1}+1}\dots w_{n-2})^{\rho_{n-1}B_{k_{n-1},n-1}^{-1}} \\
\notag &~\cdot~v_1\dots v_{n-2}(\lambda_{k_{n-1},n}^{-1}w_{k_{n-1}})^{\rho_{n-1}B_{k_{n-1},n-1}^{-1}} B_{k_n,n-2}B_{k_{n-1}+1,n-1}\\
\notag&=\eta\lambda_{n-3,n-2}(w_1\dots w_{k_{n-1}-1}w_{k_{n-1}+1}\dots w_{n-2})^{\rho_{n-1}B_{k_{n-1},n-1}^{-1}} B_{k_n,n-2}\\
\notag &~\cdot~(v_1\dots v_{n-2})^{B_{k_n,n-2}}(\lambda_{k_{n-1},n}^{-1}w_{k_{n-1}})^{\rho_{n-1}B_{k_{n-1},n-1}^{-1}B_{k_n,n-2}} B_{k_{n-1}+1,n-1}\\
\notag&=\eta\lambda_{n-3,n-2}(w_1\dots w_{k_{n-1}-1}w_{k_{n-1}+1}\dots w_{n-2})^{\rho_{n-1}B_{k_{n-1},n-1}^{-1}} B_{k_n,n-2}\\
\notag &~\cdot~(v_1\dots v_{n-4}v_{n-2})^{B_{k_n,n-2}}(\lambda_{k_{n-1},n}^{-1}w_{k_{n-1}})^{\rho_{n-1}B_{k_{n-1},n-1}^{-1}B_{k_n,n-2}} v_{n-3}^{B_{k_n,n-2}} B_{k_{n-1}+1,n-1},
\end{align}
where the braid $\eta\lambda_{n-3,n-2}(w_1\dots w_{k_{n-1}-1}w_{k_{n-1}+1}\dots w_{n-2})^{\rho_{n-1}B_{k_{n-1},n-1}^{-1}} B_{k_n,n-2}$ does not involve the strand $n-1$.
\item[Step 4] The braid $\varrho(\beta)_2$ has the following form
\begin{align}
 \notag\varrho(\beta)_2&=\eta\lambda_{n-3,n-2}(w_1\dots w_{k_{n-1}-1}w_{k_{n-1}+1}\dots w_{n-2})^{\rho_{n-1}B_{k_{n-1},n-1}^{-1}} B_{k_n,n-2}\\
\notag &~\cdot~(v_1\dots v_{n-4}v_{n-2})^{B_{k_n,n-2}}(\lambda_{k_{n-1},n}^{-1}w_{k_{n-1}})^{\rho_{n-1}B_{k_{n-1},n-1}^{-1}B_{k_n,n-2}} \rho_{n-2} B_{k_{n-1}+1,n-2},
\end{align}
\item[Step 5]Finally, the braid $\varrho^2(\alpha)$ follows
\begin{align}
 \notag\varrho^2(\beta)&=\eta\lambda_{n-3,n-2}(w_1\dots w_{k_{n-1}-1}w_{k_{n-1}+1}\dots w_{n-2})^{\rho_{n-1}B_{k_{n-1},n-1}^{-1}} B_{k_n,n-2}\\
\notag &~\cdot~(v_1\dots v_{n-4}v_{n-2})^{B_{k_n,n-2}\rho_{n-2}}(\lambda_{k_{n-1},n}^{-1}w_{k_{n-1}})^{\rho_{n-1}B_{k_{n-1},n-1}^{-1}B_{k_n,n-2}\rho_{n-2}} B_{k_{n-1}+1,n-2}\\
\notag&=\eta B_{k_n,n-2}\lambda_{n-3,n-2}^{B_{k_n,n-2}}(w_1\dots w_{k_{n-1}-1}w_{k_{n-1}+1}\dots w_{n-2})^{\rho_{n-1}B_{k_{n-1},n-1}^{-1}B_{k_n,n-2}} \\
\notag &~\cdot~(v_1\dots v_{n-4}v_{n-2})^{B_{k_n,n-2}\rho_{n-2}}(\lambda_{k_{n-1},n}^{-1}w_{k_{n-1}})^{\rho_{n-1}B_{k_{n-1},n-1}^{-1}B_{k_n,n-2}\rho_{n-2}} B_{k_{n-1}+1,n-2}\\
\notag&=\eta B_{k_n,n-2}\lambda_{n-2,k_n}(w_1\dots w_{k_{n-1}-1}w_{k_{n-1}+1}\dots w_{n-2})^{\rho_{n-1}B_{k_{n-1},n-1}^{-1}B_{k_n,n-2}} \\
\notag &~\cdot~(v_1\dots v_{n-4}v_{n-2})^{B_{k_n,n-2}\rho_{n-2}}\lambda_{n-2,k_n}^{-1}w_{k_{n-1}}^{\rho_{n-1}B_{k_{n-1},n-1}^{-1}B_{k_n,n-2}\rho_{n-2}} B_{k_{n-1}+1,n-2}\\
\notag&=\eta B_{k_n,n-2}(w_1\dots w_{k_{n-1}-1}w_{k_{n-1}+1}\dots w_{n-2})^{\rho_{n-1}B_{k_{n-1},n-1}^{-1}B_{k_n,n-2}} \\
\notag &~\cdot~(v_1\dots v_{n-4}v_{n-2})^{B_{k_n,n-2}\rho_{n-2}}\lambda_{n-2,k_n}\\
\notag&~\cdot~[\lambda_{n-2,k_n},(w_1\dots w_{k_{n-1}-1}w_{k_{n-1}+1}\dots w_{n-2})^{\rho_{n-1}B_{k_{n-1},n-1}^{-1}B_{k_n,n-2}} \\
\notag &~\cdot~(v_1\dots v_{n-4}v_{n-2})^{B_{k_n,n-2}\rho_{n-2}}]
\lambda_{n-2,k_n}^{-1}w_{k_{n-1}}^{\rho_{n-1}B_{k_{n-1},n-1}^{-1}B_{k_n,n-2}\rho_{n-2}} B_{k_{n-1}+1,n-2}
\end{align}
\begin{align}
\notag~~~~~~~&=\eta B_{k_n,n-2}(w_1\dots w_{k_{n-1}-1}w_{k_{n-1}+1}\dots w_{n-2})^{\rho_{n-1}B_{k_{n-1},n-1}^{-1}B_{k_n,n-2}} \\
\notag &~\cdot~(v_1\dots v_{n-4}v_{n-2})^{B_{k_n,n-2}\rho_{n-2}}\\
\notag&~\cdot~[\lambda_{n-2,k_n},(w_1\dots w_{k_{n-1}-1}w_{k_{n-1}+1}\dots w_{n-2})^{\rho_{n-1}B_{k_{n-1},n-1}^{-1}B_{k_n,n-2}} \\
\notag &~\cdot~(v_1\dots v_{n-4}v_{n-2})^{B_{k_n,n-2}\rho_{n-2}}]^{\lambda_{n-2,k_n}^{-1}}w_{k_{n-1}}^{\rho_{n-1}B_{k_{n-1},n-1}^{-1}B_{k_n,n-2}\rho_{n-2}} B_{k_{n-1}+1,n-2}
\end{align}
\end{itemize}
Therefore the braids $\varrho^2(\alpha)$ and $\varrho^2(\beta)$ are equal modulo $UVP_{n-2}^{\prime}$ and by the remark \ref{rmk2} the braids $\varrho^*(\alpha)$ and $\varrho^*(\beta)$ are equal modulo $UVP_n^{\prime}$, therefore $\overline{\varrho^*(\alpha)}=\overline{\varrho^*(\beta)}$.

\underline{Case 3.2.5.3.2.}  $\iota(\alpha)(s-1)\leq \iota(\alpha)(s)$. This case literally repeats the case 3.2.5.3.1 using the fact that in this case $k_{n-1}<k_n$ and then using the equalities
$$B_{k_{n-1},n-1}B_{k_{n},n-1}=B_{k_n-1,n-2}B_{k_{n-1},n-1}$$
instead of equality (\ref{a3}).

\underline{Case 3.2.6.}  $i=s$. By Lemma \ref{lem1} we have $\lambda_{s,s+1}B_{1,n}^{n-s}=B_{1,n}^{n-s}\lambda_{n,1}$, therefore
$\beta_1=\lambda_{n,1}^{-1}\alpha_1\lambda_{n,1}.$ By the definition $\lambda_{1,n}=B_{2,n}\lambda_{1,2}B_{2,n}^{-1}$, therefore
$$\beta_1=B_{2,n}\lambda_{1,2}^{-1}B_{2,n}^{-1}\alpha_1B_{2,n}\lambda_{1,2}B_{2,n}^{-1}=B_{2,n}\lambda_{1,2}^{-1}\delta_1\lambda_{1,2}B_{2,n}^{-1},$$
where $\delta_1=B_{2,n}^{-1}\alpha_1B_{2,n}$. By the case 3.1 the braids $\overline{\varrho^*(\alpha_1)}$ and $\overline{\varrho^*(\delta_1)}$ are conjugated by the element from $S_n$.
Since $\alpha$ is not a pure braid, then $\alpha_1$ and $\delta_1$ are not pure braids. Therefore the maximal number which is not fixed by $\delta_1$ is greater then or equal to $2$. Therefore, by the cases 3.2.1-3.2.5 the braids $\overline{\varrho^*(\delta_1)}$ and $\overline{\varrho^*(\lambda_{1,2}^{-1}\delta_1\lambda_{1,2})}$ are conjugated by the element from $S_n$. Also by the case 3.1 the braids $\overline{\varrho^*(\lambda_{1,2}^{-1}\delta_1\lambda_{1,2})}$ and $\overline{\varrho^*(B_{2,n}\lambda_{1,2}^{-1}\delta_1\lambda_{1,2}B_{2,n}^{-1})}=\overline{\varrho^*(\beta_1)}$ are conjugated by the element from $S_n$, therefore the braids $\overline{\varrho^*(\alpha_1)}=\overline{\varrho^*(\alpha)}$ and $\overline{\varrho^*(\beta_1)}=\overline{\varrho^*(\beta)}$ are conjugated by the element from $S_n$.

\underline{Case 3.2.7.}  $i\geq s+1$. In this case $\beta_1=\lambda_{j,j+1}^{-1}\alpha_1\lambda_{j,j+1}$, where
$$j=\iota(B_{1,n}^{n-s})(i)=i+n-s~(\text{mod}~n)<n.$$
Since the maximal number which is not fixed by $\alpha_1$ is equal to $n$, then $\overline{\varrho^*(\alpha_1)}$ and $\overline{\varrho^*(\beta_1)}$ are conjugated by the element from $S_n$ by the cases 3.2.1-3.2.5.
The proposition is proved.

The following main result of the paper follows immediately from Corollary \ref{cor1} and Proposition \ref{pr1}.
\begin{ttt} Let $\alpha$ and $\beta$ be unrestricted virtual braids. Then their closures $\widehat{\alpha}$ and $\widehat{\beta}$ are equivalent as fused links if and only if $\overline{\varrho^*(\alpha)}$ and $\overline{\varrho^*(\beta)}$ are conjugated by the element from $S_n<UVB_n$.
\end{ttt}

\end{document}